\RequirePackage{fix-cm}
\documentclass[smallcondensed]{svjour3}     
\smartqed  

\usepackage{graphicx}
\usepackage{amssymb}
\usepackage{amsmath}

\usepackage[normalem]{ulem}
\usepackage{amsthm}
\usepackage{enumitem}
\usepackage{latexsym}
\usepackage[english]{babel}
\usepackage[dvipsnames]{xcolor}
\newcommand{\uu}{{\bf u}}

\theoremstyle{remark} 

\newcommand\rhs{\mathrm{rhs}}
\newcommand\diag{\mathrm{diag}}
\newcommand\CFL{\mathrm{CFL}}

\newcommand\EC{\mathrm{EC}}
\newcommand\ES{\mathrm{ES}}
\newcommand\ED{\mathrm{ED}}
\newcommand\TA{\mathrm{TA}}
\newcommand\CN{\mathrm{CN}}

\newcommand{\avg}[1]{\left\{\hspace*{-3pt}\left\{#1\right\}\hspace*{-3pt}\right\}}

\newcommand{\uvec}{\underline{u}}

\newcommand{\fvec}{\underline{f}}

\newcommand{\umat}{\underline{\underline{U}}}

\newcommand{\dmat}{\underline{\underline{D}}}
\newcommand{\mmat}{\underline{\underline{M}}}

\newcommand{\bmat}{\underline{\underline{B}}}

\newcommand{\evec}{{\underline{e}}}

\newcommand\figinfo[1]{[\textit{#1}]}

\newcommand{\textcy}[1]{#1} 

\AtBeginDocument{%
  \setlength{\oddsidemargin}{\dimexpr(\paperwidth-\textwidth)/2-1in}%
  \setlength{\evensidemargin}{\oddsidemargin}%
  \setlength{\topmargin}{%
    \dimexpr(\paperheight-\textheight)/2-\headheight-\headsep-1in}%
}

\begin{document}

\title{Stability issues of entropy-stable and/or split-form high-order schemes\thanks{Gregor Gassner was supported by the European Research Council (ERC) under the European Union's Eights Framework Program Horizon 2020 with the research project Extreme, ERC Grant Agreement No. 714487.}}
\subtitle{\textcy{Analysis of linear stability}}


\author{Gregor J. Gassner \and Magnus Sv{\"a}rd \and  Florian J. Hindenlang }


\institute{G. J. Gassner \at
              Dep. for Mathematics/Computer Science; Center for Data and Simulation Science, University of Cologne, Cologne, Germany\\
              \email{ggassner@uni-koeln.de}       
           \and
           F. J. Hindenlang  \at
              Max-Planck Institute for Plasma Physics, NMPP division; Garching, Germany \\
              \email{florian.hindenlang@ipp.mpg.de}
           \and
           M. Sv{\"a}rd \at
              Dept. of Mathematics, University of Bergen, P.O. Box 7803, 5020 Bergen, Norway\\
            \email{Magnus.Svard@uib.no}
 }             

\date{Received: date / Accepted: date}

\maketitle

\begin{abstract}
The focus of the present research is on the analysis of \emph{\textcy{local energy} stability} of high-order (including split-form) summation-by-parts methods, with e.g. two-point entropy-conserving fluxes, approximating non-linear conservation laws. Our main finding is that \textcy{local energy} stability, \textcy{i.e., the numerical growth rate does not exceed the growth rate of the continuous problem,} is not guaranteed even when the scheme is non-linearly stable and that this \textcy{may have adverse} implications for simulation results. We show that entropy-conserving two-point fluxes are inherently \textcy{locally energy} unstable, as they can be dissipative or anti-dissipative.  Unfortunately, these fluxes are at the core of many commonly used high-order entropy-stable extensions, including split-form summation-by-parts discontinuous Galerkin spectral element methods (or spectral collocation methods).  For the non-linear Burgers equation, we further demonstrate numerically that such schemes cause exponential growth of errors during the simulation. Furthermore, we \textcy{encounter} a similar abnormal behaviour for the compressible Euler equations\textcy{, for a smooth exact solution of a density wave}.  Finally, \textcy{for the same case, we demonstrate numerically that} other commonly \textcy{known} split-forms, such as the Kennedy and Gruber splitting, are also \textcy{locally energy} unstable. 

\keywords{entropy-conserving two-point flux \and entropy-stable high-order scheme \and summation-by-parts \and split-form \and discontinuous Galerkin \and non-linear stability \and linear stability \and \textcy{energy stability}}
\end{abstract}

\section{Introduction}
\label{sec:prologue}
For consistent discretisations of well-posed linear partial differential equations (PDEs), it is well-known that stability (typically in $L^2$) is necessary and sufficient for a numerical solution to converge to the analytical solution. This was proven by Lax and Richtmyer in \cite{LaxRichtmyer56}.  
However, many important PDEs are not linear and the Lax-Richtmyer convergence theorem is not valid. Examples include the Burgers equation and the Euler equations of compressible gas dynamics, which are of primary interest in the present work. It is equally well-known that the equations of fluid dynamics are not known to be well-posed and it is therefore not known what stability properties they have. Nevertheless, there is not a complete lack of knowledge. For smooth solutions of non-linear PDEs, stability of the (well-posed) linearised problem is sufficient to ensure convergence of the numerical solution with very fine grids as shown in \cite{Strang64} by Strang. That is, Strang's theorem reverts the problem back to Lax-Richtmyer's convergence theorem to prove linear well-posedness and stability for subsonic flows, where high-order accurate schemes have been proven to be very efficient, see e.g. \cite{SvardLundberg10}. \textcy{(The necessity of stability for convergence will be discussed in Section~\ref{sec:gen_disc}.)}

\textcy{Next, we turn to solutions in the non-linear regime.}
For transonic or supersonic flows, shocks may appear. To capture shocks, an appropriate amount of artificial diffusion must be added to the scheme. This amounts to ensuring non-linear stability rather than linear. Nevertheless, once the scheme is stabilised, in the sense that the simulation does not blow up, it is again often assumed that the scheme is now convergent. However, for non-linear equations one bound (e.g. in $L^2$) on the solution $u$ is typically not enough to infer (even weak) convergence for the non-linear fluxes $f(u)$. For instance, in case of the Burgers equation (see equation (\ref{conlaw}) below), the flux \textcy{is} $f(u)=u^2/2$ and if a sequence of approximate solutions, $u^n\in L^2$, converges weakly to $u$, $(u^n)^2$ does not necessarily converge to $u^2$ (in $L^1$), see e.g. \cite{Evans1990}. In other words, a numerical solution may be perfectly stable and yet fail to converge to a solution. We also remark, that failure to converge may not only appear for solutions with shocks. Turbulence is a highly non-linear phenomena where linear stability may not be sufficient to infer convergence. 

Since for the flow equations there is little mathematical guidance to what object we want our numerical schemes to converge, we turn to physics. A solution of the compressible Euler equations should not only satisfy the equations but also the second law of thermodynamics. That is, the total specific entropy should only increase in a closed system. This translates into the solution satisfying an extra partial differential inequality. Mathematically this extra condition implies some control (stability) of the solution. This has led to the development of so-called entropy-conservative/stable schemes, e.g., \cite{Tadmor84,Tadmor87,FjordholmMishra12,FisherCarpenter13,MadraneFjordholm12,winters2017uniquely,Friedrich2018,Gassner_BR1,HoneinMoin04,IsmailRoe09,RayChandrashekar17,YamaleevCarpenter09,Crean2018,Parsani2016,Parsani:2015:ESD:2784985.2785151,Sandham2002,sjogreen2017skew,Barth1999,murman2016}.

Herein, we  study the \textcy{local energy} stability of entropy-conservative/stable high-order (including split-form) summation-by-parts (SBP) schemes. In the following subsections, we define the concepts of linear stability, entropy-stability and \textcy{local energy} stability. In particular, we differentiate between entropy-stability and entropy-dissipation.

\subsection{Linear stability}\label{sec:linear}

Consider a (system of) conservation law(s) on the periodic (in space) domain $[-1,1]\times (0,T]$.
\begin{equation}
  u_t + f(u)_x =0,\quad u(x,0)=u_0(x),  \label{conlaw}
\end{equation}
where $u$ is the solution vector/scalar and $f(u)$ is the vector/scalar flux function. The problem is closed with the initial datum $u_0$.

We begin by considering stability of linear problems. A linear PDE, $u_t=Pu$, on a periodic domain $\Omega$ with initial data \textcy{$u_0\in L^2([-1,1])$}, is well-posed if there exists a unique solution $\|u(\cdot,T)\|_2=K_c\exp(\alpha_c T)\,\|u_0\|_2$ where $\|u\|_2=(\int_\Omega u^2\,dx)^{1/2}$ is the $L^2$ norm. The constants $K_c$ and $\alpha_c$ are independent of the initial data. (We tacitly assume that the estimate is sharp.) Let $\uu_t=P_h \uu$ be a consistent approximation of the PDE, where $\uu=(u_0,u_1,...,u_N)^T$ is the discrete solution vector and
\begin{align}
\|\uu\|^2=h\uu^T\uu=\sum_{i=0}^{N} hu_i^2\label{disc_l2}
\end{align}
an $L^2$ equivalent norm.  We say that the scheme is \emph{stable}, if $(\|\uu(\cdot,T)\|_2)_t=K_n\exp(\alpha_n T)\|\uu_0\|_2$. We refer to \cite{GustafssonKreissOliger} for more information on linear well-posedness and stability. Clearly, this notion of stability implies convergence in the sense of Lax-Richtmyer. When $K_n$ and $\alpha_n$ are any bounded constants, we term the scheme \emph{Lax stable}. \textcy{In view of Lax-Richtmyer's Equivalence Theorem, convergence implies (Lax) stability. That is, Lax stability is \emph{necessary} for convergence of a consistent scheme.}

\begin{remark}
\textcy{\textit{We remind that the notions \emph{stability} and \emph{convergence} concern classes (e.g. $L^2$) of initial data and solutions. A scheme that is neither stable nor convergent, need neither blow up nor fail to converge for every simulation.}}
\end{remark}

To introduce some more restrictive stability notions, we consider the linear advection equation. That is, we consider (\ref{conlaw}) with $f=a\,u$ where  $a$ is the constant advection coefficient. By multiplying from the left by $u$ and integrating in space and time, using the periodic boundary conditions, we obtain, $\|u(\cdot,T)\|=\|u_0\|$. The same estimate can be derived using Fourier analysis, and it is sharp. As is well-known, all eigenvalues of the symbol are purely imaginary. That is, $\alpha_c=0$. If no low-order term ($u$ without a derivative) is present in the PDE, $\alpha_c = 0$ for any linearly well-posed PDE \cite{GustafssonKreissOliger}. 

Turning to a semi-discretisation of (\ref{conlaw}) with the linear flux function, let $u_i(t)$ be the approximation of $u(x_i,t)$, where $x_i=i\,h$, $i=0,...,N$ and $h>0$ is the grid spacing. 
Then we consider the semi-discrete finite volume scheme \textcy{with a central flux},
\begin{align}
(u_i(t))_t + \frac{ f_{i+1/2} - f_{i-1/2} }{ h }=0,\label{FV}\\
f^{\CN}_{i+1/2}=\frac{f(u_{i+1})+f(u_i)}{2}=\frac{au_{i+1}+au_i}{2}.\label{central_advec}
\end{align}
By inserting (\ref{central_advec}) into (\ref{FV}), multiplying by $hu_i$ and summing $i$ from $0$ to $N$, we obtain $\|\uu(\cdot,T)\|=\|\uu_0\|$, where $\uu=(u_0,u_1,...,u_N)^T$ with the discrete norm \eqref{disc_l2}.

One typically refers to this notion as $L^2$ stability, or energy well-posedness, or energy stability. As in the continuous case, the estimate can be arrived at by Fourier analysis and the eigenvalues of the discrete symbol are all imaginary. Either way, \textcy{for the discretization with the central flux} we have $\alpha_n=\alpha_c=0$. The growth of the numerical scheme \textcy{never exceeds} that of the continuous equation.  

Clearly, energy stability with $\alpha_n=0$ implies Lax stability (and convergence). However, \textcy{Lax stability can also hold for} a scheme with eigenvalues with positive real parts. (Positive eigenvalues, for short.) For instance, we could add an anti-dissipative term with a coefficient that vanishes as $h\rightarrow 0$. Such a scheme would have $\alpha_n\sim h$ and would converge. One could even have $\alpha_n=\text{constant}$. Since we require the scheme to be consistent, those eigenvalues have to be associated with the highest frequency modes, which are small for $L^2$ bounded data. However, Lax stable schemes \textcy{are often} impractical since they \textcy{might} require an excessive resolution to control the exponentially growing errors. \textcy{For that reason, we want energy stability to guarantee that the growth of the numerical scheme can never exceed that of the continuous problem. As energy-stability implies Lax stability, we refer to linear stability if the scheme is at least Lax stable.}

\begin{remark}
\textcy{\textit{A well-posed linear constant-coefficient PDE, with no zeroth-order term, has a spectrum with non-positive eigenvalues. For well-posed linear PDEs with variable coefficients and/or zeroth-order term the spectrum of the spatial differential operator is bounded, but possibly positive, and a Lax stable scheme may exceed the maximal growth of the continuous problem. See \cite{GustafssonKreissOliger} for further information.}}
\end{remark}


With these definitions of Lax stability, energy stability, and linear stability, the central difference scheme (\ref{FV})-(\ref{central_advec}) is the marginally \textcy{energy} stable scheme.

\subsection{Entropy-stability}

While entropy-stability analysis coincides with the linear stability analysis from the previous section in certain special cases (e.g. analysis of scalar linear problems with quadratic entropy), the main motivation for entropy-stability is the investigation of non-linear problems. We define an entropy (function) for a (non-linear) conservation law as a convex function $U(u)$ that also satisfies $U_u\,f_u=F_u$ where $F(u)$ is the corresponding entropy flux.

A weak solution, $u$, of (\ref{conlaw}) is termed \emph{entropy solution}, if it also satisfies,
\begin{align}
U(u)_t+F(u)_x\leq 0, \label{ent_ineq}
\end{align}
which is an incarnation of the second law of thermodynamics.

The idea of entropy-stability is to design a numerical scheme for (\ref{conlaw}) that also satisfies (\ref{ent_ineq}). In the review \cite{Tadmor03}, it is detailed how to achieve this for a finite volume scheme. We summarise the results here. First, we denote the vector $U_u=w^T$ as the \emph{entropy variables}. We consider the finite volume scheme (\ref{FV}) with non-linear flux functions $f(u)$. Multiplying by the entropy variables of cell $i$, $w_i^T$, and assuming that,
\begin{align}
    (w_{i+1}^T-w_{i}^T)\,f_{i+1/2}-(\Psi_{i+1}-\Psi_{i})\leq 0,\label{EC_cond}
\end{align}
one arrives at,
\begin{align}
    (U_i)_t +\frac{F_{i+1/2}-F_{i-1/2}}{h}\leq 0,\quad \forall i, \label{disc_ent_ineq}
\end{align}
where  $\Psi=w^Tf-F$ is the \emph{entropy potential} and $F_{i+1/2}=\frac{1}{2}(w_{i+1}^T+w_{i}^T)f_{i+1/2}-\frac{1}{2} (\Psi_{i+1}+\Psi_{i})$ is the numerical entropy flux. The inequality (\ref{disc_ent_ineq}) is a consistent approximation of (\ref{ent_ineq}). What remains to do is to find fluxes $f_{i+1/2}$ that satisfy (\ref{EC_cond}). 

As there is a somewhat confusing and loose definition of entropy-stability in the recent high-order literature, we introduce the following nomenclature used throughout this paper: A scheme is
\begin{itemize}
    \item[(i)] \emph{entropy-stable}, if (\ref{disc_ent_ineq})/(\ref{EC_cond}) is satisfied at all $i$ and for \emph{any} admissible entropy function $U(u)$; 
    \item[(ii)] \emph{entropy-conservative}, if (\ref{disc_ent_ineq})/(\ref{EC_cond}) is satisfied as an equality for \emph{one} entropy function at all $i$;
    \item[(iii)] \emph{entropy-dissipative}, if it satisfies (\ref{disc_ent_ineq})/(\ref{EC_cond}) but is neither entropy-stable nor en\-tro\-py-conservative. 
\end{itemize}
We will also refer to fluxes as having these properties with respect to (\ref{EC_cond}) (locally).\medbreak

\begin{remark}\textit{We refer to e.g. Tadmor \cite{Tadmor03} for a proof that both the Rusanov or Lax-Friedrichs schemes are entropy-stable discretisations. There are several ways to extend the finite volume entropy-conserving or dissipative scheme to high-order accuracy. A popular approach is to use high-order interpolation/reconstruction based on entropy-conservative fluxes within the blocks/elements, while the fluxes across the grid blocks/element interfaces can be either entropy-conservative or dissipative, see e.g. \cite{LeFlochMercier02,LeFlochRohde00,skew_sbp2,carpenter_esdg,FisherCarpenter13,gassner_skew_burgers,Carpenter_etal16}. Note that using entropy-stable fluxes in the volume terms reduces the accuracy to first or second order.}   

\end{remark}

\subsection{\textcy{Local energy} stability}\label{sec:loc_lin}

In fluid dynamics, e.g. in turbulence research \cite{Mankbadi1994}, it is common to investigate the stability of certain flow states and their transition from a linear laminar behaviour to a non-linear turbulent behaviour. The analysis proceeds by linearising the non-linear system of conservation laws around a given baseflow state, e.g. \cite{Trefethen1993}. The linearised equation reveals if the baseflow state is locally stable in space and time (the flow remains laminar), or allows a growth (local instability) and transitions to turbulence (non-linear behaviour). 

This idea is closely connected to linear stability \textcy{in the sense of Strang \cite{Strang64}}. It is reasonable to demand that the numerical scheme has the same local behaviour as the equations themselves. For linear PDEs, we can immediately state that \textcy{energy} stable schemes have this property (they do not allow nonphysical growth, \textcy{as the growth of the numerical scheme does not exceed the growth of the continuous problem}), while Lax stable schemes might violate this principle. 

Next, we apply the idea used in fluid dynamics research to extend the linear stability analysis to non-linear equations. Let us assume that we have a smooth 'baseflow' $\widetilde{u}(x,t)$ that satisfies the non-linear PDE. In hydrodynamics research, this state is typically a steady state, but here we allow a time dependence. We will also interchangeably denote the baseflow as the linearisation state throughout the manuscript. \textcy{In general,} local linear stability concerns the time evolution of small fluctuations $|u'(x,t)|<<|\widetilde{u}(x,t)|$ added to the baseflow. Thus, we make the ansatz $u(x,t) = \widetilde{u}(x,t) + u'(x,t)$ for the solution of the non-linear problem and derive a linear equation for the fluctuations $u'$. Consider Burgers equation, i.e.,  (\ref{conlaw}) with $f(u)=u^2/2$ on a periodic domain. We insert the ansatz
\begin{equation}\label{eq:linearised_variable}
\begin{split}
&(\widetilde{u}(x,t) + u'(x,t))_t + \frac{1}{2}\left((\widetilde{u}(x,t) + u'(x,t))(\widetilde{u}(x,t) + u'(x,t))\right)_x = 0,
\end{split}
\end{equation}
neglect all but the leading order \textcy{differential} terms \textcy{(zeroth-order terms do not affect linear stability, see \cite{GustafssonKreissOliger})}, and use that $\widetilde{u}$ is an exact solution to arrive at
\begin{equation}
\label{eq:linearised_burgers}
(u')_t + (\widetilde{u}(x,t)\,u')_x = 0,
\end{equation}
which is a variable coefficient advection equation in conservative form with the advection velocity $\widetilde{u}(x,t)$ for the fluctuations $u'(x,t)$. \textcy{As mentioned above, we are interested in the stricter stability property, i.e., the local energy stability} of the non-linear problem, at the linearisation state $\widetilde{u}(x,t)$, which is determined by the energy analysis of (\ref{eq:linearised_burgers}). 
Multiplying by $u'$ and integrating over the domain, we get after some further manipulations
\begin{equation}
\frac{d}{dt}\|u'\|^2(t) = \int\limits_{-1}^1 \left(\widetilde{u}(x,t)\right)_x\,(u'(x,t))^2\,dx.
\end{equation} 
Since $(\widetilde{u})_x$ is a known and integrable function, we obtain 
\begin{equation}\label{standard_est}
\|u'(\cdot,T)\|^2\leq \exp\left[T\,\frac{1}{2}\sup\limits_{x,t}|\left(\widetilde{u}(x,t)\right)_x|\right]\,\|u'(\cdot,0)\|^2.
\end{equation}
This shows that in general, a \underline{bounded} growth \textcy{$\alpha_c =\frac{1}{2}\sup\limits_{x,t}|\left(\widetilde{u}(x,t)\right)_x|>0 $} of the fluctuations is possible. Hence, the variable coefficient advection problem is linearly well-posed. \textcy{At the end of Section~\ref{sec:fd_burgers}, we present a sharper estimate than \eqref{standard_est} for a specific choice of baseflows $\widetilde{u}$, i.e., for a specific class of variable coefficient problems, that enables us to identify and rule out schemes that are only Lax stable.}\medbreak

\subsection{\textcy{General discussion on stability and convergence}}\label{sec:gen_disc}
\textcy{
We end the Introduction with a general discussion on stability and convergence. 
As already mentioned, if a problem is ($L^2$-)well-posed, a convergent scheme should converge for any $L^2$ bounded initial data. Non-convergence does not imply that solutions \emph{always} diverge. The same is true for \emph{stable} vs. \emph{unstable} schemes.}  

\textcy{Turning  to a non-linear problem, Strang's theorem states that linear stability (either the stricter energy stability or the weaker Lax stability) of the first variation is sufficient (and  necessary since a convergent sequence is bounded) for convergence to a smooth solution. The smoothness assumption guarantees that the first variation dominates the higher-order variations for sufficiently small perturbations. However, the first variation is a variable coefficient problem with zeroth-order terms, which may induce a growth. Linear stability guarantees that any growth is bounded, as seen in Section~\ref{sec:linear}.  Indeed, this growth may push the solution into a regime where the non-linear terms are no longer negligible. This is the mechanism that produces turbulence or shocks even when energy stable, such as central difference, schemes are used.  (In the same way, artificial growth that exceeds the growth of the continuous problem caused by an energy unstable operator may erroneously push the solution into the non-linear regime.) Hence, the spectrum of an energy stable difference operator approximating a variable coefficient problem may be positive (but the real part is bounded from above). Since variable coefficients and zeroth-order terms are "benign" with respect to linear stability, they are often ignored and the problem is reduced to the constant coefficient case (known as "freezing" the coefficients.) The constant coefficient case highlights the effect of the difference operators and conveniently makes the continuous spectrum non-positive. An energy stable scheme mimics this property, while a Lax stable scheme may allow growth in the constant coefficient problem. In the cases we study numerically, we do not freeze the coefficients and to be able to compare schemes, we therefore choose continuous problems without growth, i.e., a special case where $\alpha_c = 0$. (The theory underpinning these statements is found in \cite{GustafssonKreissOliger}.)\\}

The remainder of the paper is organised as follows: In Section~\ref{sec:burgers}, we present the local stability analysis of entropy-conservative/dissipative schemes for Burgers equation, including skew-symmetric high-order discontinuous Galerkin spectral element approximations. \textcy{The analysis reveals \textcy{local energy} stability issues with such schemes and that they are Lax stable at most.} In Section~\ref{sec:euler}, we show that the \textcy{local energy} stability issues carry over to the case of the compressible Euler equations in multi-dimensions and provide further numerical evidence.  We draw our conclusions in Section~\ref{sec:conclusion}. 

\section{\textcy{Local energy} stability analysis of the non-linear Burgers equation}
\label{sec:burgers}

The split-form of the non-linear Burgers equation is
\begin{equation}\label{eq:split_alpha}
u_t + \alpha\,\left(u^2/2\right)_x + (1-\alpha)\,u\,u_x= 0,
\end{equation}
where $\alpha\in[0,1]$ is the split parameter. The choice $\alpha=1$ gives the divergence form of the PDE, whereas $\alpha=0$ gives the advective quasi-linear form. It is well-known that the $2/3$ - trick can be used to derive an entropy estimate for the non-linear equation. The Burgers equation with $\alpha=2/3$ takes the form
\begin{equation}
u_t +\frac{1}{3}u\,u_x+\frac{1}{3}(u^2)_x=0.\label{onethird}
\end{equation}
Multiplying by $u$ and integrating in space gives
\begin{equation}
\frac{1}{2}(\|u\|^2)_t +\left[ \frac{1}{3}(u^3)\right]_{-1}^1=0,
\end{equation}
and with periodicity
\begin{align}
\frac{d}{dt}(\|u\|^2)=0.\label{nonlin_L2}
\end{align}
This trick is used to design entropy-conservative fluxes and entropy-conserving split-form schemes.

\subsection{The finite volume perspective}
\label{sec:fv_burgers}

By approximating \eqref{eq:split_alpha}, with central finite differences, we obtain
\begin{align}
(u_i)_t+\alpha\frac{\frac{1}{2}u_{i+1}^2-\frac{1}{2}u_{i-1}^2}{2h}+(1-\alpha)\,u_i\frac{u_{i+1}-u_{i-1}}{2h} = 0.\label{split_burgers}
\end{align}
It is straight forward to check that the approach is conservative, as it can be recast into the finite volume form (\ref{FV}) with the split-flux function
\begin{equation}
\label{eq:alpha_flux}
f^{\alpha}_{i+1/2} =  \frac{1}{2}\left(f(u_i) + f(u_{i+1})\right) - \frac{1}{2}\frac{(1-\alpha)\,(u_{i+1}-u_i)}{2}\,(u_{i+1}-u_i).
\end{equation}
The choice $\alpha = 2/3$ corresponds to the skew-symmetric form \eqref{onethird} and gives the entropy-conserving flux function
\begin{equation}
\label{eq:ec}
f^{\EC}_{i+1/2}  =  \frac{1}{2}\left(f(u_i) + f(u_{i+1})\right) - \frac{1}{2}\left(\frac{u_{i+1}-u_i}{6}\right)(u_{i+1}-u_i)= \frac{1}{6}(u_i^2 + u_i\,u_{i+1} + u_{i+1}^2).
\end{equation}

To gain some further insights into the properties of the fluxes \eqref{eq:alpha_flux} and \eqref{eq:ec}, we interpret them as a standard finite volume flux of the form,
\begin{align}
    f_{i+1/2}=\frac{1}{2}\left(f(u_i) + f(u_{i+1})\right)-\frac{1}{2}R_{i+1/2}(u_{i+1}-u_i),\label{FVflux}
\end{align}
where $R_{i+1/2}$ is a dissipation coefficient/matrix. We hence get for the fluxes \eqref{eq:alpha_flux} and \eqref{eq:ec} the dissipation coefficients 
\begin{equation}
\label{burgers_lambda}
R^{\alpha}_{i+1/2} =\frac{(1-\alpha)\,\left(u_{i+1}-u_i\right)}{2}, \quad R^{\EC}_{i+1/2} =\left(\frac{u_{i+1}-u_i}{6}\right). 
\end{equation}
Note, that the coefficients \eqref{burgers_lambda} are not necessarily definite. Only for the choice $\alpha = 1$, the dissipation coefficient is guaranteed to be zero for all values $u_{i+1},\,u_i$, and the fluxes reduce to the central flux. For general $\alpha$, in particular for the skew-symmetric entropy-conserving choice $\alpha = 2/3$, the flux can be positive and thus dissipative, but also negative and thus anti-dissipative. We note that for the entropy-conserving flux, this was already observed in \cite{Tadmor03}.

By linearising the finite volume scheme \eqref{FV} with the split-fluxes \eqref{eq:alpha_flux} around a smooth baseflow $\widetilde u$ that is assumed to satisfy the scheme, it is straightforward to deduce that we get again a finite volume form \eqref{FV} with the linearised fluxes 
    \begin{equation}
    \widetilde{f}^{\alpha}_{i+1/2} = \frac{(\widetilde{u}_{i+1}\,u'_{i+1})+(\widetilde{u}_i\,u'_i)}{2}-\frac{1}{2}\,\widetilde{R}^{\alpha}_{i+1/2}\,(u'_{i+1} - u'_i),
    \end{equation}
where the dissipation coefficient of the linearised flux is 
\begin{equation}
\widetilde{R}^{\alpha}_{i+1/2} = (1-\alpha)\,(\widetilde{u}_{i+1}-\widetilde{u}_{i}).  
\end{equation}
We note that the linearised flux $\widetilde{f}^{\alpha}_{i+1/2}$ is consistent with the linearised problem \eqref{eq:linearised_variable} for all split parameters $\alpha$. Furthermore, the dissipation coefficient of the linearised fluxes $\widetilde{R}^{\alpha}_{i+1/2}$ are for arbitrary $\alpha\ne 1$ positive or negative, depending on the slope of the baseflow $\widetilde{u}$. Clearly, the scheme is stable for all $\alpha$ in the sense of Lax for smooth $\widetilde u$, but the anti-dissipative parts will cause positive eigenvalues and we thus do not have \textcy{energy} stability (see Section~\ref{sec:linear}) and consequently the scheme is not \textcy{locally energy} stable. 

We emphasise that the \textcy{local energy} analysis simply demonstrates that the anti-dissipative character observed in \eqref{burgers_lambda}, is unchanged in the linear regime. Only for the choice $\alpha=1$ the dissipation coefficients are $R^{\alpha}_{i+1/2}=0$ and $\widetilde{R}^{\alpha}_{i+1/2}=0$ for all $i$, which is the marginal case for \textcy{energy} stability. This is well-known since the scheme turns into the central difference scheme of the divergence form, which has purely imaginary eigenvalues.

\begin{remark}\textit{
We note that for the entropy-conserving case, $\alpha=2/3$, Tadmor \cite{Tadmor03} proposed to remove the anti-dissipative character by replacing the en\-tro\-py-con\-serv\-ing dissipation coefficient $R^{\EC}_{i+1/2}$ in \eqref{burgers_lambda} with}
\begin{equation}
R^{\TA}_{i+1/2}=\max\left(\frac{1}{6}(u_{i+1}-u_i),0\right). \label{burgers_stable2}
\end{equation}
\textit{Then the corresponding numerical flux becomes entropy-dissipative with respect to the quadratic entropy $U(u)=u^2/2$.}
\end{remark}

An alternative way to construct an entropy-dissipative flux is to add an explicit viscosity term to an entropy-conserving flux. For instance, a Rusanov-type viscosity term of the form
\begin{equation}
\label{eq:es}
\begin{split}
f^{\ED}_{i+1/2} &= \frac{1}{6}(u_i^2 + u_i\,u_{i+1} + u_{i+1}^2) - \frac{1}{2}\max(|u_i|,|u_{i+1}|)\,(u_{i+1}-u_i)\\
&=  \frac{1}{2}\left(f(u_i) + f(u_{i+1})\right) - \frac{1}{2}\left(\frac{u_{i+1}-u_i}{6}+\max(|u_i|,|u_{i+1}|)\right)(u_{i+1}-u_i),
\end{split}
\end{equation}
which corresponds to $R^{\ED}_{i+1/2} =\left(\frac{u_{i+1}-u_i}{6}+\max(|u_i|,|u_{i+1}|)\right)\geq 0$. In contrast, the standard Rusanov flux function, see e.g. \cite{Toro:1999yq}, is an entropy-stable flux function
\begin{equation}
f_{i+1/2}^{\ES} = \frac{1}{2}\left(f(u_i) + f(u_{i+1})\right) - \frac{1}{2}\max(|u_{i+1}|,|u_i|)\,(u_{i+1}-u_i),
\end{equation}
which corresponds to $R^{\ES}_{i+1/2} = \max(|u_{i+1}|,|u_i|)$.

\subsection{The finite difference perspective}
\label{sec:fd_burgers}

To give an alternative perspective on the issue, we consider the  split-form \eqref{split_burgers} directly in its finite difference form. One reason of the popularity of the split-form schemes \eqref{split_burgers} is that for the choice $\alpha=2/3$ (skew-symmetric form) the continuous estimate \eqref{nonlin_L2} carries over 
\begin{equation}\label{global_est}
\frac{d}{dt}(\|\uu\|^2)=0.
\end{equation}
Hence, the skew-symmetric scheme is neutrally stable in $L^2$ and (once again) we conclude that it is \emph{entropy-conservative} with respect to the quadratic entropy $U(u)= u^2/2$. This appears to be at odds with the conclusion above that the scheme is not \textcy{locally energy} stable. However, \textcy{local energy} stability refers to the local behaviour in space and time deduced from the linearised problem. The estimate \eqref{global_est} shows that there is an upper limit to the local (un-)physical growth. Below, however, we will demonstrate numerically, that this global bound allows for a significant growth of local fluctuations.\medbreak

To analyse the finite difference approximation \eqref{split_burgers} we first rewrite it in the following equivalent way 
\begin{equation}
(u_i)_t + \frac{\frac{1}{2}u_{i+1}^2 - \frac{1}{2}u_{i-1}^2}{2h} - \frac{(1-\alpha)}{2}\,h^2\,\frac{(u_{i+1}-u_{i-1})}{2\,h}\,\frac{(u_{i+1} - 2\,u_i+u_{i-1})}{h^2} =0. 
\end{equation}
We can reinterpret this form as a consistent discretisation of the following continuous problem 
\begin{equation}
u_t + \left(\frac{u^2}{2}\right)_x - \nu\,u_{xx} = 0,    
\end{equation}
where the last term corresponds to dissipation with the coefficient
\begin{equation}
\nu = \frac{(1-\alpha)}{2}\,h^2\,u_x,
\end{equation}
that analogously to the finite volume perspective gets anti-dissipative for negative slopes $u_x<0$, if $\alpha\ne 1$ \textcy{and $h>0$}. 

As before, we linearise the scheme around a baseflow $\widetilde{u}_i(t)$ (which solves \eqref{split_burgers}) and insert the ansatz $u_i = \widetilde{u}_i+u'_i$ into the scheme (\ref{split_burgers}). Neglecting high-order terms of the small amplitude fluctuations, such as $(u'_i)^2$, we get
\begin{equation}
\label{eq:discrete_linear_burgers}
(u'_i)_t + \alpha\,\frac{\widetilde{u}_{i+1} u'_{i+1}-\widetilde{u}_{i-1}u'_{i-1}}{2h}+(1-\alpha)\left[\widetilde{u}_i\frac{u'_{i+1}-u'_{i-1}}{2h}+u'_i\frac{\widetilde{u}_{i+1}-\widetilde{u}_{i-1}}{2h}\right] = 0,
\end{equation}
which is a consistent approximation of 
\begin{equation}
(u')_t + \alpha\,(\widetilde{u}\,u')_x + (1-\alpha)\,\left(\widetilde{u}\,(u')_x + u'\,\widetilde{u}_x\right) = 0,
\end{equation}
which in turn is a split-form of the continuous linearised variable coefficient advection problem \eqref{eq:linearised_burgers}. Again, it is possible to equivalently rewrite this discretisation into 
\begin{equation}
\begin{split}
(u'_i)_t + \frac{\widetilde{u}_{i+1} u'_{i+1}-\widetilde{u}_{i-1}u'_{i-1}}{2h} 
-&\frac{(1-\alpha)}{2}h^2\,\Bigg[\frac{(\widetilde{u}_{i+1}-\widetilde{u}_{i-1})}{2h}\frac{(u'_{i+1}-2\,u'_i+u'_{i-1})}{h^2}\\
&+\frac{(u'_{i+1} - u'_{i-1})}{2h}\frac{(\widetilde{u}_{i+1}-2\,\widetilde{u}_{i}+\widetilde{u}_{i-1})}{h^2}\Bigg] = 0,
\end{split}
\end{equation}
which we can reinterpret as a consistent discretisation of the following continuous problem
\begin{equation}
\begin{split}
&u_t + (\widetilde{u}\,u')_x - \frac{(1-\alpha)}{2}\,h^2\,(\widetilde{u}_x\,u'_{xx}+u'_x\,\widetilde{u}_{xx}) =0\\ &u_t + (\widetilde{u}\,u')_x - (\widetilde{\nu}\,u'_x)_x =0,\\  
\end{split}    
\end{equation}
where the last term corresponds to dissipation with the coefficient
\begin{equation}
\widetilde{\nu} = \frac{(1-\alpha)}{2}\,h^2\,\widetilde{u}_{x}.   
\end{equation}
Analogously to the finite volume perspective, the coefficient gets anti-dissipative for negative slopes of the baseflow $\widetilde{u}_x<0$, if $\alpha\ne 1$ \textcy{and $h>0$.} \textcy{Note that this may lead to a significant growth for a finite $h$. By taking $h$ sufficiently small, the growth can be reduced but it may require extra resolution.}

\begin{remark}\textit{
The analysis was performed for general $\alpha$ and shows that \underline{any} split-form is at most Lax stable, except for the central approximation with $\alpha=1$. In particular, the skew-symmetric split-form, that is entropy-conserving, is at most Lax stable.} 
\end{remark}

The effect of split-form approximations for variable coefficient problems has also been studied in Manzanero et al. \cite{Manzanero_etal17}. They showed that with the additional assumption $\widetilde{u}(x,t)~>~0,\,\forall x,t$, a sharper estimate than \eqref{standard_est} can be derived
\begin{equation}\label{manzanero}
\min\limits_{x\in[-1,1]}\,\left\{\widetilde{u}(x)\right\}\,\|u'(.,T)\|^2\leq \max\limits_{x\in[-1,1]}\,\left\{\widetilde{u}(x)\right\}\,\|u'_0(.)\|^2.
\end{equation}
The estimate is interesting, as it shows that for the special choice of positive linearisation states there is \underline{no} exponential growth of the solution, i.e. $\alpha_c = 0$ in this particular case. In addition, Manzanero et al. showed that for positive linearisation states the central approximation based on the divergence form ($\alpha=1$) of the PDE with central fluxes is \emph{the only scheme} that has the correct behaviour, i.e. there is no artificial exponential growth or damping. \emph{The central scheme \eqref{FV} with the numerical flux $f^{\CN}_{i+1/2} = \frac{1}{2}(f(u_i) + f(u_{i+1}))$ is thus the least diffusive \textcy{energy} stable scheme for such a problem.} This is exactly the same conclusion that we can draw based on the above analysis, where we interpreted the scheme as a central scheme plus (anti-)dissipation. The estimate \eqref{manzanero} is very convenient for the later numerical experiments \textcy{as it will allow us to identify and rule out Lax stable approximations} and we will refer to it as the \emph{MRFVK-estimate} (the initials of the authors).

\begin{remark}
\textit{We note that it is possible to design either energy stable schemes for the linear problem, or entropy-dissipative/conservative schemes for the non-linear problem. However, it seems that the core of the issue is that discretisation of the non-linear PDE and linearisation of the resulting scheme is not the same as linearisation of the PDE and then discretising the linearised problem.}
\end{remark}

In view of the above observations, we will use the lack of \textcy{local energy} stability as an indicator that a scheme is anti-dissipative. There are two options to study the linearised equations: 1) Derive an energy estimate. Only anti-diffusive terms can cause a growth of the norm in a periodic linear scheme in the absence of low-order terms. 2) Compute the eigenvalues of the spatial operators. \textcy{We consider problems that satisfy the MRFVK-estimate, i.e., with no exponential growth and $\alpha_c=0$, such that we know a priori that the spectrum must have non-positive eigenvalues.} Hence, we can determine if the scheme has anti-dissipative terms by examining if the eigenvalues have positive real parts.

\subsection{Analysis of high-order diagonal-norm split-form summation-by-parts operators}
\label{sec:dgsem_burgers}
In this subsection we focus on high-order diagonal-norm split-form summation-by-parts approximations with different entropic properties for the Burgers equation. In particular, we analyse the skew-symmetric discontinuous Galerkin spectral element method \cite{gassner_skew_burgers} that is strongly related to split-form summation-by-parts finite difference approximations with simultaneous approximation terms, e.g. \cite{skew_sbp2}.
We note that the entropy-conserving flux function is a key building block for the construction of high-order entropy-conserving/dissipative schemes, see e.g., \cite{LeFlochRohde00,skew_sbp2,carpenter_esdg,FisherCarpenter13,gassner_skew_burgers,Carpenter_etal16}.\medbreak

Details on the high-order discretisation can be found e.g. in \cite{gassner_skew_burgers,Gassner_BR1}. We give a short summary of the necessary building blocks. We split the computational domain into elements of equal size $h$. We map the elements to a unit reference element $[-1,1]$. In each reference element, the solution is approximated by a polynomial of degree $N$. We use $N+1$ Legendre-Gauss-Lobatto nodes to represent the polynomials with Lagrangian basis functions. Thus, the solution is represented as a vector of nodal values $\uvec$. For convenience we also introduce the diagonal injection of the nodal values to simplify notation, i.e. $\umat = \diag(\uvec)$. The corresponding main operators are the mass matrix $\mmat=\diag([\omega_0,...,\omega_N])$ (or integration matrix), the differentiation matrix $\dmat$ and the boundary matrix $\bmat=\diag([-1,0,...,0,1])$. The $\omega_i$ are the corresponding $N+1$ Legendre-Gauss-Lobatto integration weights in the reference element. The operators satisfy the summation-by-parts property
\begin{equation}
\mmat\,\dmat + (\mmat\,\dmat)^T=\bmat,
\end{equation}
where the mass matrix is interpreted as the diagonal-norm operator. We consider the split-form discontinuous Galerkin approximation
\begin{equation}
\label{eq:dgsem_burgers_split}
\frac{h}{2}\,\uvec_t + \alpha\frac{1}{2}\dmat\,\umat\,\uvec + (1-\alpha)\,\umat\,\dmat\,\uvec = -\mmat^{-1}\,\bmat\left[\fvec^* - \frac{1}{2}\,\umat\,\uvec\right],
\end{equation}
which follows directly from the continuous split-form \eqref{eq:split_alpha} of the Burgers equation. The numerical flux vector is
\begin{equation}
\fvec^* = [f_{i-1/2},0,...,0,f_{i+1/2}]^T,
\end{equation}
where $f_{i+1/2}$ is a consistent and unique numerical flux function, e.g. as used in the finite volume approach, see Section~\ref{sec:fv_burgers}.

We get an entropy-conserving semi-discretisation with the choice $\alpha=2/3$ (skew-symmetric form) and $\fvec^*$ taken to be the entropy-conserving numerical flux \eqref{eq:ec} for the quadratic entropy $U(u) = u^2/2$. We obtain entropy-dissipation when switching the numerical flux from entropy-conserving to entropy-dissipative, e.g. \eqref{eq:es}. Thus, these versions of the high-order discretisation have a global bound on the discrete $L^2$-norm, see e.g. \cite{gassner_skew_burgers}. 

Having defined the scheme, we proceed with the \textcy{local energy} stability analysis. As before, we consider a perturbation of a smooth baseflow $\uvec =\underline{\widetilde{u}} + \uvec'$ and assume that  $\underline{\widetilde{u}}$ solves \eqref{eq:dgsem_burgers_split}.  The analysis of the numerical flux function $f_{i+1/2}$ appearing in the surface integral, is similar to the finite volume discussion in Section~\ref{sec:fv_burgers}. Similarly, if we linearise the volume terms by inserting the ansatz with the fluctuations, neglecting high-order fluctuation terms, we obtain 
\begin{equation}
\begin{split}
&\left[\alpha\frac{1}{2}\dmat\,\umat\,\uvec + (1-\alpha)\,\umat\,\dmat\,\uvec\right] - \left[\alpha\frac{1}{2}\dmat\,\widetilde{\umat}\,\widetilde{\uvec} + (1-\alpha)\,\widetilde{\umat}\,\dmat\,\widetilde{\uvec}\right]\\[0.1cm]
\Rightarrow\,\,&\alpha\left[\dmat\,\widetilde{\umat}\,\uvec'\right] + (1-\alpha)\,\left[\widetilde{\umat}\,\dmat\,\uvec'+\umat'\,\dmat\,\widetilde{\uvec}\right].
\end{split}
\end{equation}
The analysis of these volume terms is analogous to the finite difference discussion in Section~\ref{sec:fd_burgers}. As mentioned already before, for $\widetilde{\underline{u}} > 0$ (for all elements), the MRFVK-estimate \eqref{manzanero} applies. The corresponding discrete estimate is only obtained when $\alpha = 1$, which is the central approximation with divergence form volume terms. All other choices, including the $\alpha=2/3$ entropy-conserving/dissipative choice, allow local artificial exponential growth of fluctuations, i.e., \emph{they are not \textcy{locally energy} stable}. 


We note that while the central approximation is \textcy{locally energy} stable, its deficiency is in the non-linear regime as the central approximation does not have a non-linear global bound. Hence, the central approximation is not preferable either.

\subsection{Numerical investigations for the non-linear Burgers equation}
\label{sec:results_burgers}

To numerically investigate our theoretical findings of the \textcy{local energy} stability issues, we consider the non-linear Burgers equation with periodic boundary conditions. \textcy{To apply the MRFVK-estimate, we choose a smooth initial data that is positive in the domain} 
\begin{equation}
\label{eq:linear_state_burgers}
\widetilde{u}(x) = \sin(\pi\,x-0.7)+2,\,\,x\in[-1,1],
\end{equation}
which will also serve as the linearisation state. \textcy{As the corresponding linearised variable coefficient problem satisfies the MRFVK-estimate, there is no (local in time) exponential growth of the solution for this test problem and hence the maximum real part of the eigenvalues is zero, i.e. $\alpha_c = 0$.} As discussed in the Sections \ref{sec:fv_burgers} and \ref{sec:fd_burgers}, positive eigenvalues of the spectrum of the linearised spatial operator \textcy{for this particular test case indicate exponential growth exceeding the growth of the continuous problem (artificial anti-dissipation) and that the scheme is Lax stable at most.}

\begin{figure}[!htbp]
  \centering
 \includegraphics[width=0.55\textwidth,trim=0 0 0 0,clip]{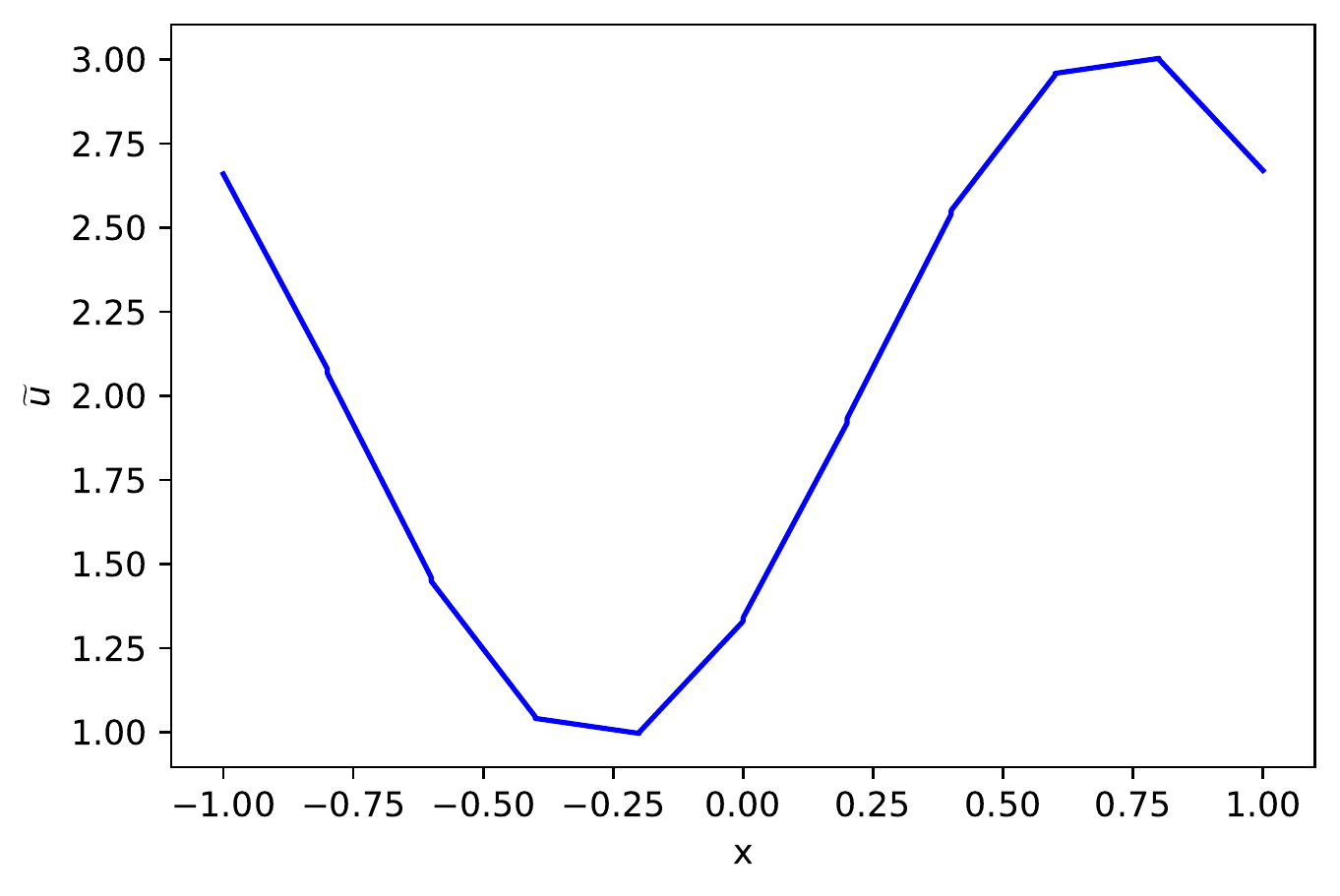}
\caption{Plot of the baseflow \eqref{eq:linear_state_burgers} for a discretisation with $10$ elements and polynomial degree $\widetilde{N}=1$ in every element. Note that there are small jumps at the element interfaces.}
\label{fig:linearisation_state_burgers}     
\end{figure}

Figure~\ref{fig:linearisation_state_burgers} shows an illustration of the baseflow for the investigations used below. To mimic mild under-resolution of the baseflow, we project the linearisation state \eqref{eq:linear_state_burgers} onto piece-wise linear polynomial basis functions ($\widetilde{N}=1$) for all Burgers simulation results presented below. As can be seen in Figure~\ref{fig:linearisation_state_burgers}, the baseflow is discontinuous across element interfaces, although the jumps are relatively small as the baseflow is only mildly under-resolved. 

\subsubsection{Comparison of the spectra of the semi-discrete operator}
 
To compute the spectrum of a spatial operator, we use the typical finite difference approach to approximate the linearisation and determine the $j$-th column of the Jacobi matrix of the semi-discrete residual as
\begin{equation}
A\underline{e}_j \approx  \frac{\underline{\rhs}(\widetilde{\underline{u}}+\underline{e}_j\,\epsilon) - \underline{\rhs}(\widetilde{\underline{u}}-\underline{e}_j\,\epsilon)}{2\,\epsilon}
\end{equation}
where $\epsilon = 10^{-8}$ and $\evec_j$ is the unit vector in the $j$-th direction. With the finite difference approach, the accuracy of the eigenvalues are limited to about single machine precision and hence, we interpret eigenvalues within this accuracy constraint as zero, even when formally the number is positive. The computed eigenvalues of the approximative Jacobi matrix $A$ give an approximate spectrum of the semi-discrete high-order operator. We are particularly interested in assessing the maximum real part across all eigenvalues as positive eigenvalues imply exponential growth of their eigenmodes. 

\begin{figure}[!htbp]
  \centering
 \includegraphics[width=0.48\textwidth,trim=0 0 0 0,clip]{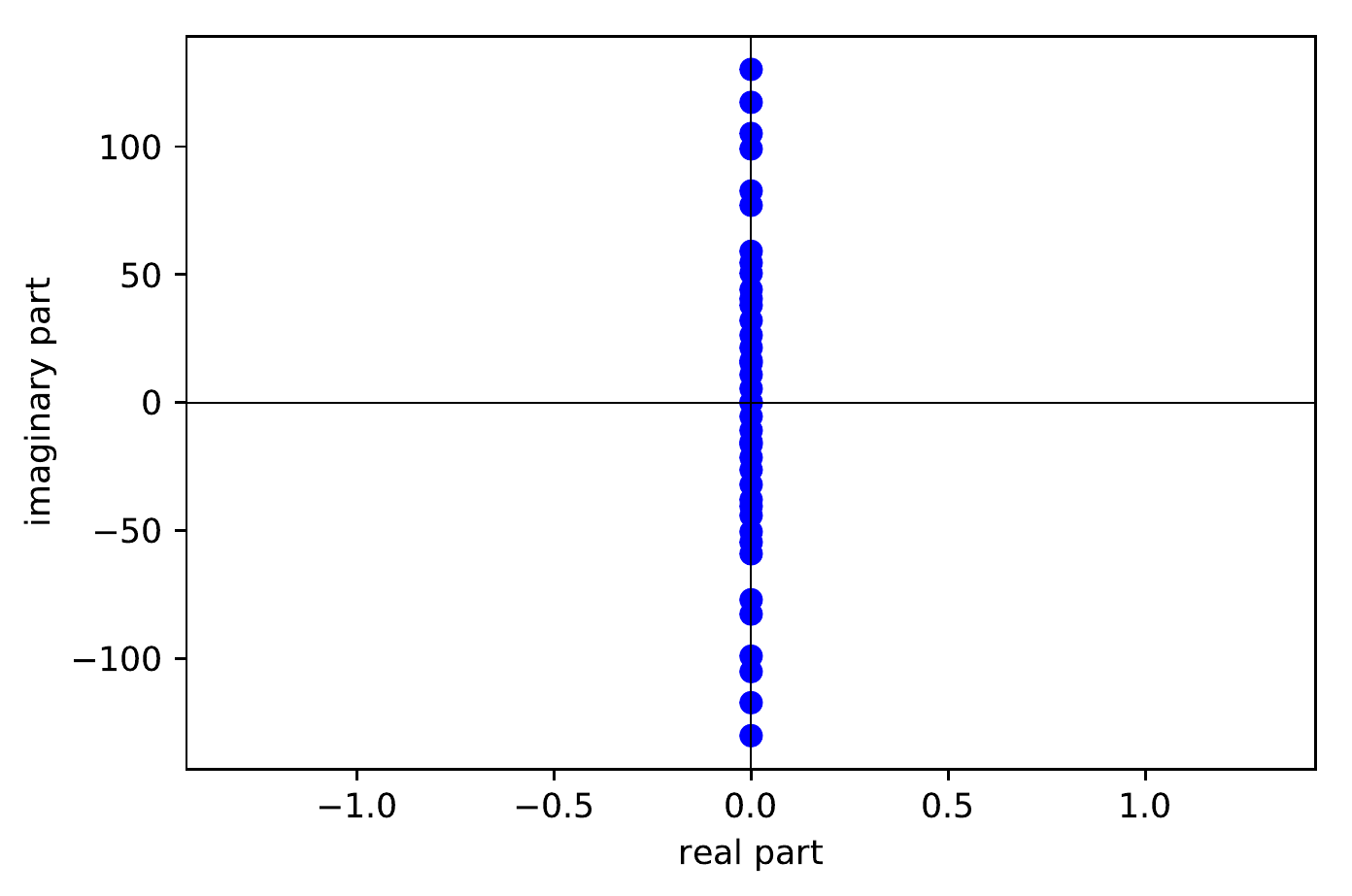} \includegraphics[width=0.48\textwidth,trim=0 0 0 0,clip]{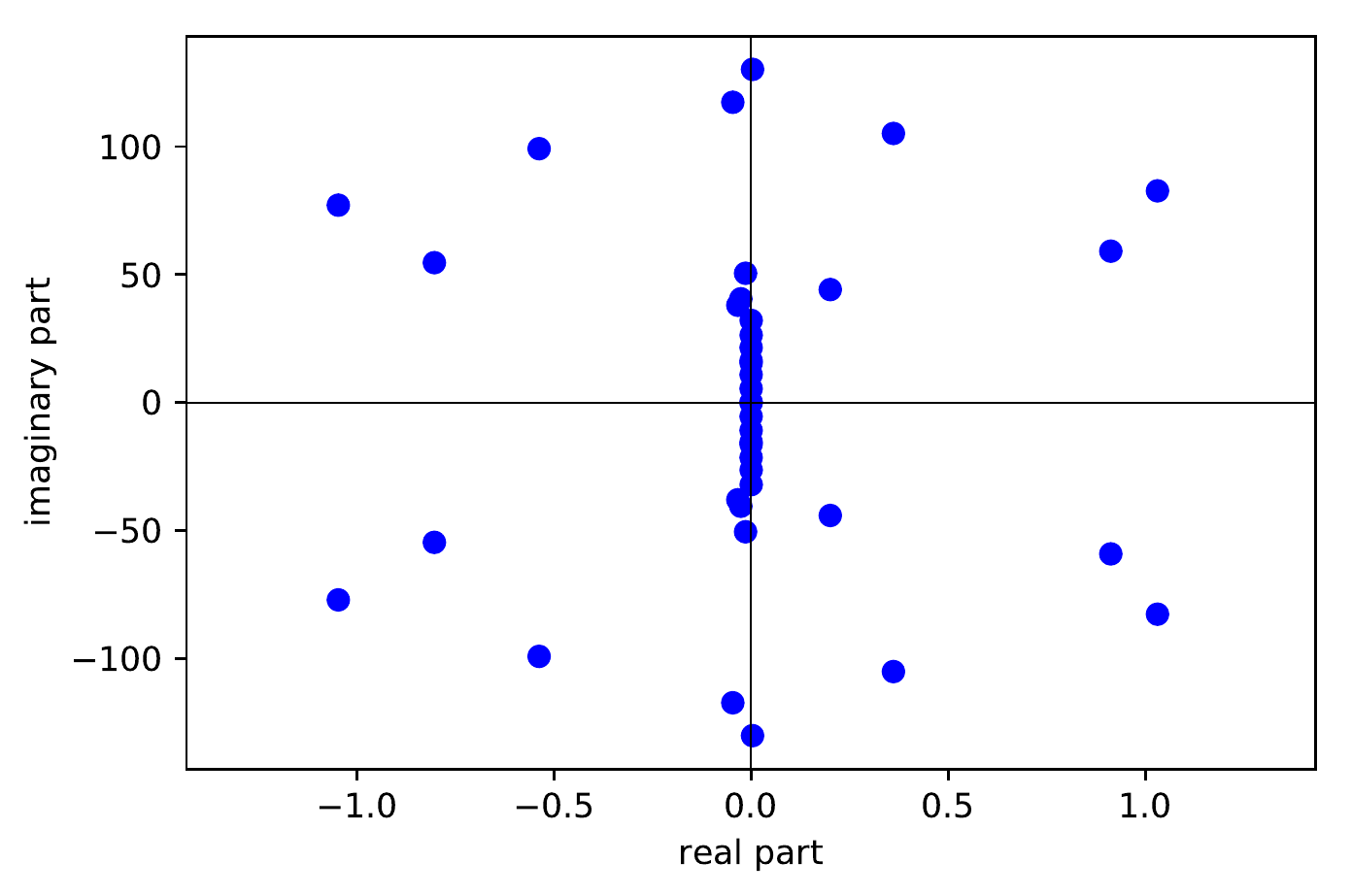}
\caption{
Spectra of spatial operators. Left: Central approximation with the divergence form ($\alpha=1$) and the central numerical flux. Right: Entropy-conserving approximation with skew-symmetric form ($\alpha=2/3$) and the entropy-conserving flux \eqref{eq:ec}. \figinfo{Discretisation with polynomial degree $N=3$ and $10$ elements with the baseflow shown in Figure~\ref{fig:linearisation_state_burgers}. The maximum real part in of the approximate spectrum is about $8.8\times10^{-8}$ for the central scheme and about $1.0307$ for the entropy-conserving scheme.}
}
\label{fig:spectrum_burgers_1}     
\end{figure}

For the approximation \eqref{eq:dgsem_burgers_split}, we plot the results of the fully central approximation ($\alpha=1$) and the entropy-conserving approximation ($\alpha=2/3$) in Figure~\ref{fig:spectrum_burgers_1}. The plots reveal a significant positive real part of $1.0307$ for the entropy-conserving approximation, whereas the spectrum of the central approximation is neutrally stable and almost purely imaginary in accordance with \eqref{manzanero}.\medbreak

\begin{figure}[!htbp]
  \centering
 \includegraphics[width=0.55\textwidth,trim=0 0 0 0,clip]{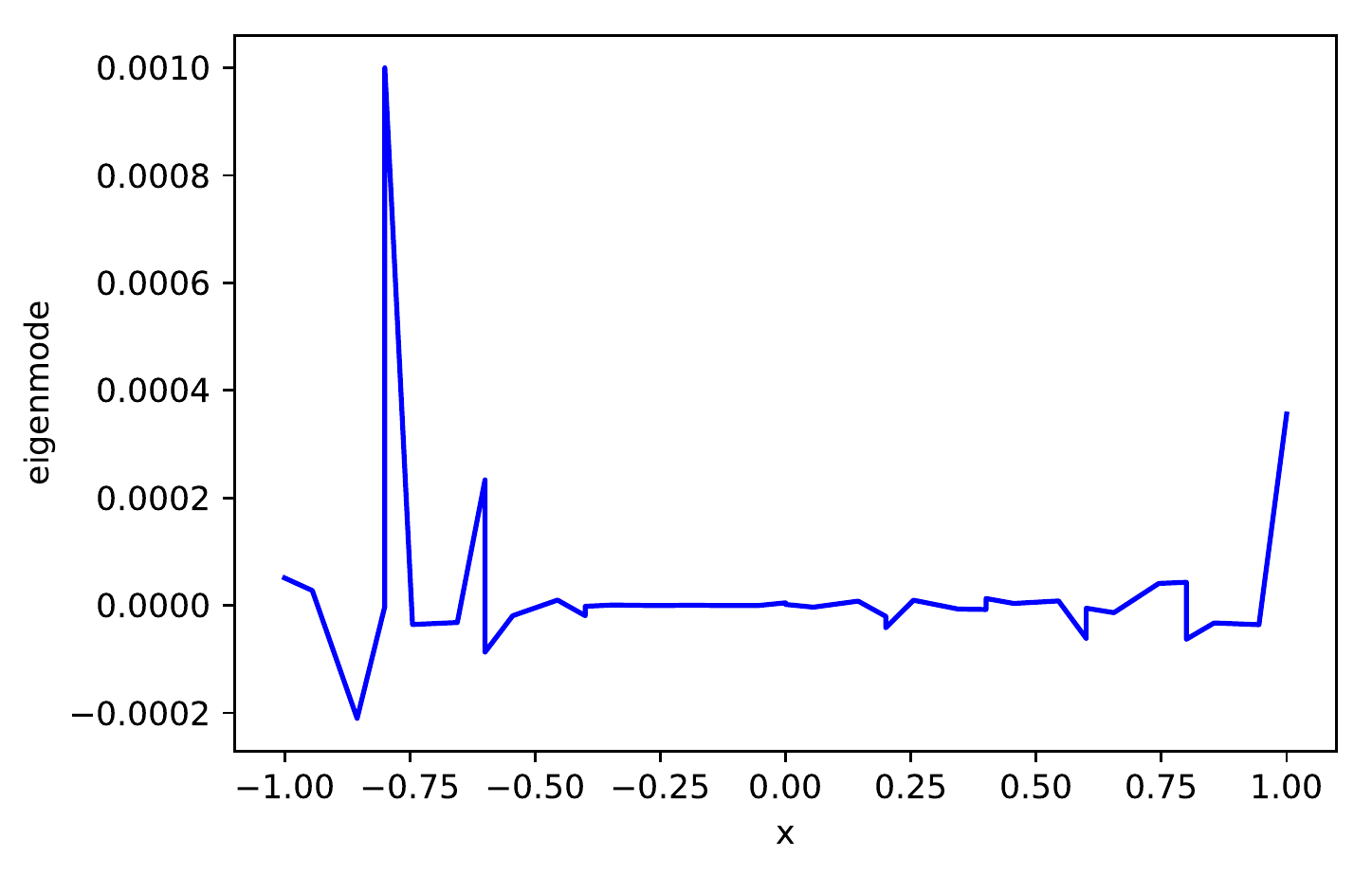}
\caption{Plot of the eigenmode corresponding to the eigenvalue with the largest real part (about $1.0307$) for the entropy-conserving discretisation with skew-symmetric form ($\alpha=2/3$) and the entropy-conserving flux \eqref{eq:ec}. \figinfo{The discretisation uses $10$ elements and a polynomial degree $N=3$, with the baseflow \eqref{eq:linear_state_burgers}. The eigenmode is scaled such that its largest peak is $0.001$.}} 
\label{fig:critical_eigenmode_burgers}     
\end{figure}

With a maximum real part of about $1.0307$ in the defective spectrum, we should expect a strong exponential growth of the corresponding eigenmode. We plot the corresponding eigenmode in Figure~\ref{fig:critical_eigenmode_burgers}. We note that the eigenmode is active in the part of the domain, where the slope of the baseflow is negative. Although unclear if there is a direct connection, it is interesting to note that this coincides with the discussion on the entropy-conserving split-form in the previous Section~\ref{sec:fd_burgers}, where the dissipation coefficient gets negative, when the slope of the linearisation state is negative.\medbreak
 
\begin{figure}[!htbp]
  \centering
 \includegraphics[width=0.48\textwidth,trim=0 0 0 0,clip]{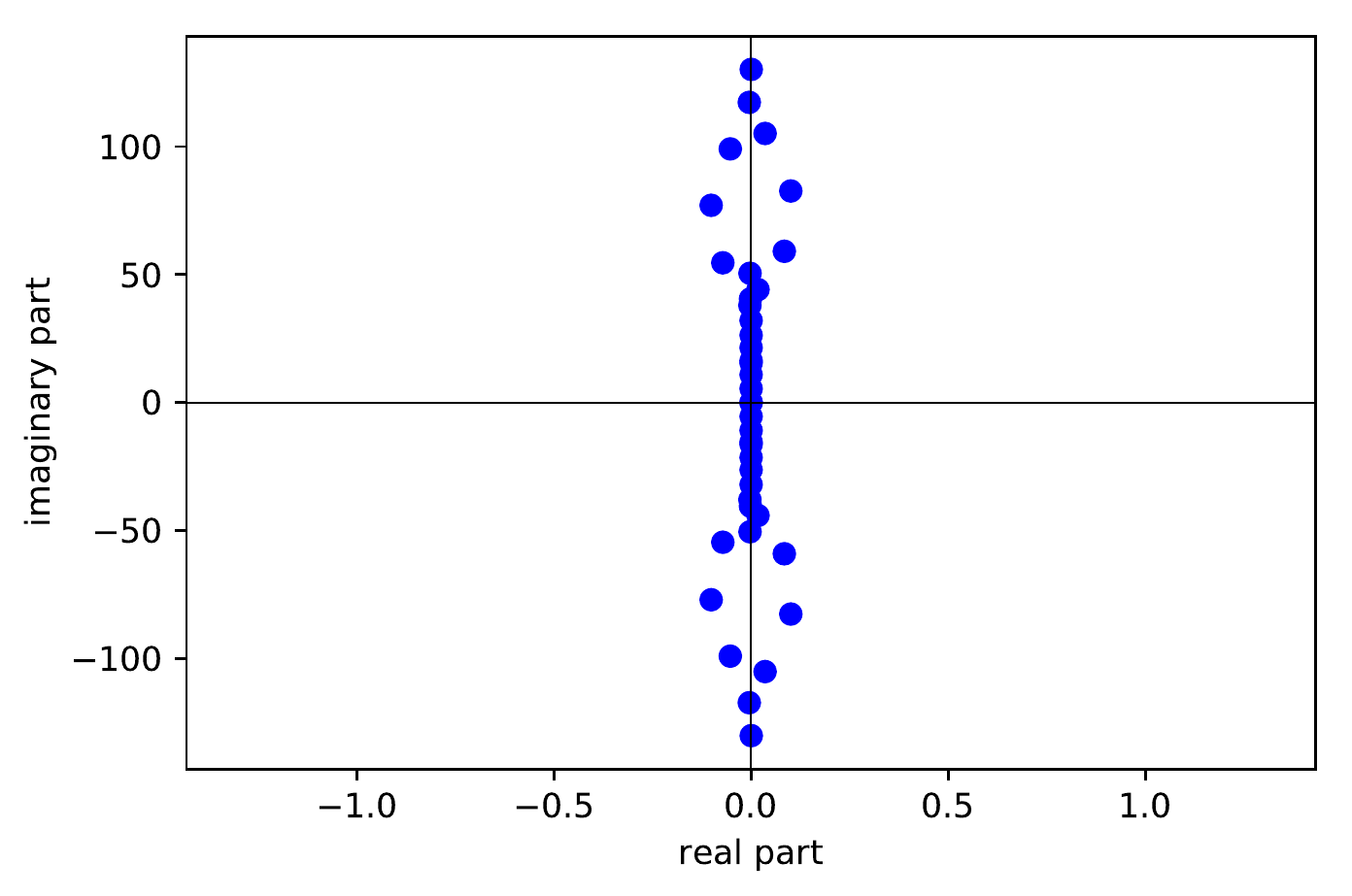} \includegraphics[width=0.48\textwidth,trim=0 0 0 0,clip]{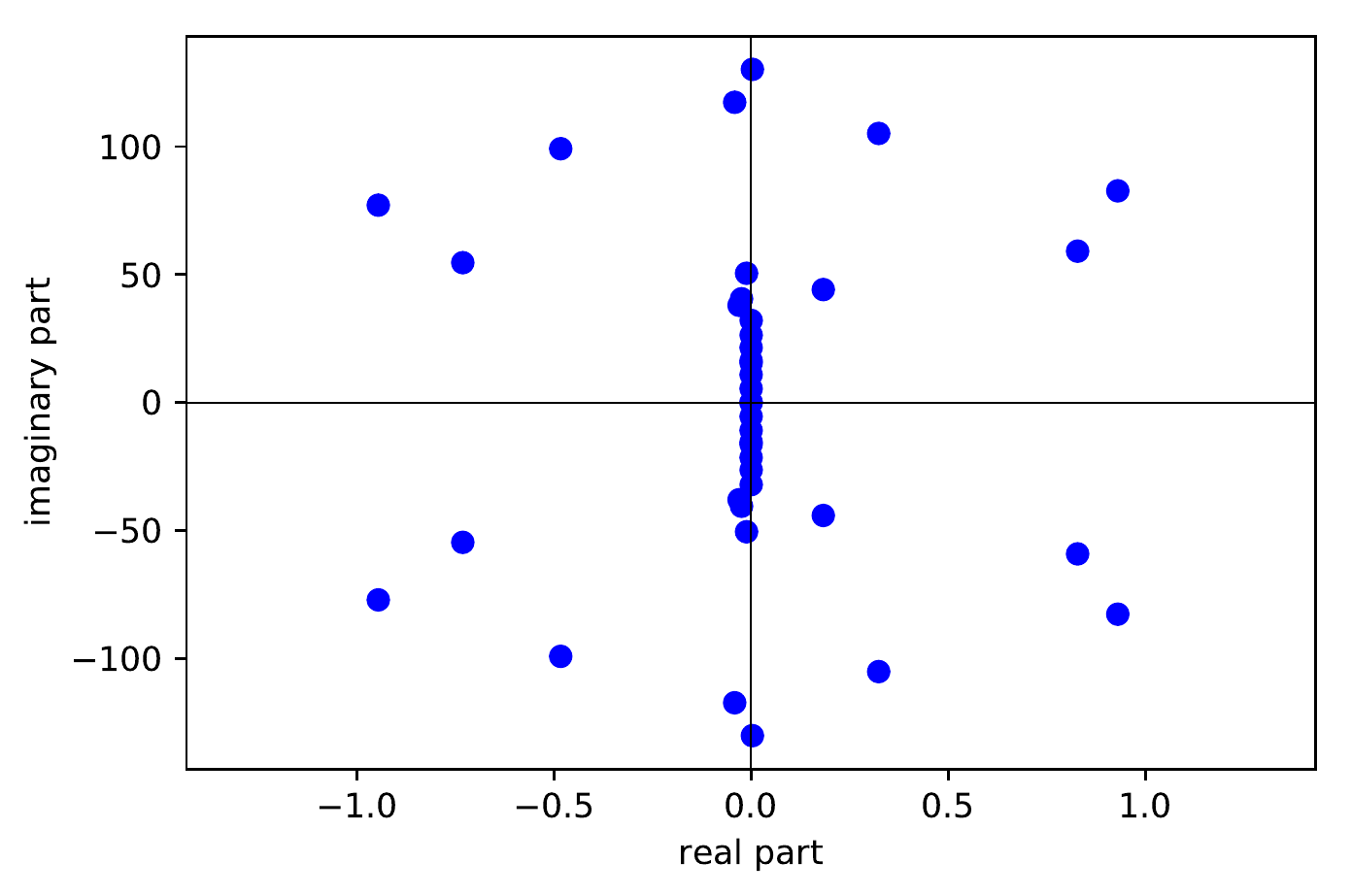}
\caption{
Spectra of spatial operators. Left: Volume term in divergence form ($\alpha=1$), surface term uses entropy-conserving flux \eqref{eq:ec}. Right: Volume term in skew-symmetric form ($\alpha=2/3$), surface term uses central flux.
\figinfo{Discretisation with polynomial degree $N=3$ and $10$ elements with the baseflow shown in Figure~\ref{fig:linearisation_state_burgers}.  The maximum real part is about $0.1006$ (left) and  $0.9300$ (right). }
}
\label{fig:spectrum_burgers2}     
\end{figure}

Next, we investigate if either the volume terms or the surface terms are the cause for this faulty and undesired behaviour of the discretisation. Figure~\ref{fig:spectrum_burgers2} shows the result for the divergence form of the volume terms ($\alpha=1$) combined with the entropy-conserving flux \eqref{eq:ec} in the left plot. In the right plot, the discretisation ingredients are flipped with skew-symmetric ($\alpha=2/3$) volume terms and the central numerical flux for the surface terms. It is clearly seen that the spectra still have eigenvalues with large positive real parts, here the maximum is $0.1006$ in the left plot and $0.93$ in the right plot. The exponential growth rate with the entropy-conserving surface flux is smaller than for the fully entropy-conserving scheme. The reason is, that the baseflow is relatively well resolved, with only small jumps at element interfaces. As the anti-dissipation scales directly with this jump, its effect is decreased. The plot on the right also suggests that the volume terms in entropy-conserving form are indeed a major contributor to the lack of \textcy{local energy} stability. \medbreak

\begin{figure}[!htbp]
  \centering
 \includegraphics[width=0.48\textwidth,trim=0 0 0 0,clip]{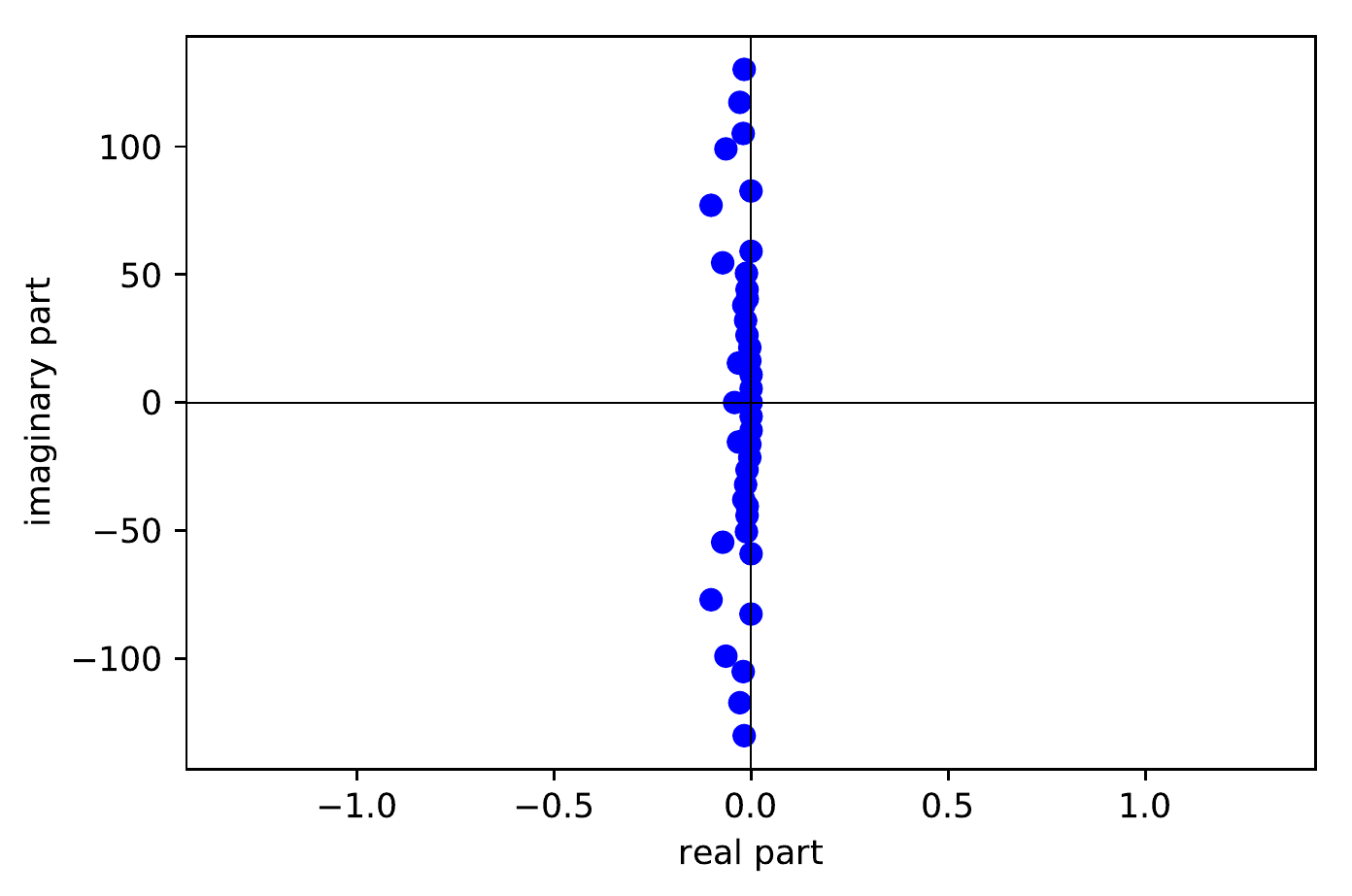} \includegraphics[width=0.48\textwidth,trim=0 0 0 0,clip]{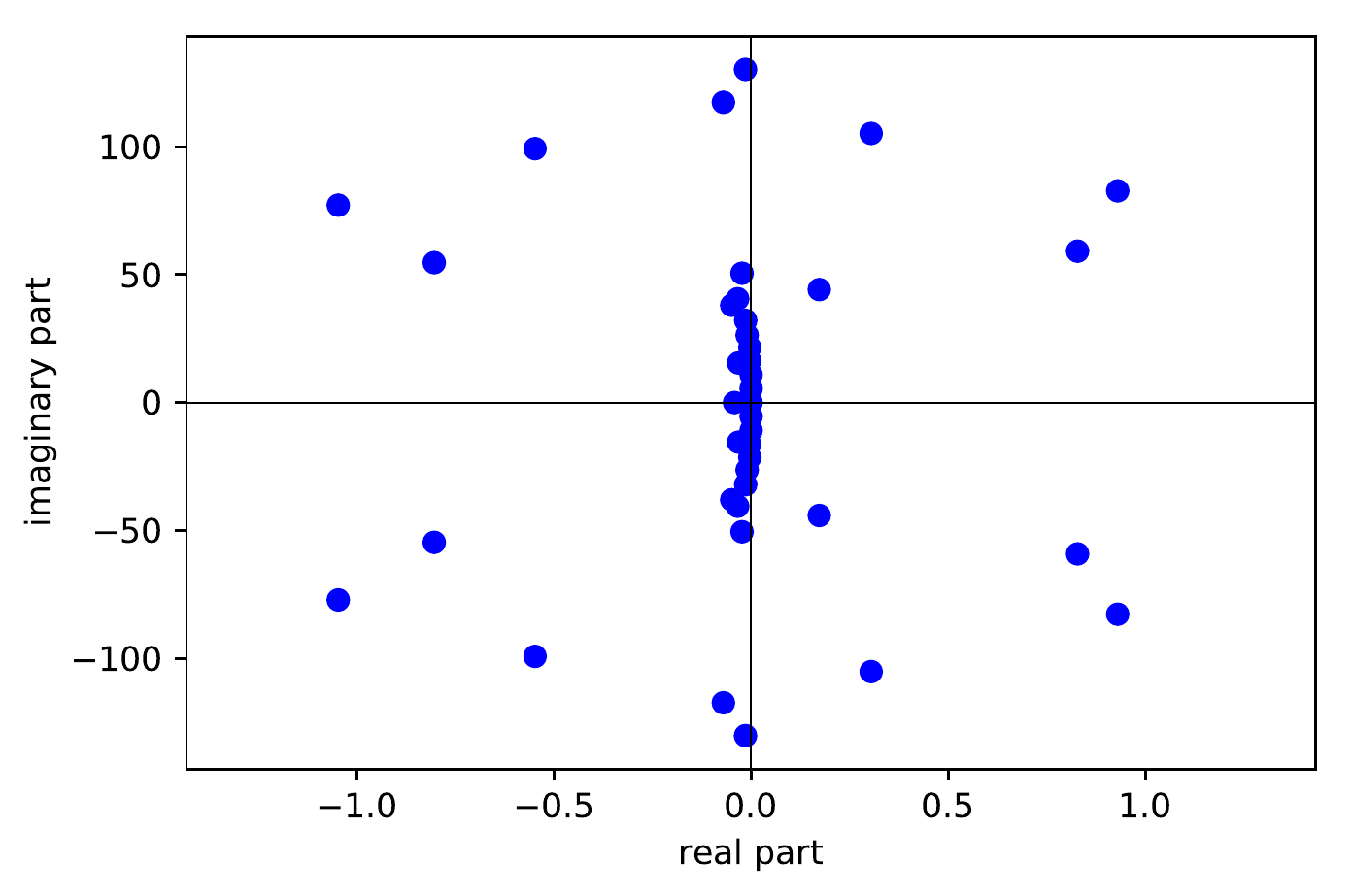}
\caption{
Spectra of spatial operators.  Left: Volume term in divergence form ($\alpha=1$), surface term uses the modified entropy-stable flux that excludes anti-dissipation \eqref{burgers_stable2}. Right: The volume terms are in skew-symmetric form ($\alpha=2/3$), the surface term uses the modified entropy-dissipative flux that excludes anti-dissipation \eqref{burgers_stable2}. 
\figinfo{Discretisation with polynomial degree $N=3$ and $10$ elements with the baseflow shown in Figure~\ref{fig:linearisation_state_burgers}. The maximum real part is about $-9.01\times10^{-8}$ (left) and  $0.9298$ (right).} 
}
\label{fig:spectrum_burgers3}     
\end{figure}

We further investigate the impact of the modified entropy-conserving flux by Tadmor \eqref{burgers_stable2} that excludes the anti-dissipation. To demonstrate that it has a stabilising effect on the scheme, we redo the plots of Figure \ref{fig:spectrum_burgers2}. The results are illustrated in Figure~\ref{fig:spectrum_burgers3} and suggest that indeed Tadmor's modification seems to get rid of the faulty influence of the anti-dissipation (left plot), when combined with divergence form volume terms. This corroborates our theoretical findings that the positive part of the spectrum is associated with the anti-diffusive character of the scheme. However, it is also clear, that this particular modification of the surface term is not enough to fix the \textcy{locally energy} unstable behaviour of the skew-symmetric volume terms, as can be seen in the right plot, which still has a significant positive real part, due to the anti-diffusion in the volume terms.

\begin{remark}\label{carpe}\textit{
It is interesting to note that Carpenter et al. considered a modification of the volume terms to introduce entropy-dissipation, where they use entropy generation of the central approximation compared to the entropy-conserving approximation to adjust the volume terms (see \cite{carpenter_esdg}, Section~5, equation (5.6)). Unfortunately, our numerical investigations show only little influence on the spectra for this modification, as can be seen in Figure~\ref{fig:spectrum_burgers4}}. 
 \begin{figure}[!htbp]
  \centering
 \includegraphics[width=0.5\textwidth,trim=0 0 0 0,clip]{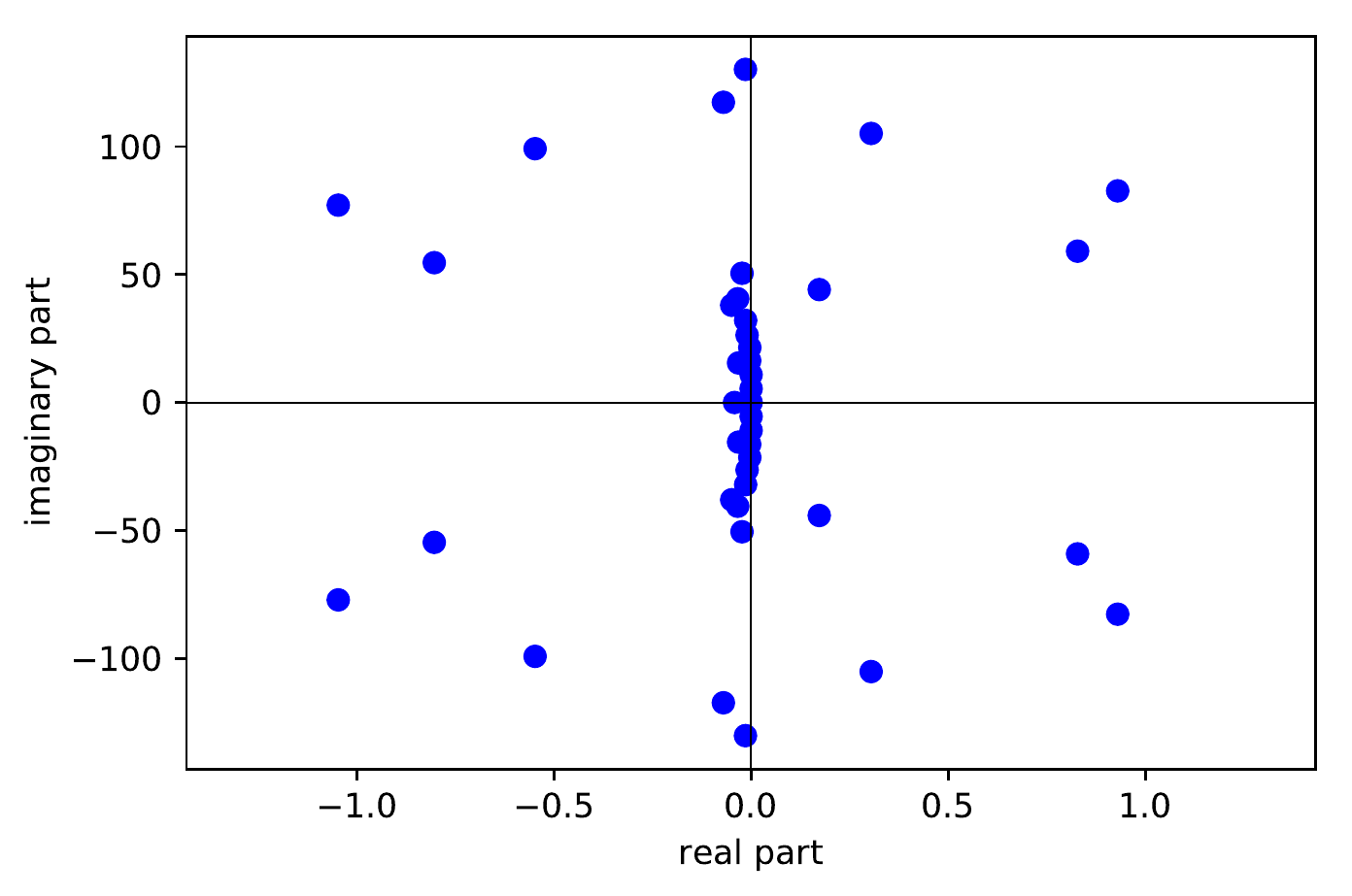}
 \caption{
 Spectrum of spatial operator. The volume terms are entropy-dissipative according to the modification of Carpenter et al. (see  \cite{carpenter_esdg}, Section~5, equation (5.6)). The surface flux is the modified version of Tadmor \eqref{burgers_stable2}, that excludes anti-dissipation. 
 \figinfo{Discretisation with polynomial degree $N=3$ and $10$ elements with the baseflow shown in Figure~\ref{fig:linearisation_state_burgers}. The maximum real part is about $0.9297$.}
 }
\label{fig:spectrum_burgers4}     
\end{figure}
\textit{In particular, we still observe significant positive real parts of about $0.9297$, similar to the case with full skew-symmetric volume terms. This highlights that a sufficient amount of dissipation is necessary to pull the spectrum to the left-hand side.}
\end{remark}

\begin{figure}[!htbp]
  \centering
 \includegraphics[width=0.48\textwidth,trim=0 0 0 0,clip]{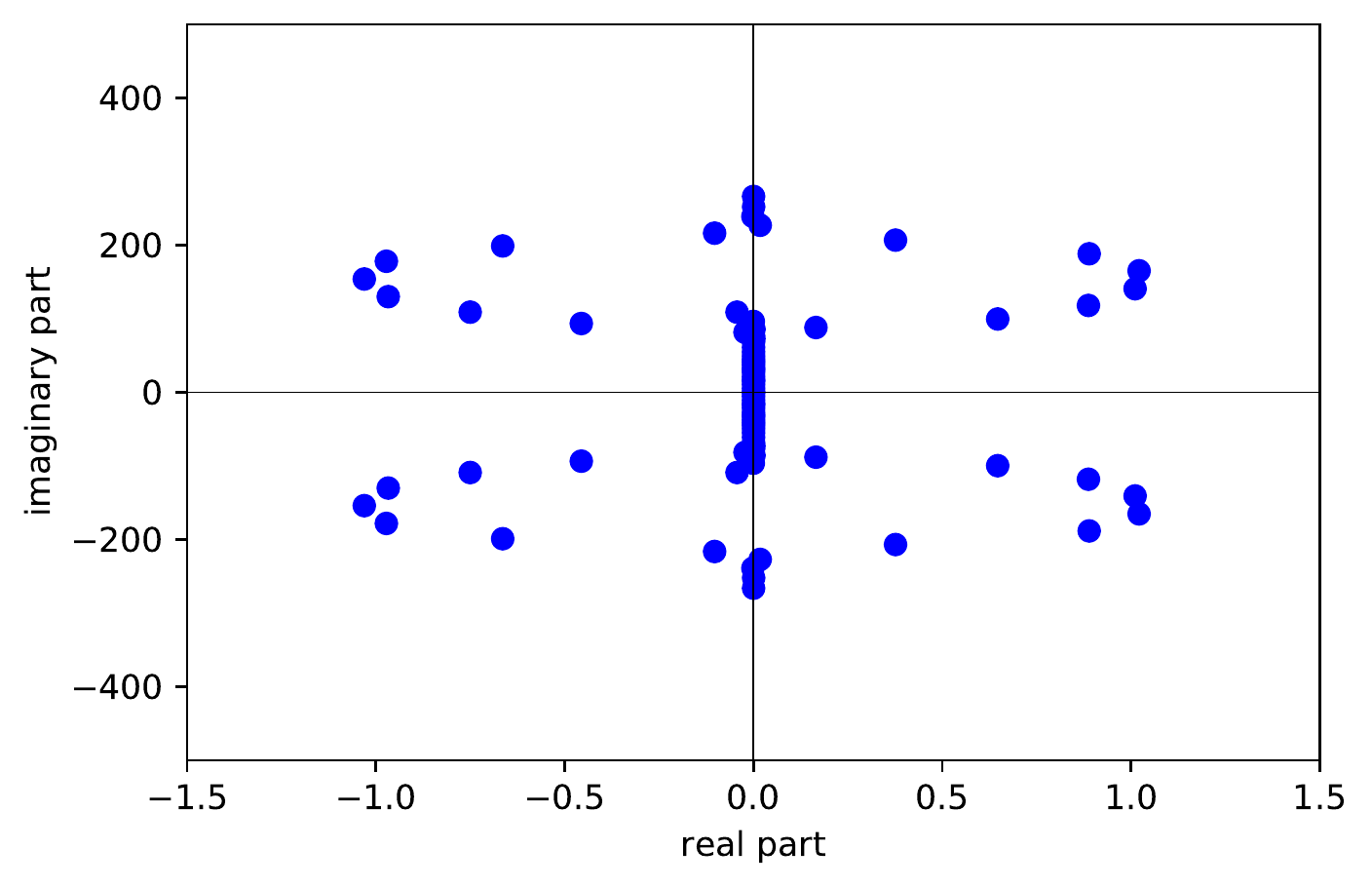} \includegraphics[width=0.48\textwidth,trim=0 0 0 0,clip]{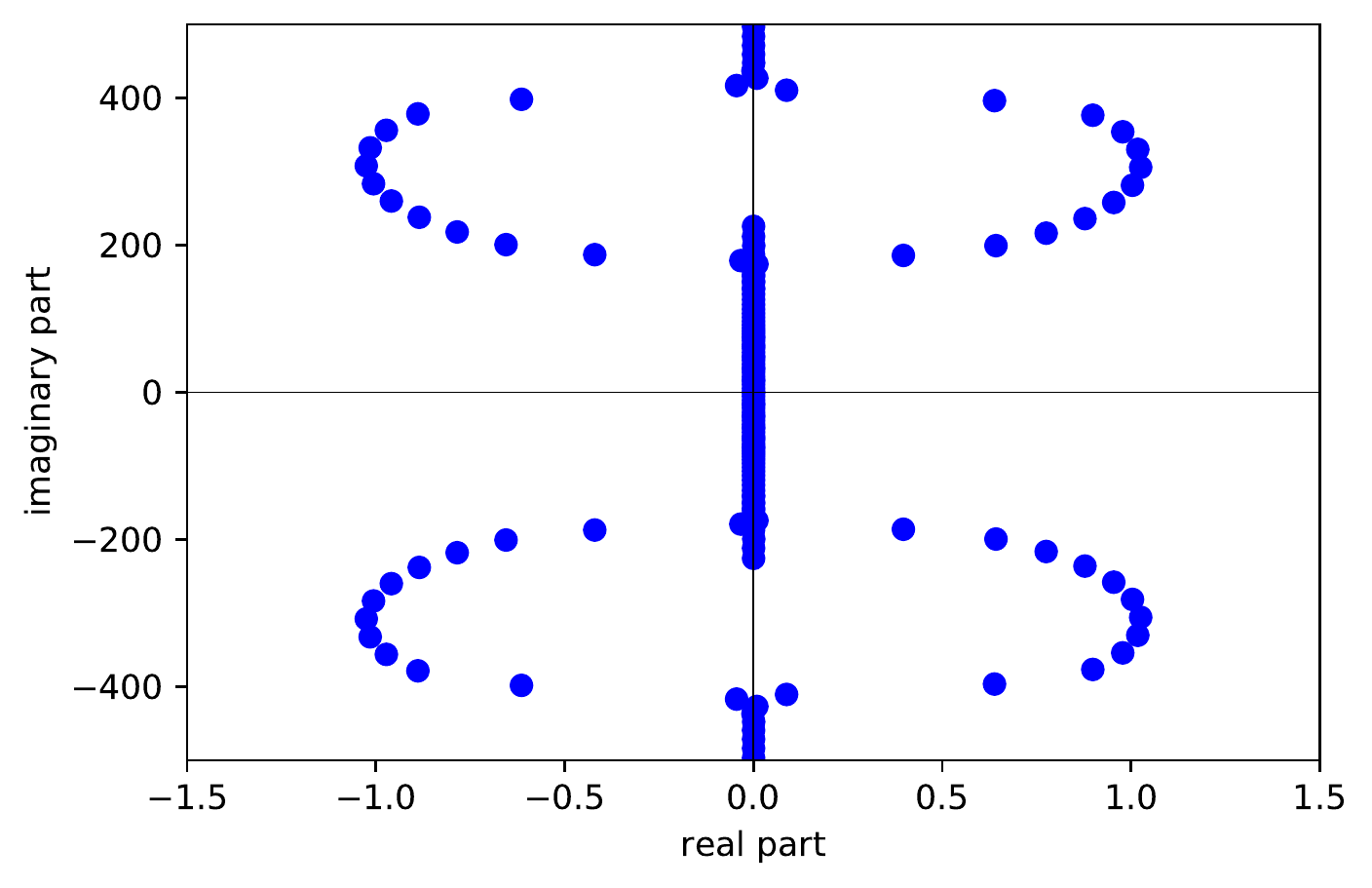}
\caption{
Spectra of spatial operator, for entropy-conserving approximation with skew-symmetric form ($\alpha=2/3$) and the entropy-conserving flux \eqref{eq:ec}.  Left: Number of elements is $20$. Right: Number of elements is $40$.
\figinfo{Discretisation with polynomial degree $N=3$. Linearisation with the baseflow shown in Figure~\ref{fig:linearisation_state_burgers}.  Maximum real part is $1.021$ (left) and $1.025$ (right). }
}
\label{fig:spectrum_burgers_111}     
\end{figure}

Next, we investigate the impact of resolution on the behaviour of the spectrum. As discussed in the previous Section~\ref{sec:fd_burgers}, we are able to interpret the skew-symmetric forms as central discretisations with dissipation. The dissipation coefficient can be positive or negative and scale with the mesh size $h$, which suggests that refining the mesh (decreasing $h$) should improve the faulty behaviour. We plot in Figure~\ref{fig:spectrum_burgers_111} the resulting spectra of the entropy-conserving skew-symmetric form when increasing the number of elements from $10$ (compare right part of Figure~\ref{fig:spectrum_burgers_1}) to $20$ and $40$ respectively. One can clearly see that the maximum positive real part does not shrink and stays at the value of about $1$. However, we can also clearly observe that the faulty eigenmodes get pushed to higher frequencies, as the corresponding imaginary parts increase with decreasing mesh size. We note that we observe a similar behaviour for all discretisations discussed in this work: the faulty eigenmodes always get pushed to higher frequencies when refining the mesh, whereas the maximum real part stays about the same. \textcy{These numerical investigations suggest that the high-order schemes might be Lax stable. Depending on the initial condition}, we expect that for very well resolved smooth problems, the content of high frequency eigenmodes decays quickly. Thus, the contribution of the faulty high frequency eigenmodes should also be reduced when refining the mesh, compare the discussion at the end of Section~ \ref{sec:linear}, \textcy{and convergence could be achieved}. We refer to \textcy{Section~\ref{sec:burgers_sim_growth}} for further numerical results that support this expectation, see the results in Figure~\ref{fig:longtime_burgers3} and Figure~\ref{fig:longtime_burgers4}.\medbreak

\begin{figure}[!htbp]
  \centering
 \includegraphics[width=0.48\textwidth,trim=0 0 0 0,clip]{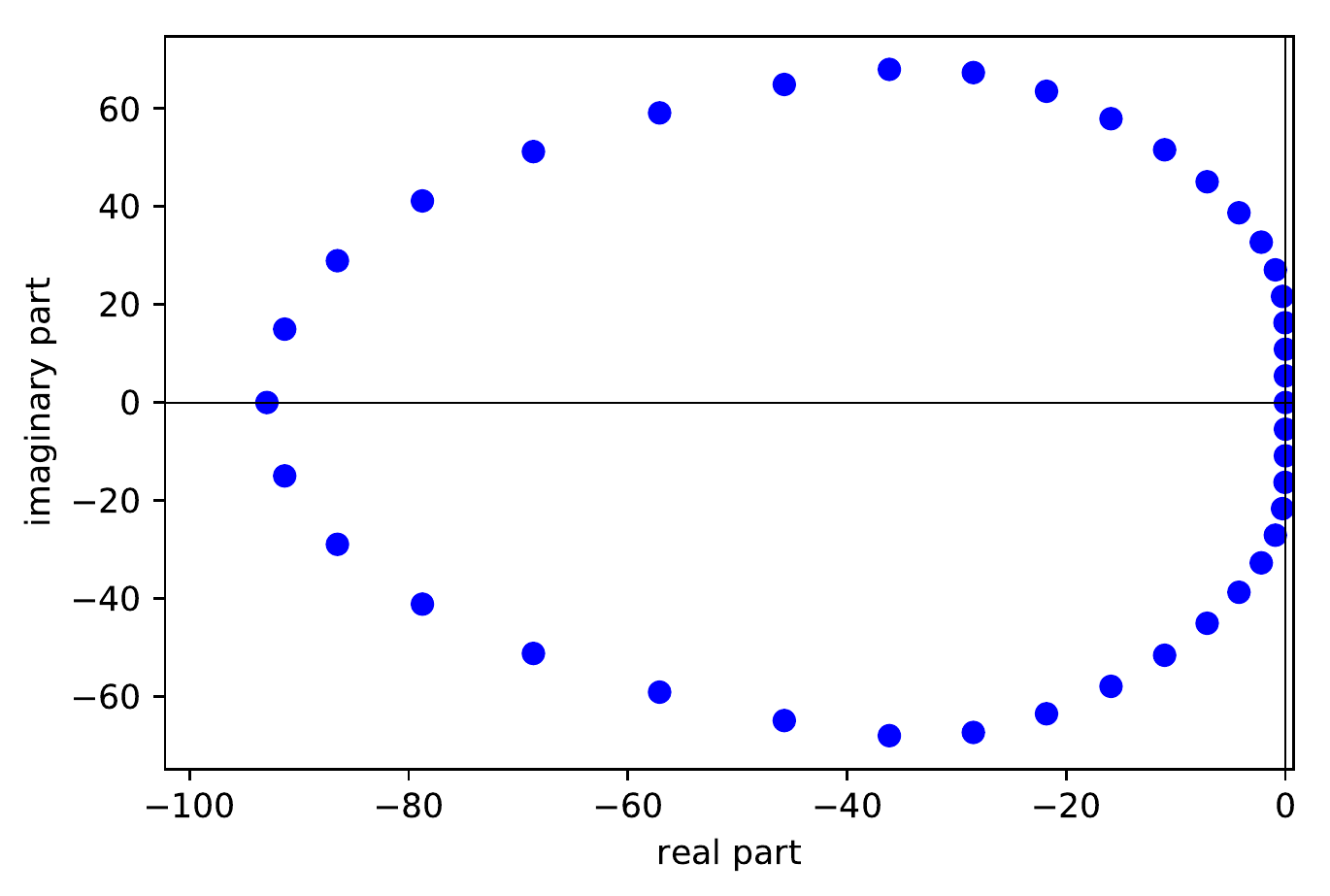} \includegraphics[width=0.48\textwidth,trim=0 0 0 0,clip]{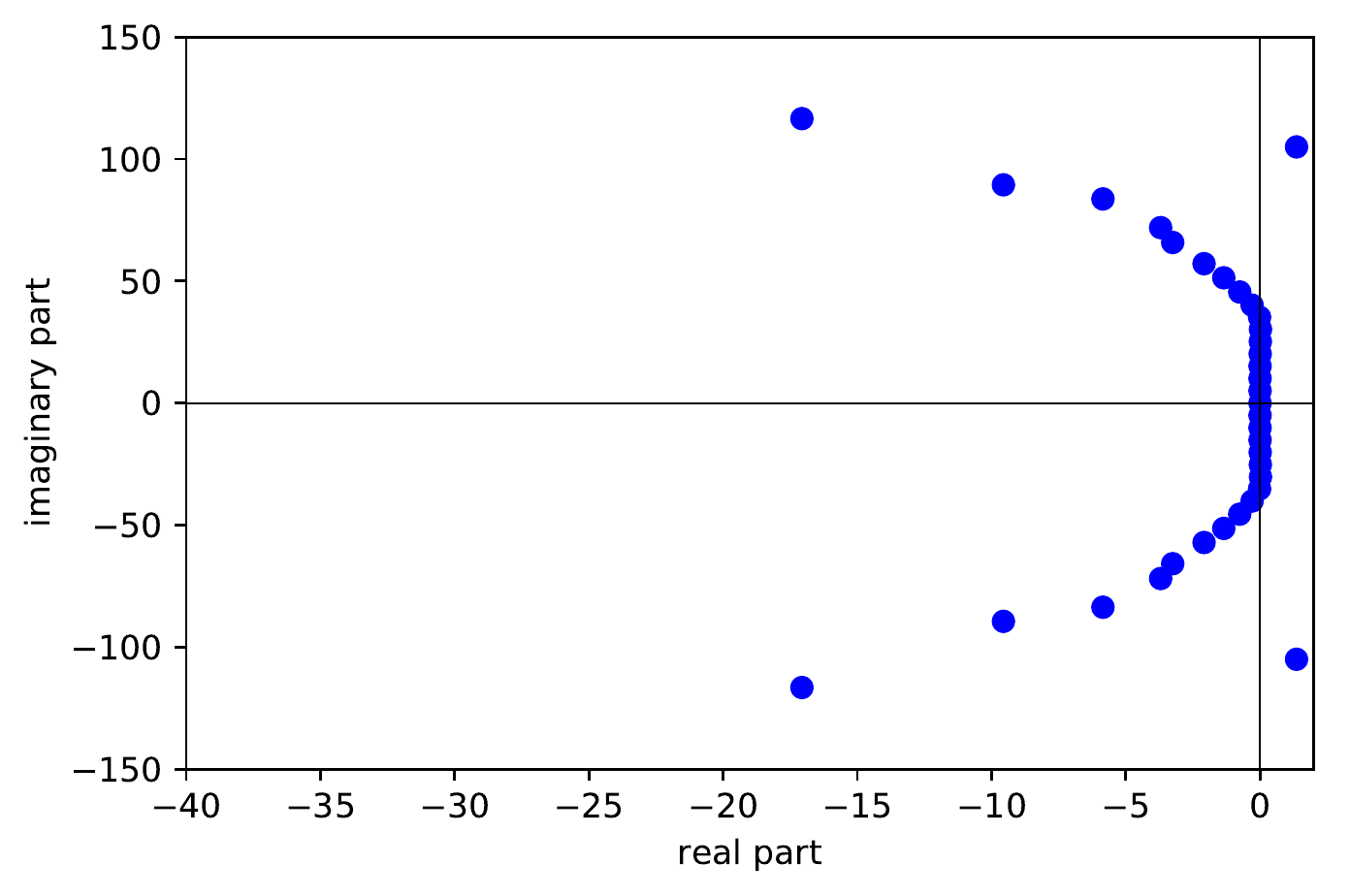}
\caption{
Spectra of spatial operator, the surface term uses entropy-dissipative Rusanov-type flux \eqref{eq:es}, the volume term is in skew-symmetric form ($\alpha=2/3$). Left: Discretisation with polynomial degree $N=3$ and $10$ elements and linearisation with baseflow shown in Figure~\ref{fig:linearisation_state_burgers}. Maximum real part is $-1.06\times10^{-7}$. Right: Discretisation with polynomial degree $N=15$ and $3$ elements, and an increased baseflow frequency from $\pi$ to $4\,\pi$, to mimic severe under-resolution. Zoom in on the spectrum, to show the eigenvalues with maximum real parts of $1.359$.
 }
\label{fig:spectrum_burgers5}     
\end{figure}

Lastly, we investigate the effect of a strongly dissipative numerical surface flux function. As in most practical simulations, the skew-symmetric volume terms are combined with a dissipative numerical surface flux, such as e.g. the Rusanov-type entropy-dissipative flux \eqref{eq:es}. From a theoretical point of view, it is not obvious how a surface dissipation term can apply control over volume type anti-dissipation. However, as can be seen in the left part of Figure~\ref{fig:spectrum_burgers5}, the spectrum gets shifted towards the negative real axis, with the maximum real part of now $-1.06\times10^{-7}$.\medbreak

Unfortunately, it is possible to design a discretisation setup where the surface dissipation is not enough to stabilise the skew-symmetric volume integrals. To provoke such a behaviour, we first increase the polynomial degree from $N=3$ to $N=15$ and decrease the number of elements from $10$ to $3$. This gives about the same number of degrees of freedom (48 in comparison to 40), however the number of element interfaces is smaller, which should obviously decrease the stabilising effect of the surface terms. We note that it is especially the high-order polynomial approximations that are attractive for under-resolved turbulence simulations, e.g. \cite{Gassner:2013qf}. Second, we mimic severe under-resolution of the baseflow by increasing the frequency of the baseflow \eqref{eq:linear_state_burgers} from $\pi$ to $4\,\pi$. In combination with the very low number of elements, this results in large jumps across the element interfaces. A zoom in of the resulting spectrum is plotted in the left part of Figure~\ref{fig:spectrum_burgers4}, where the maximum real part of the eigenvalue is $1.359$, which means that the entropy-dissipative scheme in this setup is again \textcy{locally energy} unstable and allows for exponential growth.\medbreak

As an intermediate conclusion, the numerical results confirm the theoretical findings and the faulty behaviour of the entropy-conserving numerical flux and the entropic split-form approach with diagonal-norm summation-by-parts operators. The schemes allow for a nonphysical exponential growth of local fluctuations. Adding a dissipative mechanism through the numerical fluxes at the interfaces does not necessarily pull the spectrum to the negative half-plane. Neither does necessarily adding dissipation to the volume terms. Both these "fixes", at best only mask the inherent issue. We suspect that sufficient additional artificial dissipation in the volume terms in combination with truly entropy-dissipative surface fluxes might produce locally correct behaviour. However, that comes at the cost of being over-dissipative and reducing the benefit of a high-order, low dissipation scheme.

On the other hand, the central approximation in these particular investigations is \textcy{locally energy} stable. However, the central approximation of course lacks a non-linear (entropy) stability estimate and hence is not suitable for the simulation of highly non-linear (under-resolved) problems. 

\subsubsection{Simulation of the growth of the fluctuations}
\label{sec:burgers_sim_growth}

The goal of this section is to solidify the numerical results from the previous section and perform actual simulations to monitor the growth of the fluctuations. The semi-discretisation \eqref{eq:dgsem_burgers_split} of the non-linear Burgers equation gives a system of weakly coupled ordinary differential equations in time that we integrate with a three-stage third-order strong stability preserving Runge-Kutta scheme with a time step computed as 
\begin{equation}
\Delta t = \CFL\,\frac{h}{(N+1)\,\lambda_{\max}},
\end{equation}
where $\lambda_{\max}=\max\limits_{u}\,f'(u)$ and $\CFL=0.05$ to keep the time integration errors low.\medbreak

We particularly focus on validating the previous results and are hence interested in the evolution of the fluctuations. 
It is clear that the function shown in Figure~\ref{fig:linearisation_state_burgers} is not a solution of the non-linear homogeneous Burgers equation. Instead, we consider the inhomogeneous Burgers equation, where we use the residual computed with the linearisation state as a discrete source term. For the inhomogenous source term, we first compute the fully non-linear semi-discrete residual with the linearisation state $\underline{\rhs}(\widetilde{u})$ as a pre-processing step. Then, during the actual simulation, we substract in every single Runge-Kutta stage the residual of the linearisation state. It is important to note, that we project the linearisation state to piecewise linear polynomials, to guarantee that for a DGSEM approximation with $N=3$ or higher, the volume terms in the residual computations are equal for all split-forms including the divergence form and the skew-symmetric form. By solving the inhomogenous Burgers equation with this discrete source term, we make sure that a simulation with initial condition equal to the linearisation state would give zero growth down to machine precision. This is important, as we are now able to add fluctuations $u'_0(x)$ in a controlled way around the linearisation state to get our actual initial conditions
\begin{equation}
u(x,t=0) = \widetilde{u}(x)+u'_0(x),
\end{equation}
such that we can monitor the growth of the fluctuations in comparison to the linearisation state. 

\begin{figure}[!htbp]
  \centering
\includegraphics[width=0.48\textwidth,trim=0 0 0 0,clip]{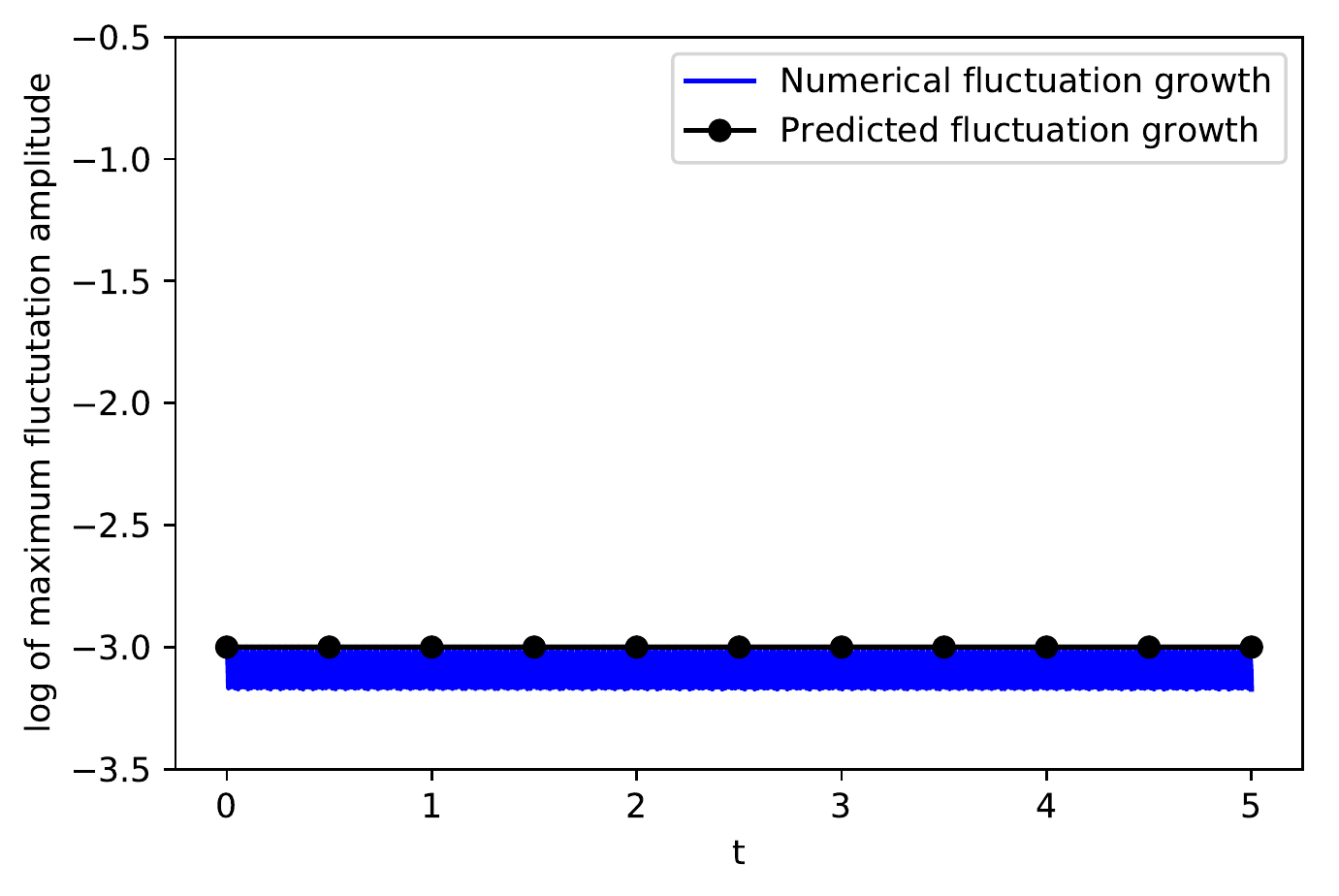} \includegraphics[width=0.48\textwidth,trim=0 0 0 0,clip]{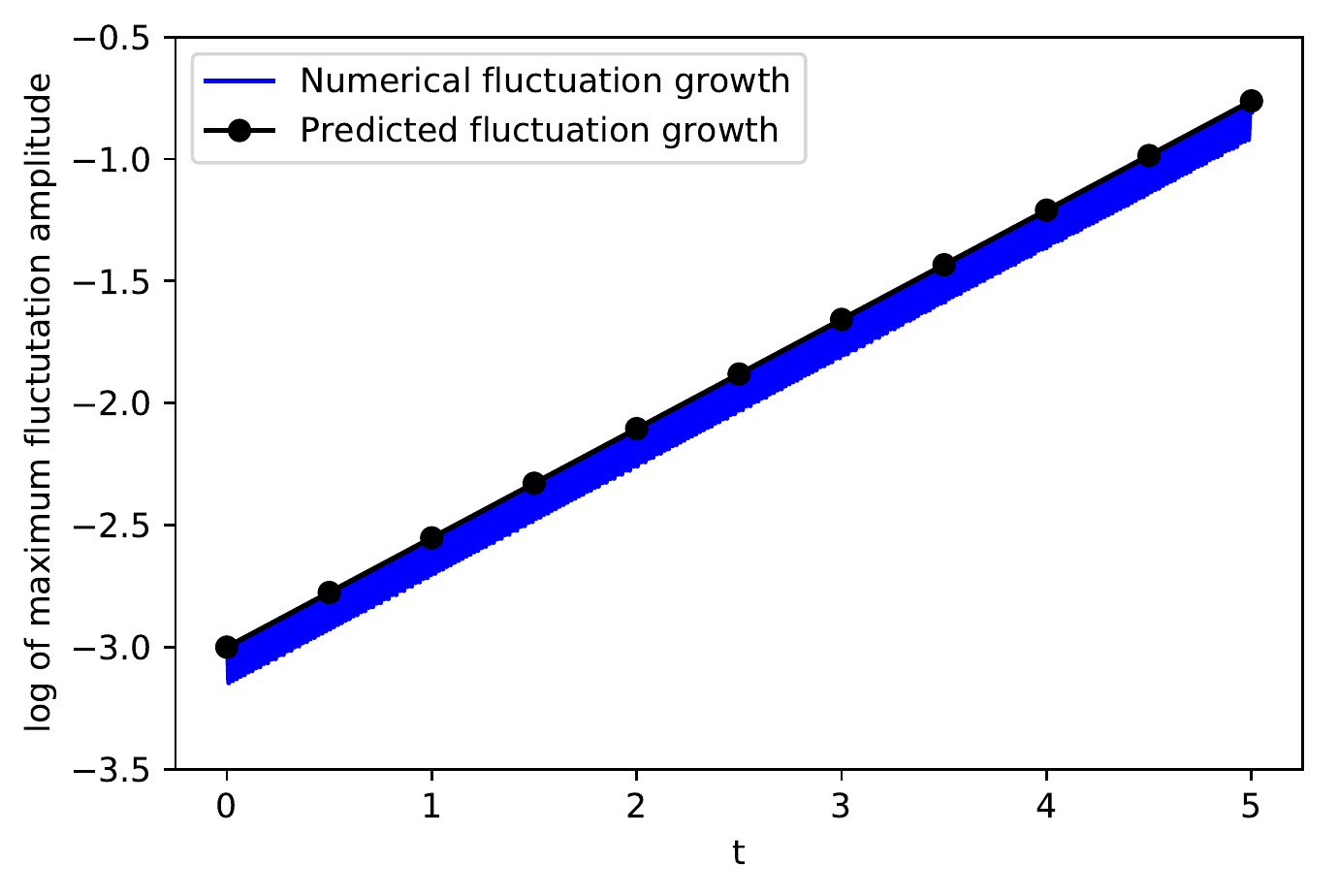} 
\caption{
Evolution of the maximum amplitude of the fluctuations of the inhomogeneous simulation. Left: Central approximation, volume term in divergence form ($\alpha=1$) and  central numerical flux for the surface term. Right: Entropy-conserving scheme with the skew-symmetric volume terms ($\alpha=2/3$) and the entropy-conserving flux \eqref{eq:ec}. \figinfo{Discretisation with $10$ elements and polynomial degree $N=3$, using baseflow from Figure~\ref{fig:linearisation_state_burgers} for initialisation and  $\rhs(\widetilde{u})$, and an added initial fluctuation being the scaled eigenmode shown in Figure~\ref{fig:critical_eigenmode_burgers} with the maximum initial amplitude $10^{-3}$. End time of the simulation is $T=5$. For comparison, we take the maximum real part of the corresponding spectra and use it to predict the exponential growth. The maximum real part is $8.8\times10^{-8}$ (left) and  $1.0307$ (right). Spectra shown in Figure~\ref{fig:spectrum_burgers_1}.  }
 }
\label{fig:fluctuations_burgers}     
\end{figure}

To bring the simulation results in context with the investigation of the spectra from the previous subsection, we choose the critical eigenmode corresponding to the eigenvalue with the largest real part shown in Figure~\ref{fig:critical_eigenmode_burgers} as our initial fluctuation $u'_0(x)$. Note that the eigenmode is scaled such that its maximum peak is $10^{-3}$, thus we have an initial perturbation with amplitude $10^{-3}$. Throughout the simulation, we monitor the maximum peak of the approximative solution in comparison to the linearization state $\widetilde{u}$ after every time step and plot the result in log-scale over time in Figure~\ref{fig:fluctuations_burgers}. To compare these numerical results to the previous analysis, we take the maximum real part of the spectra and use this to estimate the growth rate and to generate a prediction of the corresponding exponential growth, starting with the perturbation amplitude $10^{-3}$. The spectrum for the central approximation shows indeed no exponential growth, as the amplitude of the fluctuations stay at about the initial amplitude $10^{-3}$. In stark contrast, the entropy-conserving discretisation shows significant exponential growth, with a growth rate matching the prediction from the previous analysis of the spectrum.\medbreak

\begin{figure}[!htbp]
  \centering
\includegraphics[width=0.48\textwidth,trim=0 0 0 0,clip]{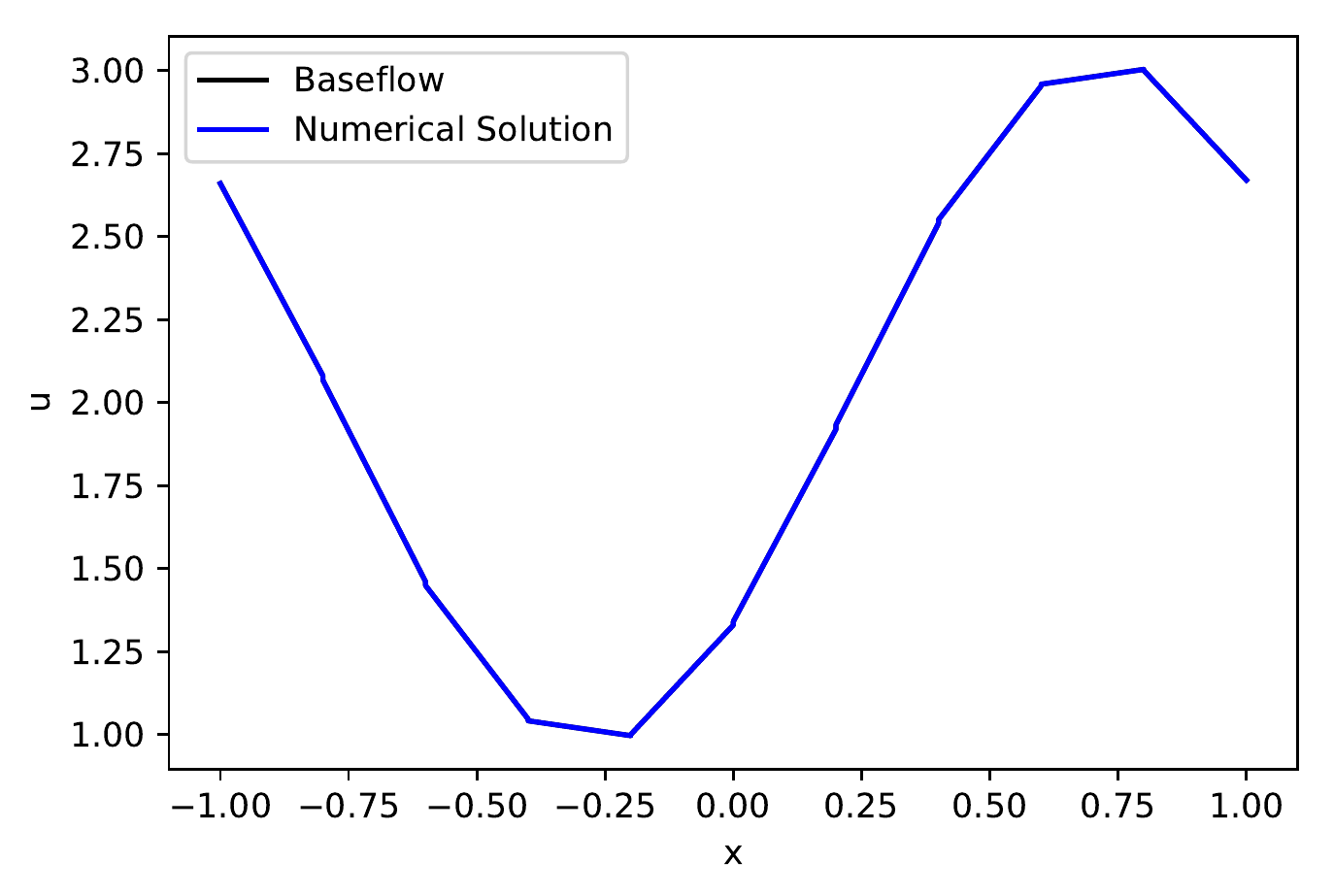} \includegraphics[width=0.48\textwidth,trim=0 0 0 0,clip]{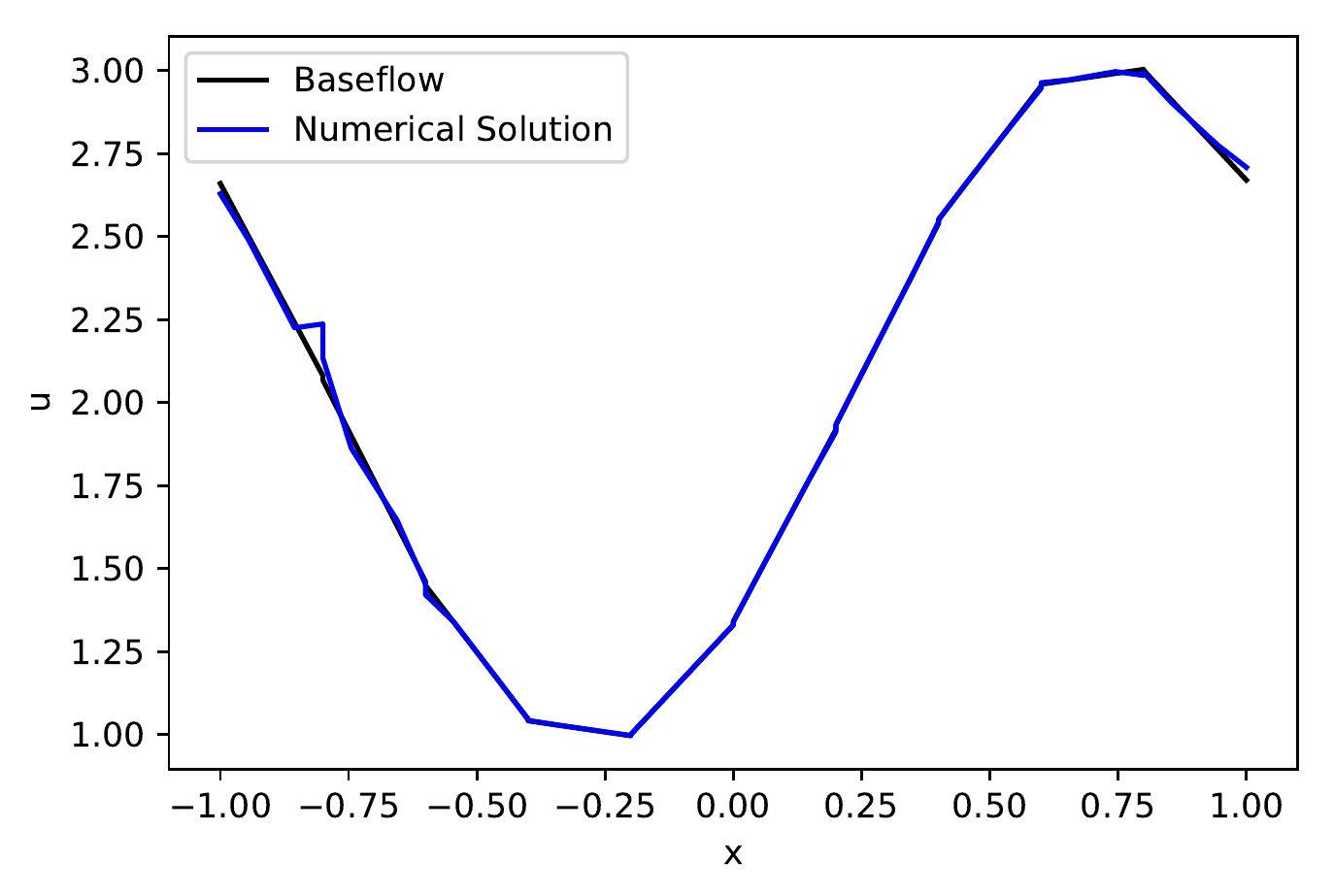} 
\caption{
Numerical solution of the inhomogeneous simulation at $T=5$ and comparison to baseflow. Left: Central approximation with divergence form volume term ($\alpha=1$) and the central numerical flux for the surface term. Right: Entropy-conserving scheme with the skew-symmetric volume terms ($\alpha=2/3$) and the entropy-conserving flux \eqref{eq:ec}. \figinfo{Discretisation with $10$ elements and polynomial degree $N=3$, using baseflow from Figure~\ref{fig:linearisation_state_burgers} for initialisation and  $\rhs(\widetilde{u})$, and an added initial fluctuation being the scaled eigenmode shown in Figure~\ref{fig:critical_eigenmode_burgers} with the maximum initial amplitude $10^{-3}$. Spectra shown in Figure~\ref{fig:spectrum_burgers_1}. } 
 }
\label{fig:solution_burgers}     
\end{figure}

For completeness, we plot the resulting numerical solutions in Figure~\ref{fig:solution_burgers}, where we can clearly observe how the faulty exponential growth of the small scale fluctuations start to get visible in case of the entropy-conserving scheme. The amplitude $10^{-3}$ of the initial fluctuations grows up to more than $0.1$, such that this can already be seen in the simulation results. It is interesting to note again, that the eigenmode is active in the part of the domain with negative slopes, hence the solution first starts to deviate from the baseflow in this part of the domain.

\begin{remark}\textit{
We did several investigations with different versions of the scheme, different configurations of volume and surface terms, and different initial fluctuation distributions. In every case, a faulty spectrum with a positive maximum real part leads to an exponential growth in the actual simulation and shows that the scheme is \textcy{locally energy} unstable. For discretisations where the maximum real part is zero (in relation to the accuracy of the approximative Jacobian), all simulations keep the initial condition with only negligible variations for all times and are \textcy{locally energy} stable. 
}
\end{remark}

\begin{figure}[!htbp]
  \centering
\includegraphics[width=0.48\textwidth,trim=0 0 0 0,clip]{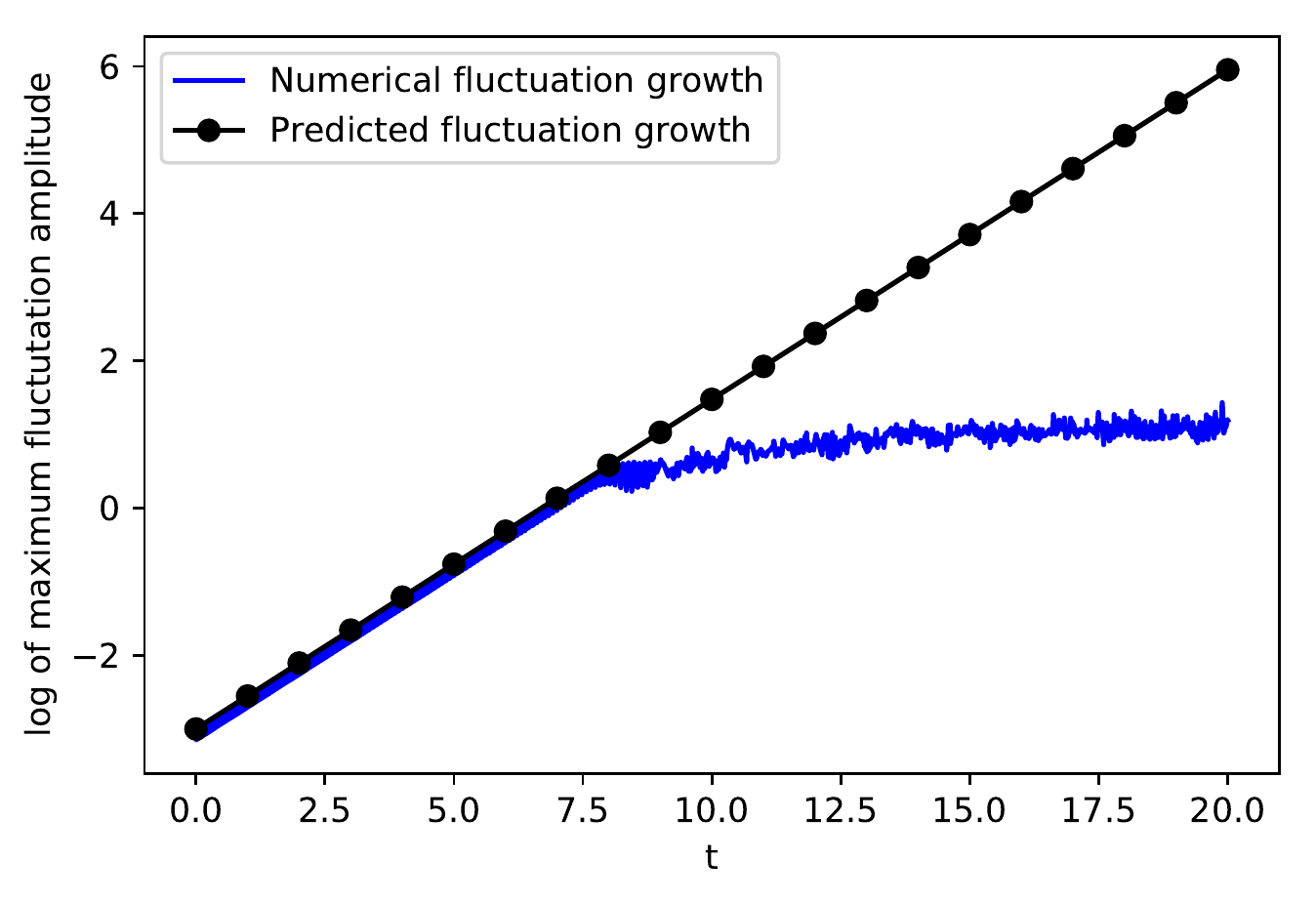} \includegraphics[width=0.48\textwidth,trim=0 0 0 0,clip]{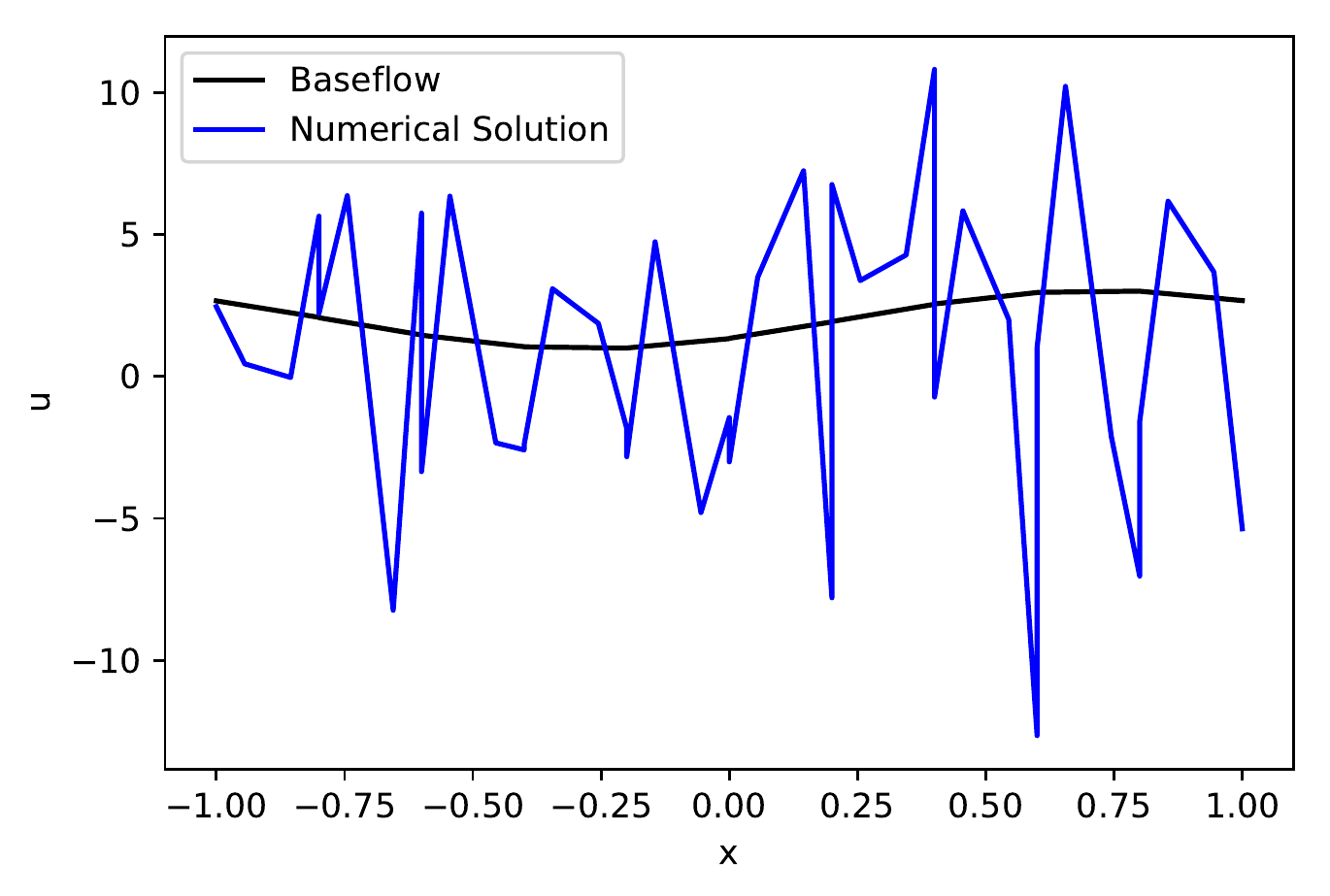} 
\caption{
Long time behaviour of the inhomogeneous simulation, using the entropy-conserving scheme with skew-symmetric volume terms ($\alpha=2/3$) and the entropy-conserving flux \eqref{eq:ec}. Left: Evolution of the maximum fluctuation amplitude over time, compared to the predicted exponential growth with the maximum real part $1.0307$. Right: Numerical solution at final time $T=20$ and comparison with the initial baseflow. 
\figinfo{Discretisation with $10$ elements and polynomial degree $N=3$, using baseflow from Figure~\ref{fig:linearisation_state_burgers} for initialisation and  $\rhs(\widetilde{u})$, and an added initial fluctuation being the scaled eigenmode shown in Figure~\ref{fig:critical_eigenmode_burgers} with the maximum initial amplitude $10^{-3}$.} 
 }
\label{fig:longtime_burgers}     
\end{figure}

Next, we consider a long term simulation of the entropy-conserving discretisation. The goal is to investigate the behaviour of the method when the amplitude of the fluctuations are so large, that non-linear effects cannot be neglected. We thus increase the end time from $T=5$ to $T=20$ and show the results in Figure~\ref{fig:longtime_burgers}. It is interesting to see that at about $t=7.5$, non-linear behaviour seem to emerge and start to dominate the scheme. It is important to note that the entropy-conserving scheme is neutrally stable in the $L^2$-norm, which is clearly visible in the left plot. The right plot shows the numerical solution at the final time, and reveals a wild collection of large scale oscillations that completely dominate the solution and overshadow the baseflow. It is important to recall, the the physical prediction is that the initial condition, the baseflow (black line), stays steady for all times!\medbreak

\begin{figure}[!htbp]
  \centering
\includegraphics[width=0.48\textwidth,trim=0 0 0 0,clip]{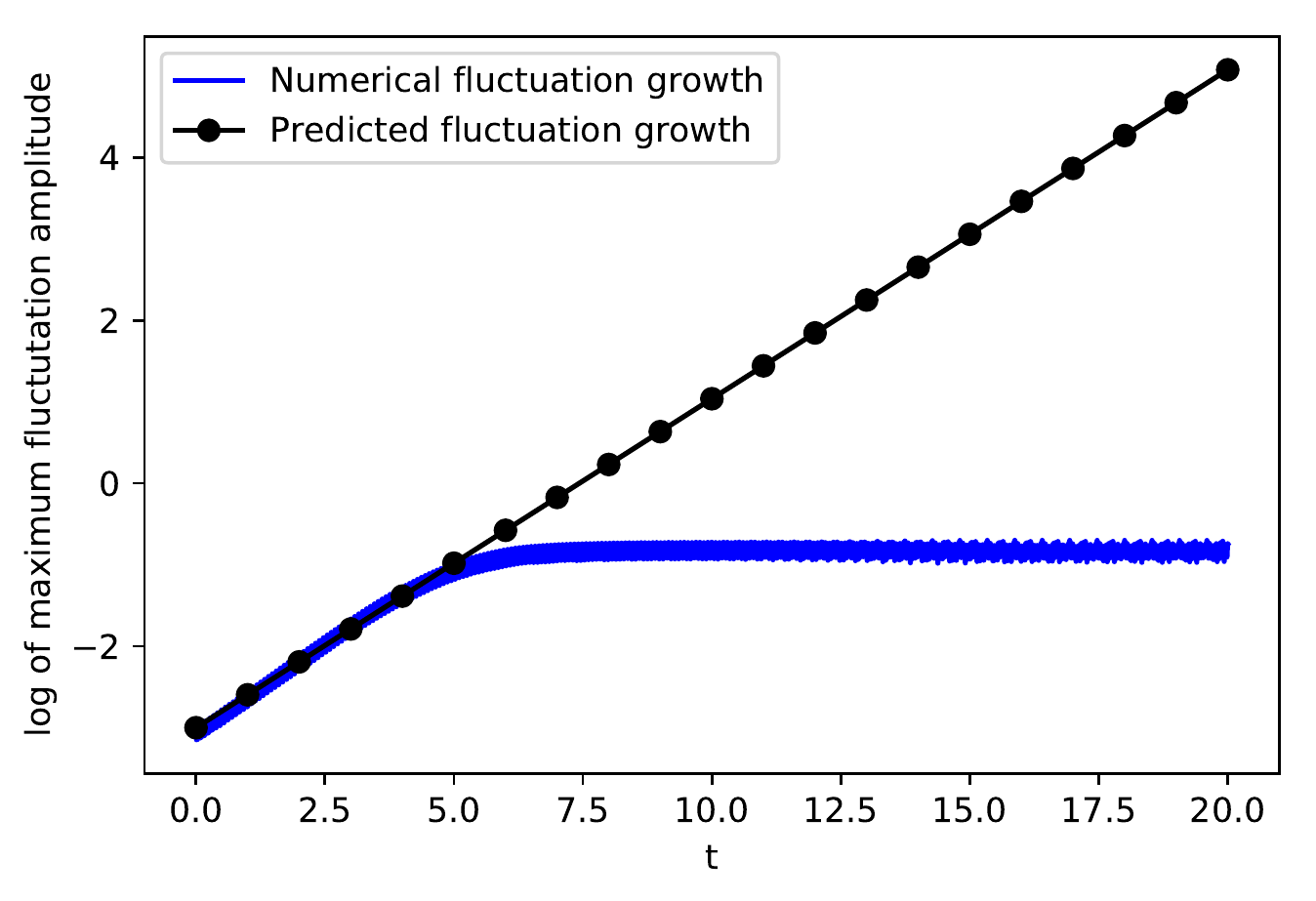} \includegraphics[width=0.48\textwidth,trim=0 0 0 0,clip]{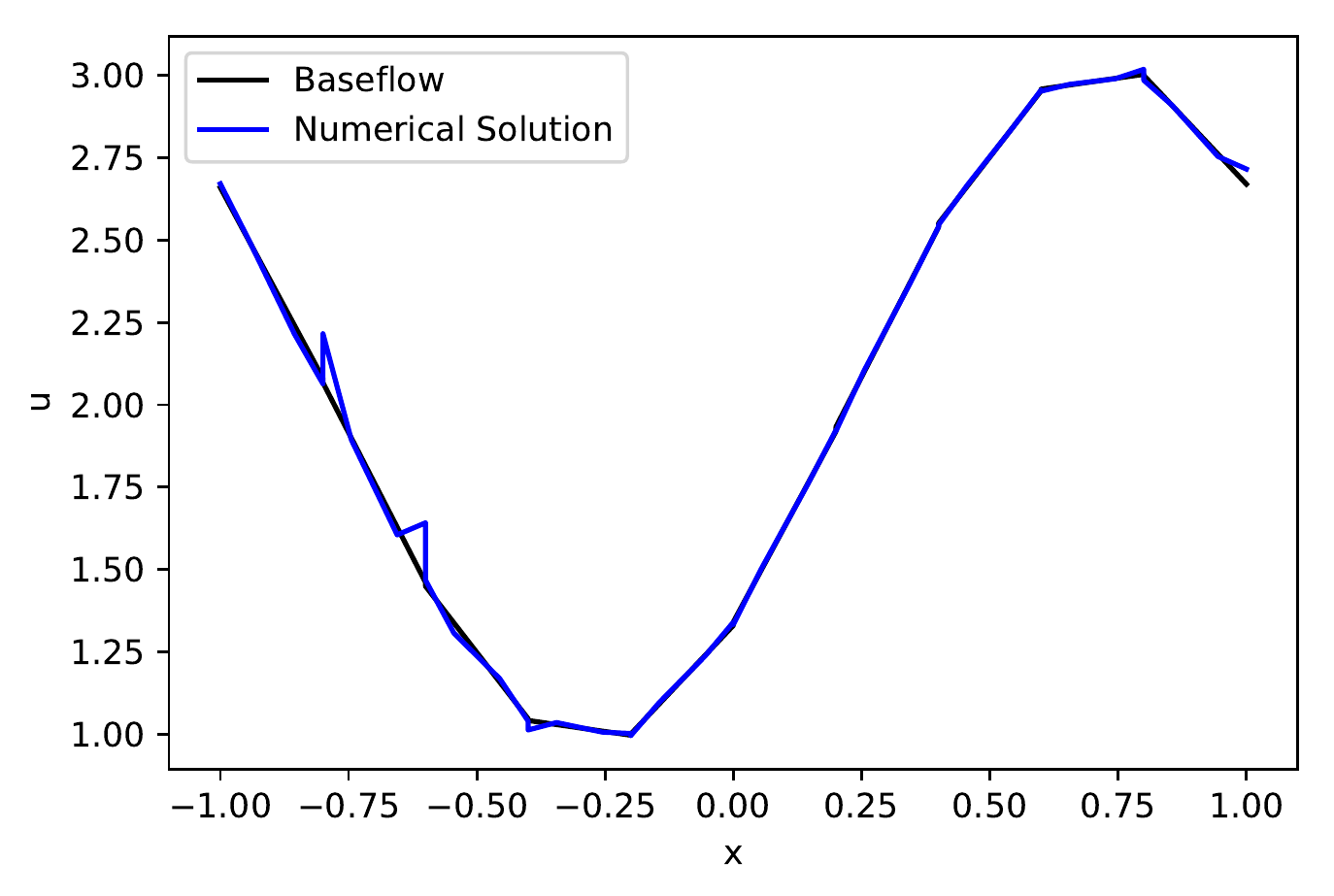} 
\caption{
Long time behaviour of the inhomogeneous simulation, using the modified entropy-dissipative scheme of Tadmor \eqref{burgers_stable2} and the volume modification of Carpenter et al. discussed in remark \ref{carpe}. Left: Evolution of the maximum fluctuation amplitude over time, compared to the predicted exponential growth with the maximum real part $0.929$. Right:  Numerical solution at final time $T=20$ and comparison with the initial baseflow. 
\figinfo{Discretisation with $10$ elements and polynomial degree $N=3$, using baseflow from Figure~\ref{fig:linearisation_state_burgers} for initialisation and  $\rhs(\widetilde{u})$, and an added initial fluctuation being the scaled eigenmode shown in Figure~\ref{fig:critical_eigenmode_burgers} with the maximum initial amplitude $10^{-3}$.} 
 }
\label{fig:longtime_burgers2}     
\end{figure}

We further show an interesting result for the dissipative modification of Tadmor \eqref{sec:fv_burgers} and the volume modification of Carpenter et al. discussed in Remark \ref{carpe}. As shown in the previous section, the scheme is not \textcy{locally energy} stable, despite the added dissipation. The left plot of Figure~\ref{fig:longtime_burgers2} shows this behaviour with a similar faulty exponential growth, until the non-linear effects kick in at about $t=5$ and stabilise the solution via its global $L^2$ estimate.
Due to the added dissipation, the global stability mechanism kicks in at an earlier time $t=5$ compared to the entropy-conserving scheme at time $t=7.5$ in Figure~\ref{fig:longtime_burgers}. Consequently, the maximum fluctuation amplitude is smaller in comparison as well and the numerical solution plotted in the left part of Figure~\ref{fig:longtime_burgers2} looks much 'nicer' compared to the result of the entropy-conserving scheme in Figure~\ref{fig:longtime_burgers}. Still, although it is better than the entropy-conserving result, this should not deflect the attention from the fact that it also displays faulty behaviour.\medbreak

\begin{figure}[!htbp]
  \centering
\includegraphics[width=0.48\textwidth,trim=0 0 0 0,clip]{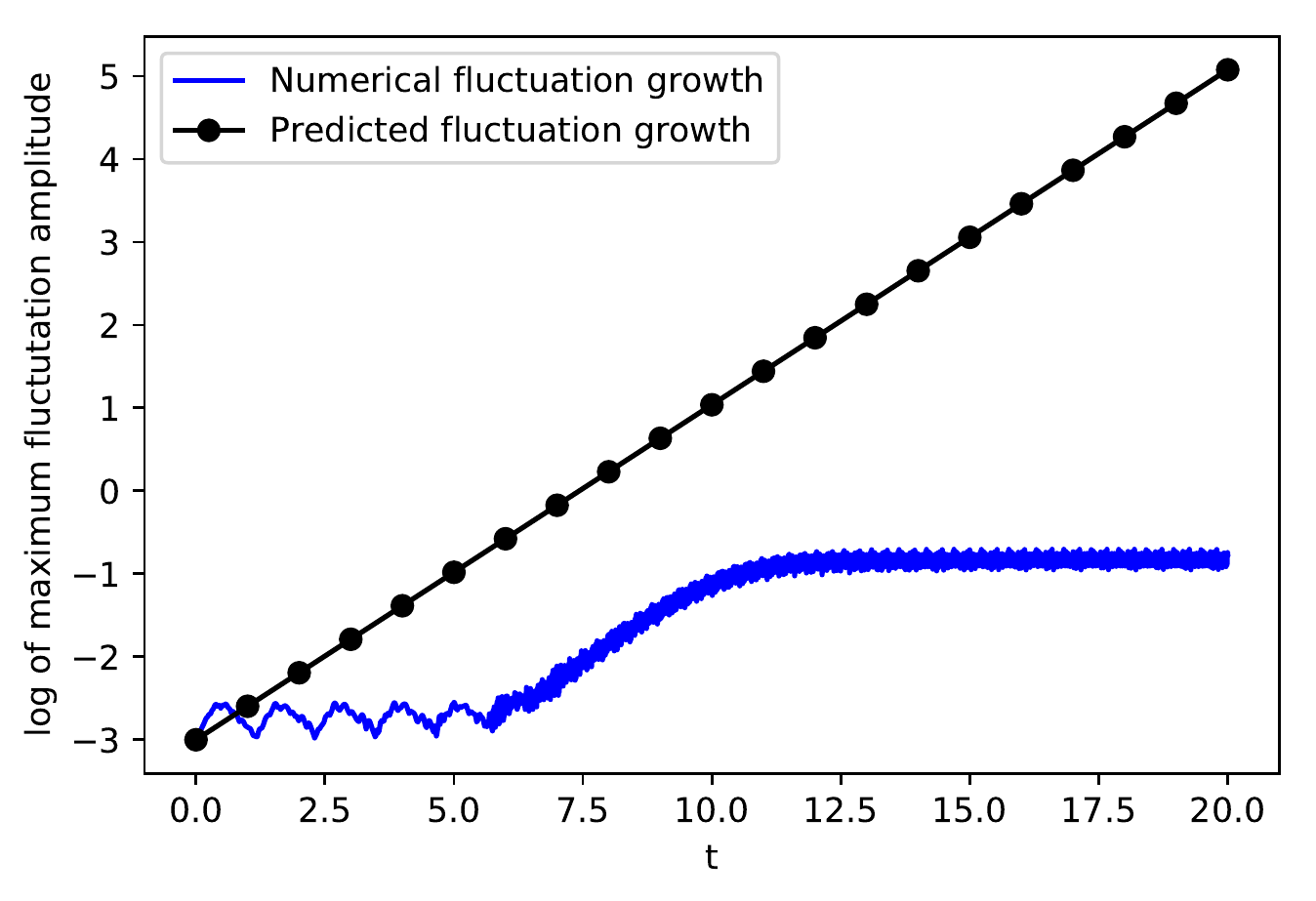} \includegraphics[width=0.48\textwidth,trim=0 0 0 0,clip]{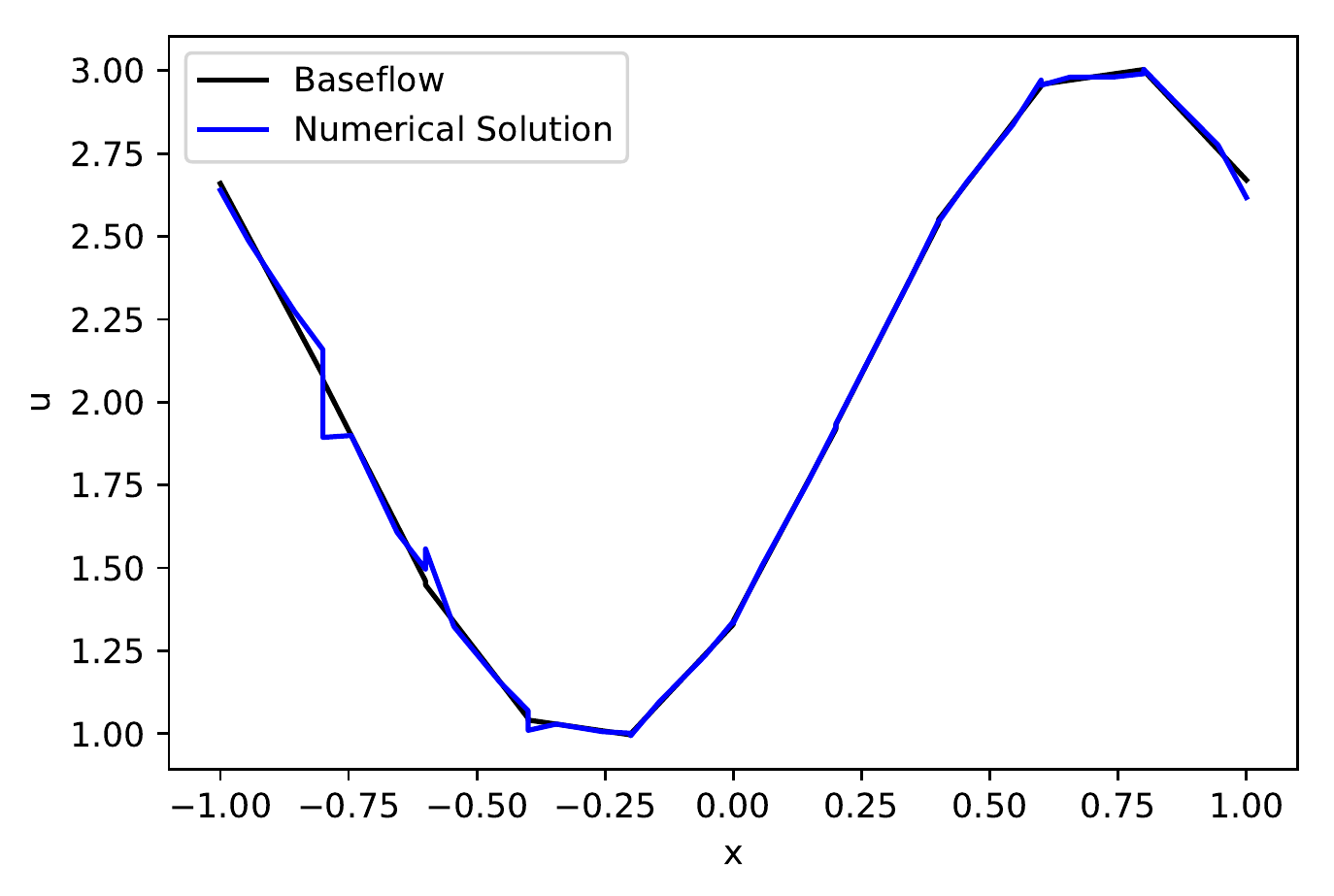} 
\caption{
Long time behaviour of the inhomogeneous simulation, with a smooth initial fluctuation, using the modified entropy-dissipative scheme of Tadmor \eqref{burgers_stable2} and the volume modification of Carpenter et al. discussed in remark \ref{carpe}. Left: Evolution of the maximum fluctuation amplitude over time, compared to the predicted exponential growth with the maximum real part $0.929$. Right:  Numerical solution at final time $T=20$ and comparison with the initial baseflow. 
\figinfo{Discretisation with $10$ elements and polynomial degree $N=3$, using baseflow from Figure~\ref{fig:linearisation_state_burgers} for initialisation and  $\rhs(\widetilde{u})$, and an added smooth initial fluctuation $u'_0(x) = 0.001\,\cos(\pi\,x)$.} 
 }
\label{fig:longtime_burgers3}     
\end{figure}

\begin{figure}[!htbp]
  \centering
\includegraphics[width=0.48\textwidth,trim=0 0 0 0,clip]{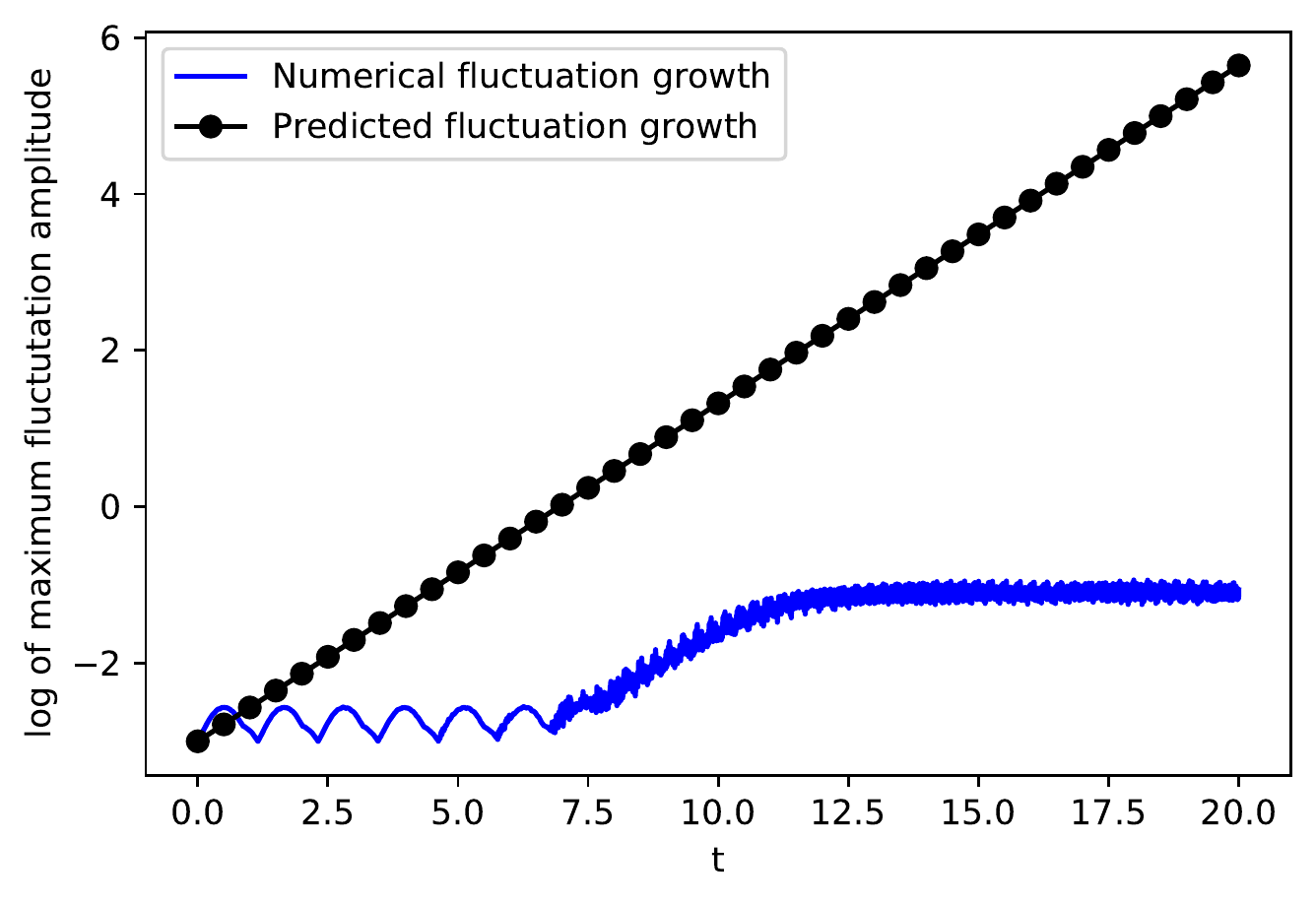}\includegraphics[width=0.48\textwidth,trim=0 0 0 0,clip]{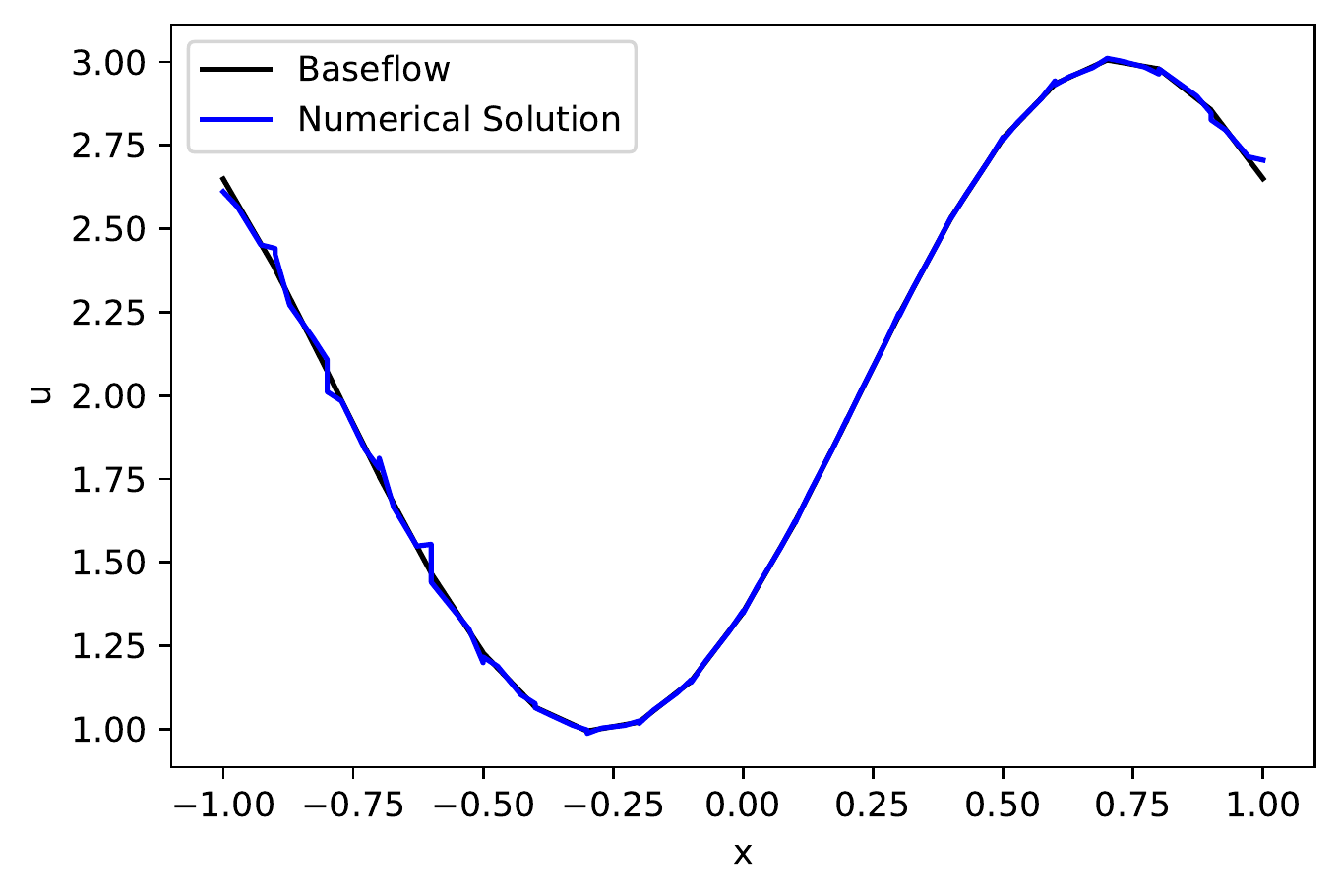} 
\caption{
Long time behaviour of the inhomogeneous simulation, with a smooth initial fluctuation and $20$ elements, using the modified entropy-dissipative scheme of Tadmor \eqref{burgers_stable2} and the volume modification of Carpenter et al. discussed in remark \ref{carpe}. Left: Evolution of the maximum fluctuation amplitude over time, compared to the predicted exponential growth with the maximum real part $0.929$. Right:  Numerical solution at final time $T=20$ and comparison with the initial baseflow. 
\figinfo{Discretisation with $20$ elements and polynomial degree $N=3$, using baseflow from Figure~\ref{fig:linearisation_state_burgers} for initialisation and  $\rhs(\widetilde{u})$, and an added smooth initial fluctuation $u'_0(x) = 0.001\,\cos(\pi\,x)$.}
 }
\label{fig:longtime_burgers4}     
\end{figure}

Finally, we document another interesting numerical result, where the initial fluctuations are not chosen as the particular eigenmode corresponding to the critical eigenvalue with large maximum real part. Instead, we choose a smooth distribution of the initial fluctuation $u'_0(x) = 0.001\,\cos(\pi\,x)$. This fluctuation is well resolved by the chosen discretisation. \textcy{The purpose of this particular setup is to investigate the behaviour of the scheme for smooth initial data under grid refinements, to further investigate if the schemes might be Lax stable.} We choose again the scheme with the entropy-dissipative surface modification of Tadmor and volume modification of Carpenter et al. for this assessment. It can be seen in the left plot of Figure~\ref{fig:longtime_burgers3} that until about $t=5$, the solution behaves very accurately and keeps the amplitude of the initial fluctuations in check at about $10^{-3}$. If one would choose the final time of the simulation smaller than $t=5$, one might get tricked into believing that the scheme works very well. However, the smooth distribution of the fluctuations also contain parts of the eigenmodes that correspond to eigenvalues with positive maximum real part. It comes as no surprise that these eigenmodes grow throughout the simulation and get amplified until they finally dominate the initial smooth fluctuations. Once, the eigenmode dominates, we can observe the exponential growth with the rate $0.9297$ until the non-linear behaviour takes over. 
If we increase the resolution from $10$ to $20$ elements, we discussed in the previous section, compare Figure \ref{fig:spectrum_burgers_111}, that the faulty eigenmodes get shifted to higher frequencies. Thus, we expect that for smooth (well-resolved) fluctuations the amplitude of the faulty eigenmodes decrease, as higher frequencies parts containing these modes should decay quickly. As can be seen in Figure~\ref{fig:longtime_burgers4}, we can observe that indeed the approximative solution stays stable for longer times, up to $t=7.5$ in comparison to $t\approx 6$. Consequently, the overall maximum growth of the fluctuation amplitude is smaller compared to the lower resolution before the non-linear stability estimate kicks in. However, we can also clearly observe that increasing the resolution only shifts the problem to later times and does not fully resolve it, as we still get artificial exponential growth of the fluctuations at later simulation times. \textcy{However, these results further strengthen the suggestion, that the schemes might be Lax stable and that they might be convergent for fine enough grids.}\medbreak

The results of the longterm simulations of the entropy-conserving or entropy-dissipative discretisations shown in Figure~\ref{fig:longtime_burgers},   \ref{fig:longtime_burgers2},  \ref{fig:longtime_burgers3} and \ref{fig:longtime_burgers4} are highly worrisome. They demonstrate that the discretisations might generate results that appear stable, but the simulations are nonphysical as they allow for an exponential growth of small scale fluctuations until the non-linear stability estimate kicks in and provides a global bound. \textcy{We note that} the physical behaviour of the linearised problem with the given positive baseflow is known and does not allow exponential growth. 

The worrisome aspect is that the scheme artificially generates solution structures, that are caused by the lack of \textcy{local energy} stability of the numerical scheme and not by actual physics. We have also shown that grid refinement might shift the issues to later simulation times \textcy{and that the schemes might converge for smooth solutions, again, hinting towards Lax stability.} At first, one might lean towards rating this behaviour as acceptable in practice, as one could just refine the grid until the simulation time can be reached safely. However, at second sight, this behaviour is at least as worrisome as we do not want to artificially increase the resolution of our discretisation, just to make it stable. This is certainly not a feasible strategy for three-dimensional turbulence. Furthermore, while we can easily assess the behaviour and judge the numerical results for the Burgers equation, simulation results for fluid dynamics equations in multiple spatial dimensions are much more complicated and harder to assess. It seems to be very difficult to decide, when the resolution in a complex application is enough to suppress the artificial exponential growth (assuming that it is even feasible to increase the resolution). In particular, one may then ask the critical question, if those artificial structures generated by the scheme might get transformed by the high-order scheme into something that resembles a meaningful fluid dynamics structure. One might even interpret the extended range of medium and small scales appearing as a positive aspect of the "low numerical diffusion of the high-order approximation". 

To further investigate this hypothesis, we have to consider the true flow equations, which is the topic of the next section.

\section{Local stability analysis for the compressible Euler equations}
\label{sec:euler}

In this section, we show that the issues of the entropy-conserving fluxes carry over to the compressible Euler equations. To simplify the analysis, we first consider the one-dimensional version 
\begin{equation}
\label{eq:euler}
u_t + f(u)_x = \begin{pmatrix}\rho\\\rho\,v\\\rho\,E\end{pmatrix}_t + \begin{pmatrix}\rho\,v\\\rho\,v^2+p\\\rho\,E\,v+p\,v\end{pmatrix}_x  = \begin{pmatrix}0\\0\\0\end{pmatrix},
\end{equation}
where $\rho$ is the density, $v$ the velocity, $E$ the specific total energy, $p=(\gamma-1)\,\rho\,(E-v^2/2)$ the pressure with the perfect gas assumption, and $\gamma=1.4$ the adiabatic coefficient. For our analysis, we consider the entropy pair
\begin{equation}
U(u) = -\frac{\rho\,s}{\gamma-1},\quad F(u) = -\frac{\rho\,v\,s}{\gamma-1}, \quad s(u) = \ln\left(\frac{p}{\rho^{\gamma}}\right)
\end{equation}
where $s$ is the thermodynamic entropy. The corresponding entropy variables are
\begin{equation}
w = \frac{\partial U}{\partial u} = \left( \frac{\gamma - s}{\gamma-1}-\frac{\rho}{p}\,\frac{v^2}{2}\;,\; \frac{\rho}{p}v\;,\;-\frac{\rho}{p}\right)^T.
\end{equation}

We consider first the finite volume discretisation \eqref{FV} of the compressible Euler equations \eqref{eq:euler}. There are several numerical flux functions available in literature, e.g., \cite{IsmailRoe2009,Chandrashekar2012,Ranocha2017} that give discrete entropy-conservation. We consider for instance the variant presented by Chandrashekar \cite{Chandrashekar2012} that is given by
\begin{equation}
\label{eq:ec2}
f^{\EC}_{i+1/2} = \begin{pmatrix}\{\rho\}^{\ln}\,\avg{v}\\\{\rho\}^{\ln}\,\avg{v}^2+\hat{p}\\\{\rho\}^{\ln}\,\avg{v}\,\hat{h}\end{pmatrix},
\end{equation}
where $\avg{.}$ is the arithmetic mean value, the logarithmic mean is 
\begin{equation}
\{\rho\}^{\ln} = \frac{\rho_{i+1} - \rho_i}{\ln(\rho_{i+1}) - \ln(\rho_i)},
\end{equation}
and 
\begin{equation}
\hat{p} = \frac{\avg{\rho}}{\avg{\rho/p}},\quad \hat{h} = \frac{1}{\{\rho/p\}^{\ln}(\gamma-1)}+\avg{v}^2-\frac{1}{2}\,\avg{v^2} +\frac{\hat{p}}{\{\rho\}^{\ln}}.
\end{equation}

We consider a particular simple test problem with periodic boundary conditions in the domain $[-1,1]$, where the initial condition is given by a density wave
\begin{equation}
\label{eq:entropywave}
\rho(x) = 1+0.98\,\sin(2\,\pi\,x),\quad v = 0.1,\quad p=20,
\end{equation}
with a minimum density of $\rho_{\min}=0.02$. Plugging this solution into the compressible Euler equations \eqref{eq:euler}, the problem reduces to a set of constant coefficient advection equations. Thus, the exact solution of this problem is directly given by the characteristics theory, i.e. the wave is transported without changing its shape or form. As this test problem reduces the compressible Euler equations to a simple advection problem, we know that again the finite volume scheme with the central numerical flux 
\begin{equation}
\label{eq:central}
f^{\CN}_{i+1/2} = \frac{1}{2}(f(u_i) + f(u_{i+1}))
\end{equation}
is the least diffusive scheme that is \textcy{energy} stable. Moreover, all three equations reduce to exactly the same advection equation when using the central flux which is well-known to produce accurate solutions. Thus, we are again able to use the central flux as the baseline scheme for comparisons with other numerical flux functions.

We turn to the analysis of the mass equation and compare the entropy-con\-serv\-ing flux \eqref{eq:ec2} to the central flux \eqref{eq:central} and get the relationship
\begin{equation}
(f_1)^{\EC}_{i+1/2} = (f_1)^{\CN}_{i+1/2} - \frac{1}{2}\,(R_1)^{\EC}_{i+1/2}\,(\rho_{i+1} - \rho_i),
\end{equation}
where the dissipation coefficient is 
\begin{equation}
 (R_1)^{\EC}_{i+1/2} = \frac{2(\avg{\rho v} -\{\rho\}^{\ln}\avg{v})}{\rho_{i+1}-\rho_i}=(\avg{\rho}-\{\rho\}^{\ln})\frac{2\avg{v}}{\rho_{i+1}-\rho_i}+\frac{(v_{i+1}-v_i)}{2},
\end{equation}
which shows that again anti-dissipation could occur for certain configurations, as $(R_1)^{\EC}_{i+1/2}$ is indefinite. In particular, if we consider the test problem \eqref{eq:entropywave} with the constant positive velocity,  and since $\{\rho\}^{\ln}\leq\avg{\rho}$, the dissipation coefficient $(R_1)^{\EC}_{i+1/2}$ is negative in the part of the domain where the density slope is negative, i.e. $\rho_x < 0$. 
\begin{remark}\textit{
Unlike the central scheme, the entropy-conservative flux does not reduce to the same scheme for the advection equation for all three equations. Hence, the numerical solutions of the three equations will not be the same but interact with each other. Consequently, the scheme does not only introduce negative artificial dissipation, but also contains a mechanism for generating perturbations.
}
\end{remark}

This short analysis for the simple test case already demonstrates that the issues of anti-dissipation carry over from the discussion on the Burgers equation to the compressible Euler equations. We note that the entropy-conserving numerical flux is again the main building block to extend the finite volume approach to high-order discretisations \cite{Fisher12,Gassner:2016ye,carpenter_esdg}. Thus, it is expected that the issues with anti-dissipation and lack of \textcy{local energy} stability most likely carry over to the high-order framework.\medbreak

In the following numerical results, we do not simulate the evolution of the fluctuations with respect to a baseflow by using a inhomogeneous test case. Instead, we perform a standard simulation with the non-linear homogeneous compressible Euler equations. We are interested in numerical investigations of the actual behaviour of the high-order framework. We choose again a high-order discontinuous Galerkin spectral element method with Legendre-Gauss-Lobatto nodes, which belongs to the class of diagonal-norm summation-by-parts methods. In contrast to the Burgers equation, no explicit split-form to obtain entropy-conservation is known up to this point. Instead, Fisher \cite{Fisher12} and Carpenter et al. \cite{carpenter_esdg,Carpenter_etal16} extended the high-order entropy-conserving reconstruction of LeFloch et al. \cite{lefloch_rhode_2000} to diagonal-norm SBP operators on finite domains, including e.g. the spectral collocation operator with Legendre-Gauss-Lobatto nodes. For a detailed description of the algorithms in multiple dimensions and additional assessment of their properties, we further refer to \cite{Gassner_BR1,Gassner:2016ye}.

\begin{figure}[!htbp]
  \centering
 \includegraphics[width=0.41\textwidth,trim=0 0 0 0,clip]{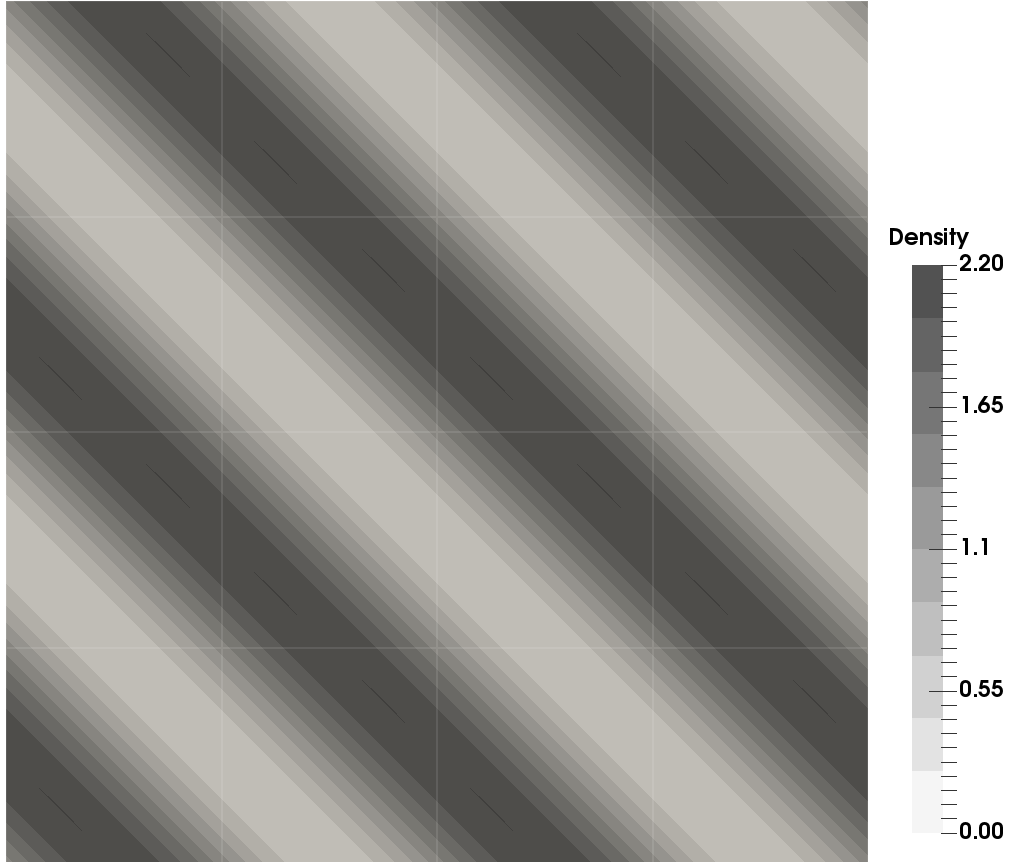}
\caption{Visualisation of the discrete density distribution of the initial condition \eqref{eq:entropywave2D} used for the numerical simulations, on a $4\times4$ grid with polynomial degree $N=5$. } 
\label{fig:initial_condition}     
\end{figure}

The numerical investigations are carried out with the open source three-di\-men\-sion\-al high-order simulation code Fluxo (github.com/project-fluxo) and for cross-evaluation with the open source Julia language based simulation framework Trixi.jl (github.com/trixi-framework/Trixi.jl). The assessment is done with $4\times 4$ elements in 2D in the domain $[-1,1]^2$ with polynomial degree $N=5$.
The initial condition \eqref{eq:entropywave} is extended to multiple spatial dimension and reads
\begin{equation}
\label{eq:entropywave2D}
\rho(x) = 1+A\,\sin(2\,\pi\,(x+y)),\quad A=0.98,\quad (v_1,v_2,v_3) = (0.1,0.2,0),\quad p=20.
\end{equation}
Figure~\ref{fig:initial_condition} shows the visualisation of the initialised density. We observe that this particular test problem is very well resolved. Furthermore, this initial distribution is also used as a linearisation state to compute the spectrum of the non-linear operators.\medbreak

The semi-discretisation is explicitly integrated in time with a low-storage 5-stage Runge-Kutta method of 4th order accuracy \cite{Kennedy:2003cr}. The time step is updated during the simulation with $\CFL=0.05$. The end time of the simulation is $T=5$, but we will add the final crash times for simulations that did not reach the end time. For the investigation of the spectra, we use the initial conditions as our linearisation state $\widetilde{u}(x,y,z) = u_0(x,y,z)$. The spectra are generated with Trixi.jl, which features automatic differentiation. Thus, the Jacobians are computed exactly, the in-build linear algebra package of Julia is then used to compute the eigenspectrum.

\begin{figure}[!htbp]
  \centering
\includegraphics[width=0.61\textwidth,trim=0 0 0 0,clip]{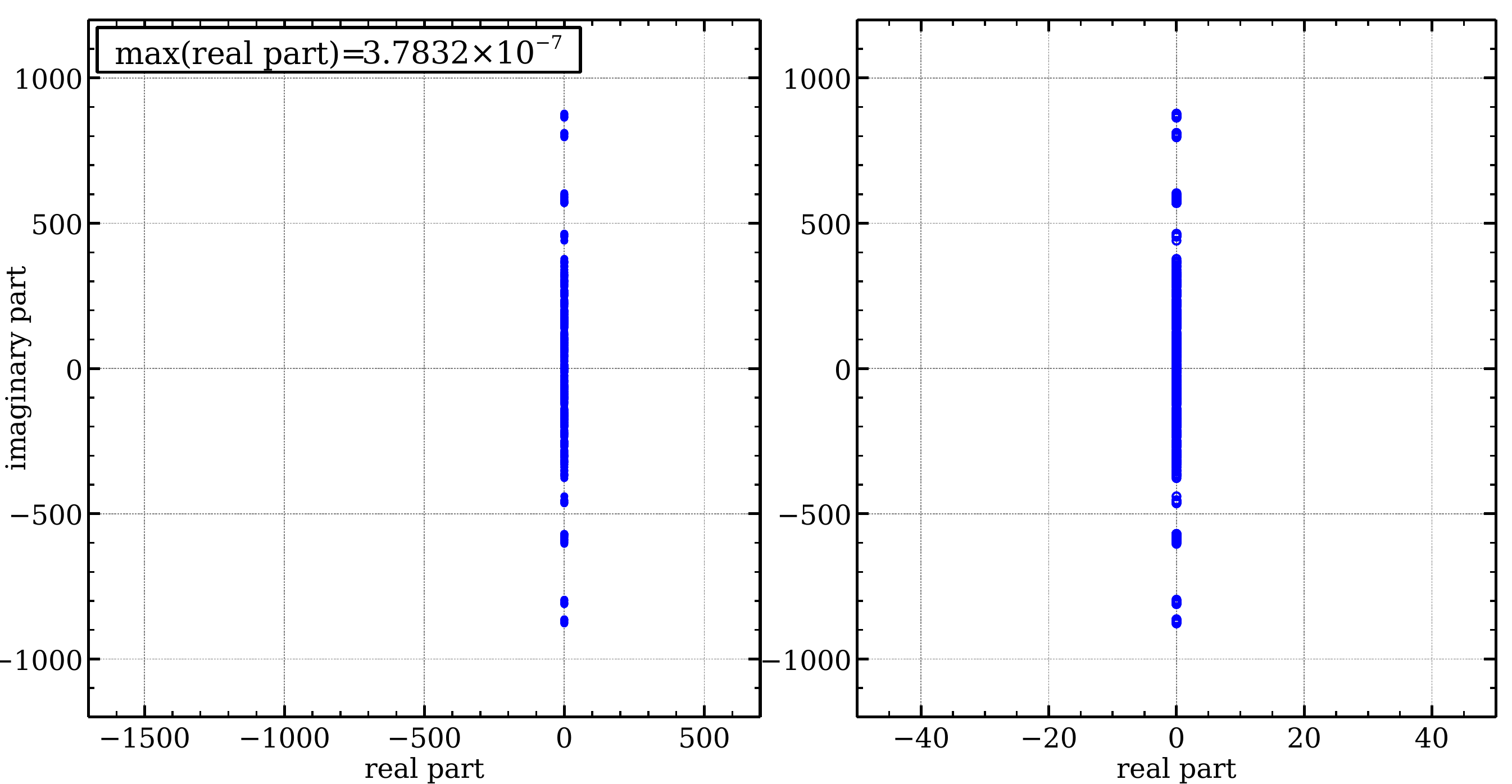} \hspace{0.02\textwidth}\includegraphics[width=0.36\textwidth,trim=0 0 0 0,clip]{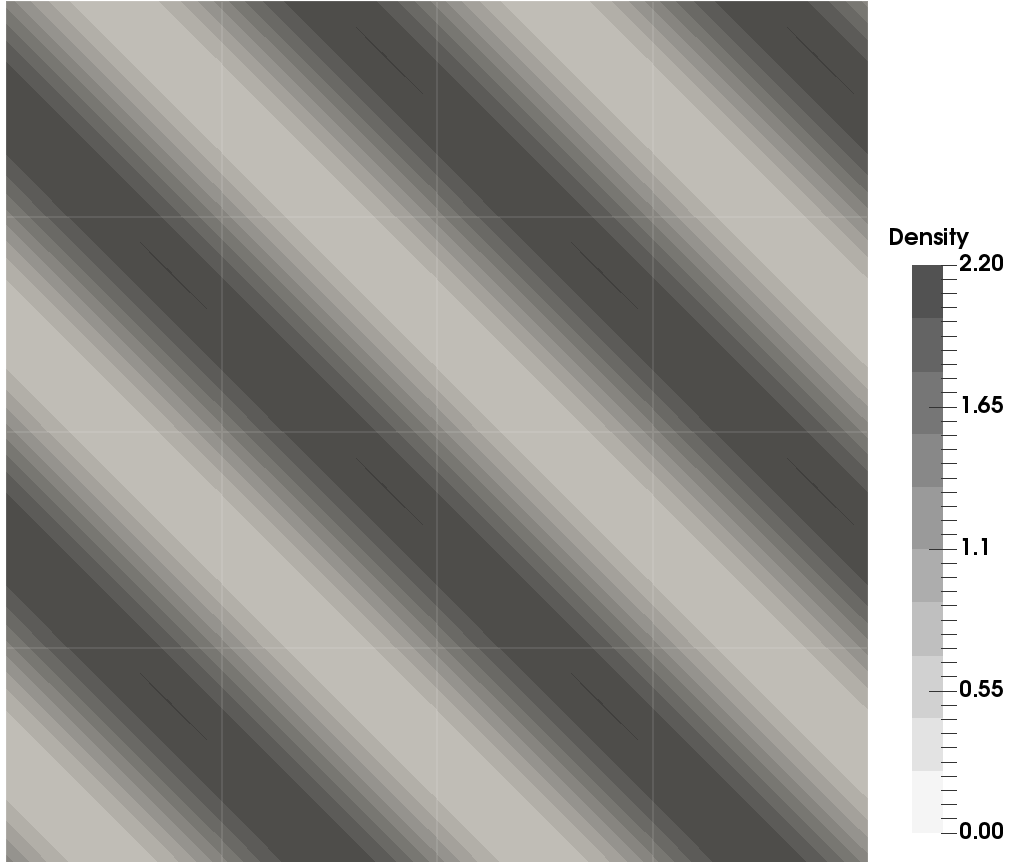} 
\caption{
Numerical results with the central approximation, volume term in divergence form and the central numerical flux function for the surface term.  Left: Spectrum of the spatial operator at $t=0$, maximum real part is $3.783\times10^{-7}$. Middle: Zoom in view of the spectrum. Right: Visualisation of the density distribution at the final simulation time $T=5$.
\figinfo{The grid is $4\times 4$ with polynomial degree $N=5$. The initial condition is shown in Figure~\ref{fig:initial_condition}.}
}
\label{fig:central_simulation}     
\end{figure}

As a first result, we verify in Figure~\ref{fig:central_simulation} the behaviour of the central approximation, where we use the volume terms in divergence form and the numerical surface flux as the central flux. As can be seen in the left plot, the spectrum is essentially imaginary, with the largest positive real part being $3.783\times10^{-7}$. To strengthen this assumption we performed a long time simulation with the central approximation up to $T=200$, which translates to over 1.6 million RK time steps without observing any spurious growth at all.


The right plot shows the corresponding simulation results at the final time $T=5$, and confirms that the solution is well resolved and properly advected by the dissipation-free purely central high-order approximation. \medbreak

\begin{figure}[!htbp]
  \centering
\includegraphics[width=0.61\textwidth,trim=0 0 0 0,clip]{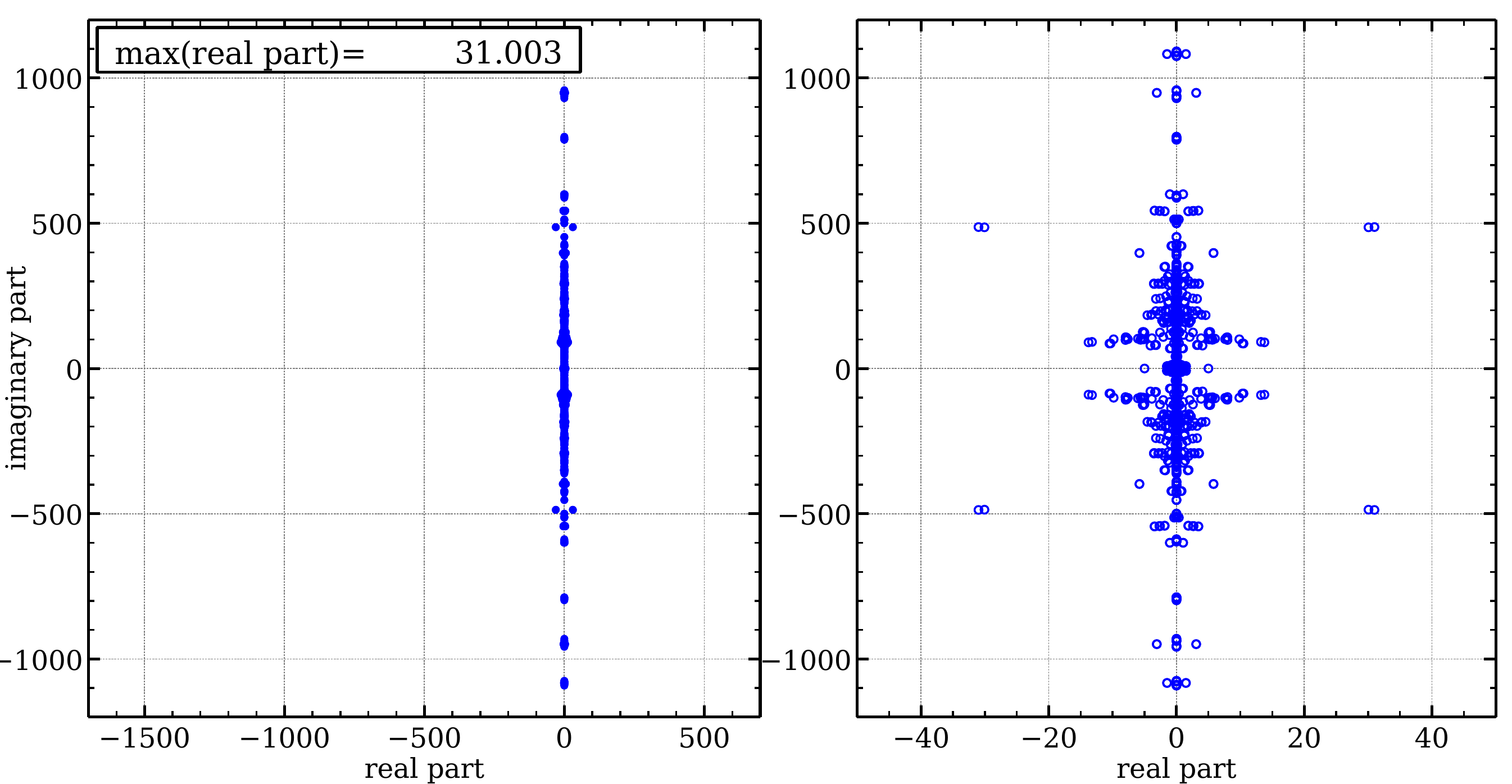} \hspace{0.02\textwidth}\includegraphics[width=0.36\textwidth,trim=0 0 0 0,clip]{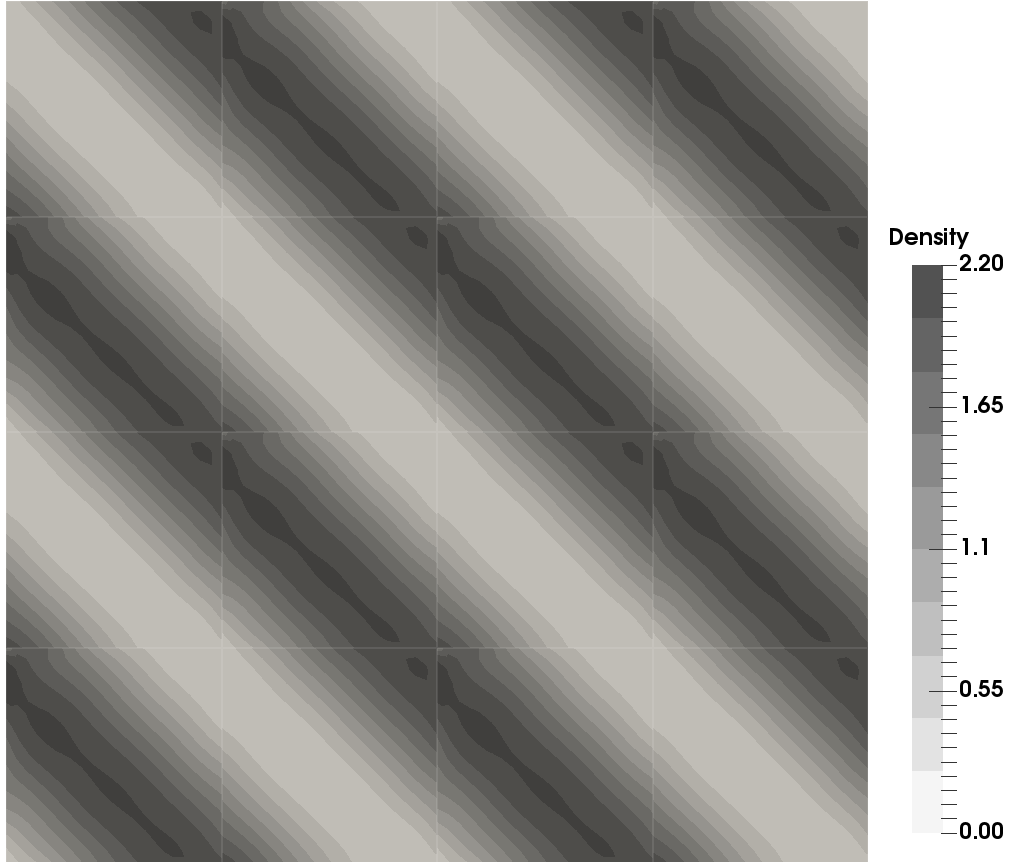} 
\caption{
Numerical results with the entropy-conserving approximation, using entropic flux-difference volume terms \cite{carpenter_esdg} with the entropy-conserving flux of Chandrashekar \cite{Chandrashekar2012}, which is also used for the surface term.
Left: Spectrum of the spatial operator at $t=0$, maximum real part is $31.003$. Middle: Zoom in view of the spectrum. Right: Visualisation of the density at time $t=0.55$, right before the code crashes because of negative density.
\figinfo{The grid is $4\times 4$ with polynomial degree $N=5$. The initial condition is shown in Figure~\ref{fig:initial_condition}.}
}
\label{fig:ec_simulation}     
\end{figure}

Next, we use the entropy-conserving high-order approximation based on the numerical flux function by Chandrashekar \cite{Chandrashekar2012}. (We only document results with this flux, but we also tested other entropy-conservative fluxes including the one by Ismail and Roe \cite{IsmailRoe2009} with similar outcome.) Figure~\ref{fig:ec_simulation} shows the results of this investigation. The left plot shows the spectrum of the operator and reveals that there are significant positive real parts with a maximum value of about $31.003$, which confirms our theoretical analysis from above. For strong negative density gradients, the scheme produces anti-dissipation which in our terminology implies a \textcy{local energy} instability. The impact on the simulation is drastic for this particular test problem, as the simulation crashes at time $t=0.5533$ before reaching the final time $T=5$. We note that for illustration, the test case is setup on purpose such that severe faulty behaviour in the discrete density evolution can cause negative density. However, we stress again that the simple central approximation is able to comfortably reach the final time $T=5$, and much longer end times without problems.  This is a simple well resolved test problem and it is thus even more worrisome that the entropy-conserving high-order discretisation fails. Furthermore, we also tested with other setups with less density variation and we could observe similar issues, it just requires longer simulation times for the instability to grow and crash the simulation\footnote{i.e. for a density variation $\rho=1+0.5\sin(2\pi(x+y))$,  the entropy-conserving scheme crashes at $T=26.7$.}. We emphasise that the behaviour is independent of the choice of CFL, as it is not a time integration issue but a spatial discretisation failure.\medbreak

\begin{figure}[!htbp]
  \centering
\includegraphics[width=0.61\textwidth,trim=0 0 0 0,clip]{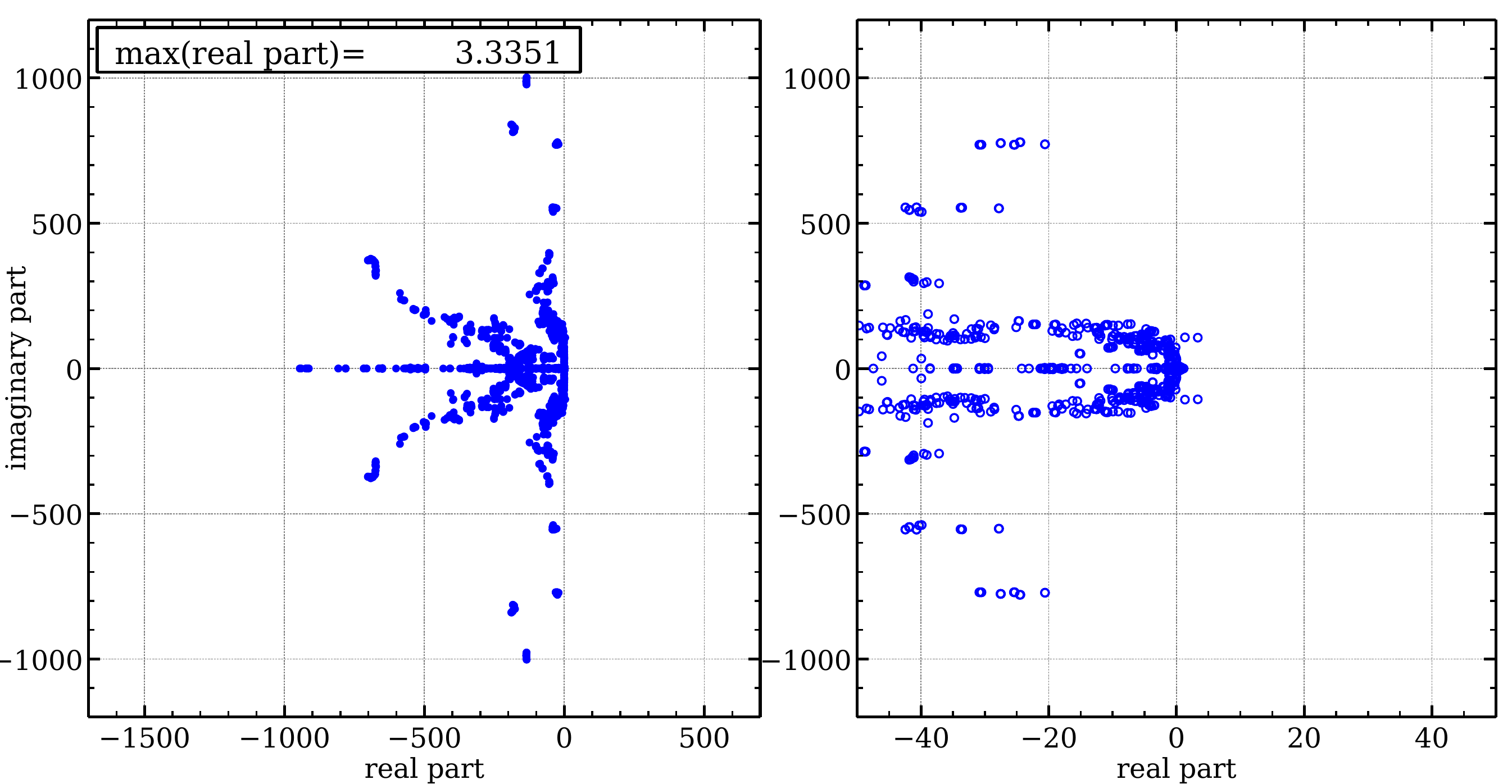}\hspace{0.02\textwidth}\includegraphics[width=0.36\textwidth,trim=0 0 0 0,clip]{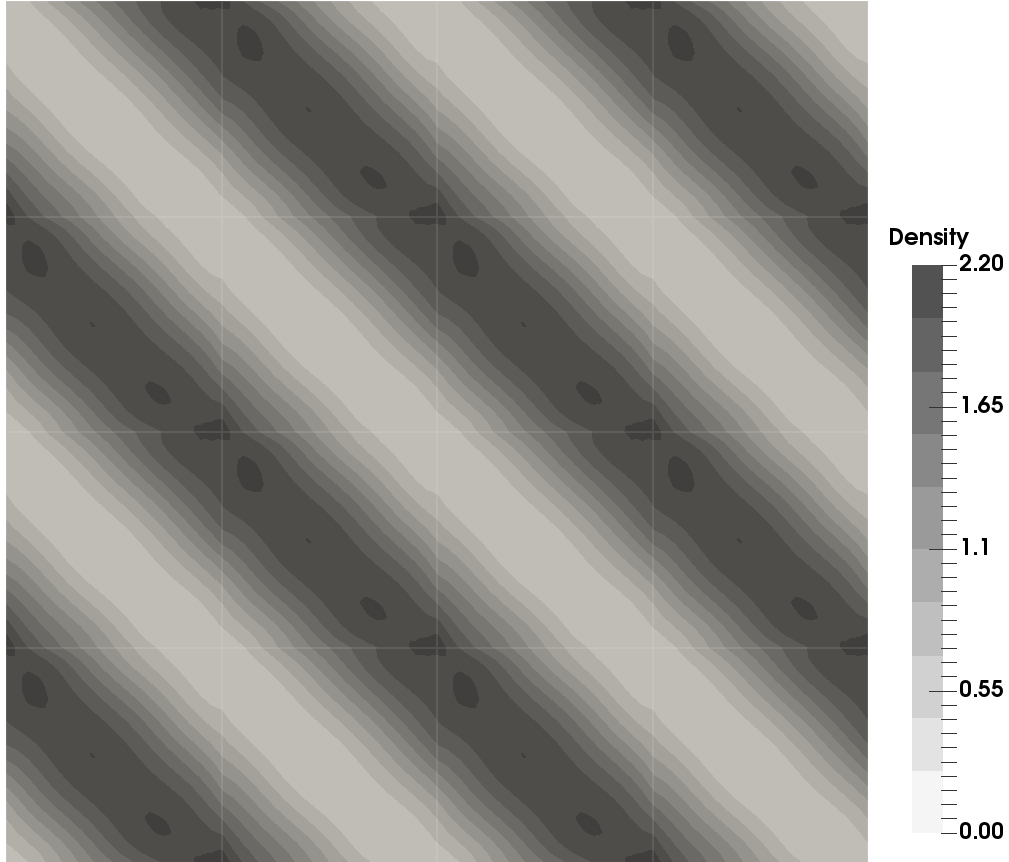} 
\caption{
Numerical results with the entropy-dissipative approximation, using entropic flux-difference volume terms \cite{carpenter_esdg} with the entropy-conserving flux of Chandrashekar \cite{Chandrashekar2012}, and guaranteed entropy-dissipative surface terms via the Rusanov flux \cite{toro2009}. Left: Spectrum of the spatial operator at $t=0$, maximum real part is $3.3351$. Middle: Zoom in view of the spectrum. Right: Visualisation of the density at time $t=0.65$, right before the code crashes because of negative density.
\figinfo{The grid is $4\times 4$ with polynomial degree $N=5$. The initial condition is shown in Figure~\ref{fig:initial_condition}.} 
}
\label{fig:es_simulation}     
\end{figure}

We next switch on the surface dissipation and use the Rusanov flux \cite{toro2009}. The high-order scheme is now formally entropy-dissipative. We note again that the volume terms are still in entropy-conserving form and it is not obvious that a surface dissipation term can fully control anti-dissipation generation through the volume. The results are plotted in Figure~\ref{fig:es_simulation}, which shows that the maximum real part is  $3.3351$ and hence the discretisation is \textcy{locally energy} unstable again. Indeed, the simulation crashes at an early time $t=0.6595$. \medbreak

\begin{figure}[!htbp]
  \centering
\includegraphics[width=0.61\textwidth,trim=0 0 0 0,clip]{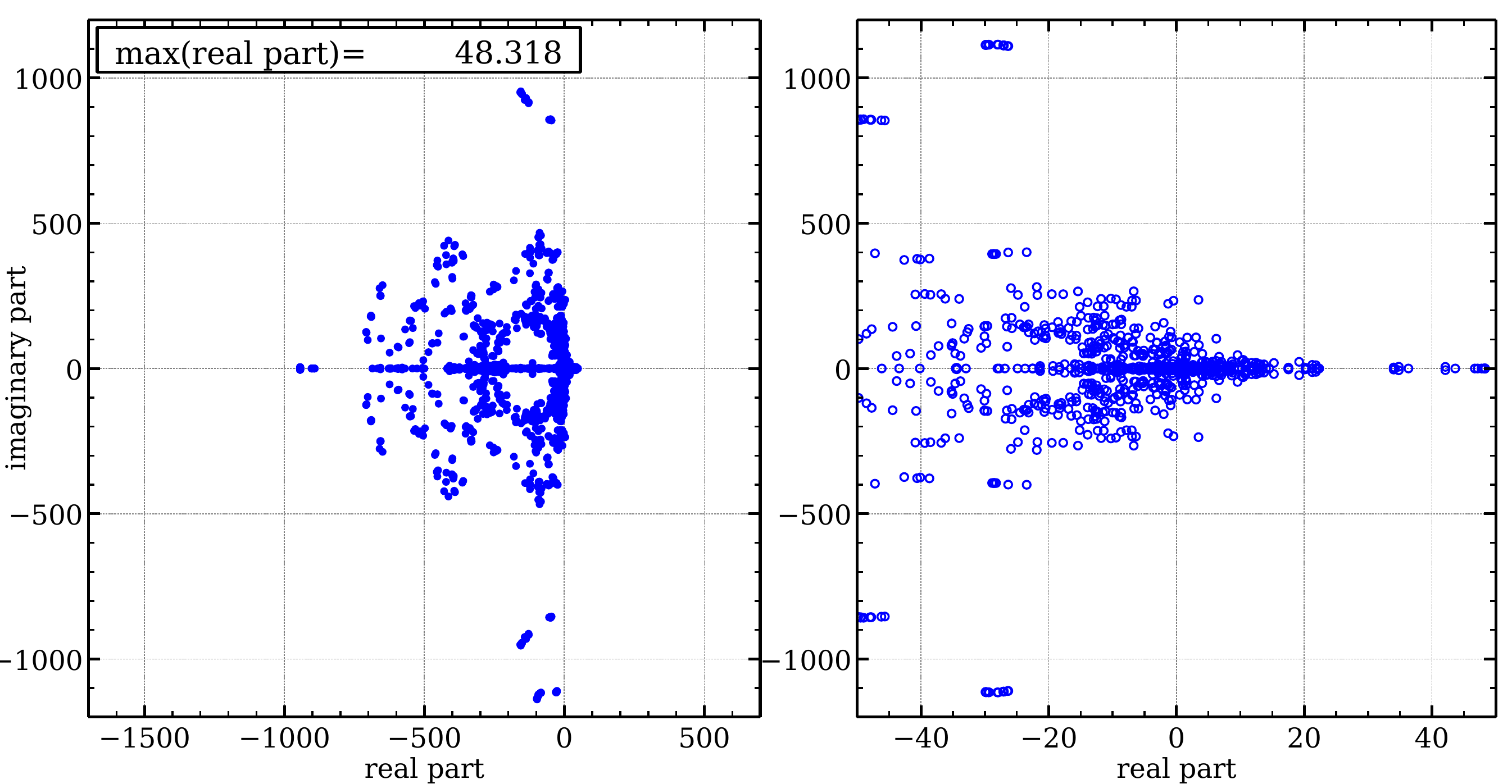}\hspace{0.02\textwidth}\includegraphics[width=0.36\textwidth,trim=0 0 0 0,clip]{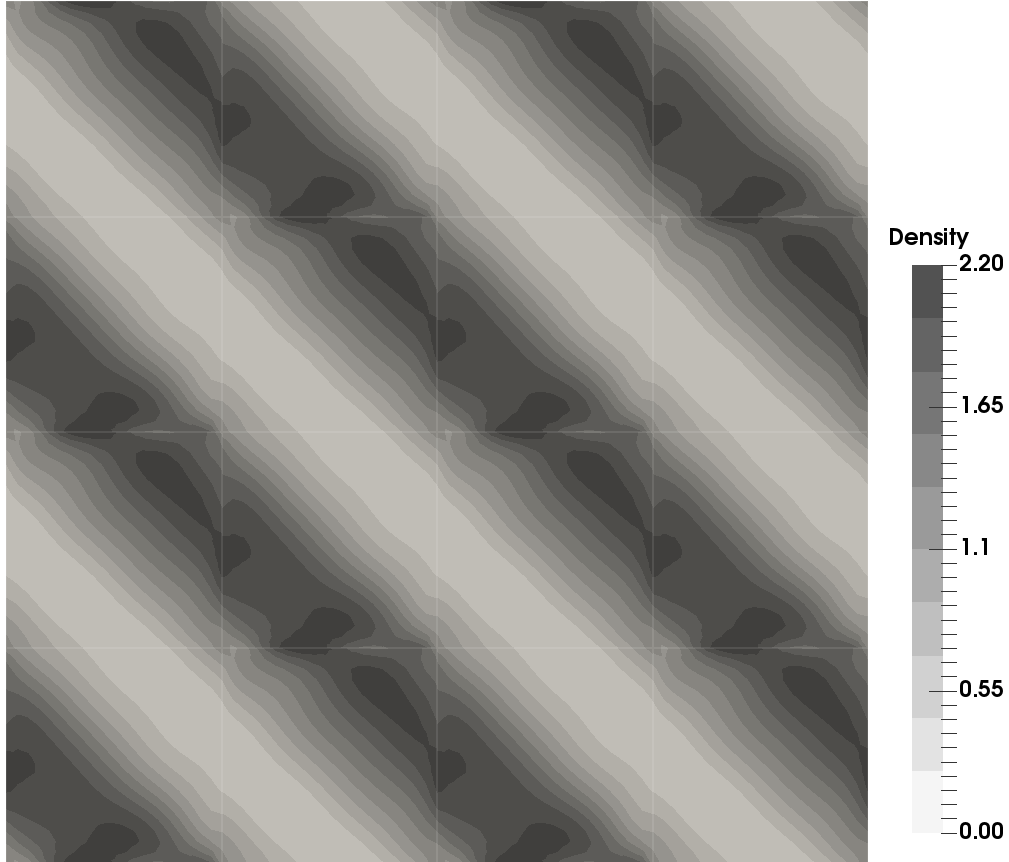} 
\caption{
Numerical results with the kinetic-energy-dissipative approximation, using split-form volume terms \cite{Gassner:2016ye} based on the Kennedy and Gruber splitting \cite{kennedy2008} and surface terms with the dissipative Rusanov flux \cite{toro2009}. Left: Spectrum of the spatial operator at $t=0$, maximum real part is about about $48.318$. Middle: Zoom in view of the spectrum. Right: Visualisation of the density at time $t=0.05$, right before the code crashes with negative density.
\figinfo{The grid is $4\times 4$ with polynomial degree $N=5$. The initial condition is shown in Figure~\ref{fig:initial_condition}. }
}
\label{fig:kg_simulation}     
\end{figure}

While we focus on the discretisations with entropic properties, it is interesting to note that there are many different split-forms for compressible fluid dynamics, with different properties such as e.g. kinetic-energy-preservation, see e.g. \cite{kennedy2008,ducros2000,Gassner:2016ye,Morinishi2010276,Coppola2019}. These split-forms are all different, but have one obvious common property: they are not the central approximation. Consequently, all these split-forms can be compared to the central approximation and have additional terms, that are either dissipative or maybe anti-dissipative. Hence, if for our considered baseflow the split-form does not reduce to the central approximation, the split-form has a similar issue with \textcy{local energy} stability as the entropy-conserving variant. As an example and numerical evidence, we consider the Kennedy and Gruber splitting \cite{kennedy2008} with the Rusanov numerical flux \cite{toro2009} for the surface terms. This gives a discretisation that is kinetic-energy-dissipative \cite{Gassner:2016ye}. The numerical results are shown in Figure~\ref{fig:kg_simulation} and clearly show faulty behaviour for this test problem. The maximum real part is about $48.318$, which again reveals that the discretisation is \textcy{locally energy} unstable. Again the impact on the actual simulation is strong, as the simulation crashes at an early time $t=0.0845$.

In appendix~\ref{sec:app}, we summarise additional simulation results obtained for the same density wave test case and the same schemes presented in this section. In particular, we compare the time evolution of the $L_2$ error in density for the different schemes, on a $4\times4$ and a $8\times8$ grid and also for two 'weaker' amplitudes $A$ of the initial density from \eqref{eq:entropywave2D}. The results underline that the undesired behaviour of the EC and KG split-form schemes is not limited to one particular setup of the density wave and \textcy{while being reduced (another hint that these schemes are Lax stable)} does not disappear \textcy{completely} on the finer grid or when decreasing the density gradients.

\section{Conclusion}
\label{sec:conclusion}

In this work, we demonstrate worrisome and \textcy{sometimes} downright faulty behaviour of the recently developed class of entropic high-order summation-by-parts schemes (including split-forms). To articulate the issue, we introduced the notion of \textcy{local energy} stability, i.e. the stability of fluctuations around a linearisation state for the non-linear problem. The core of the issue is that discretisation of the non-linear PDE and linearisation of the resulting scheme is not the same as discretising the linearised equations. Hence, non-linear stability does not automatically give \textcy{local energy} stability.

\textcy{To study the stability properties of the numerical schemes, we deliberately choose test cases, where the continuous problem has no exponential growth.} For these special test cases, we could use the central approximation based on the simple central flux as a baseline discretisation to obtain the least diffusive scheme that is \textcy{locally energy} stable. We compare entropy-conserving fluxes to the central flux and observe that the difference can either be dissipative or anti-dissipative. Artificial anti-dissipation is problematic, as it tries to artificially "push" the solution away from the linearisation state, i.e. it generates artificial exponential growth of small scale fluctuations. The theoretical findings have been verified in numerical tests, where the corresponding spectra of the spatial operator show significant positive real parts \textcy{for continuous problems with no exponential growth} and the simulations demonstrate actual artificial exponential growth of fluctuations. We observed that adding a dissipative mechanism through the numerical fluxes at the interfaces does not necessarily guarantee \textcy{local energy} stability either, and neither does adding dissipation to the volume terms. Both these "fixes" at best only mask the inherent issue. 

We further show that the faulty behaviour carries over to the compressible Euler equations with a potential negative impact on the robustness of these schemes. \textcy{Again, we choose a smooth exact solution of the compressible Euler equations that has no exponential growth}. While the simple central approximation, without added dissipation, happily runs this test case, the entropy-conserving and the entropy-dissipative schemes crash at early times due to nonphysical exponential growth of density fluctuations. However, the central approximation lacks a non-linear stability estimate and thus cannot be directly applied for strongly non-linear problems, such as e.g. under-resolved turbulence. As a final remark, we demonstrate that this faulty behaviour also appears in other split-forms used in the fluid dynamics community, including kinetic-energy-preserving variants such as the Kennedy and Gruber splitting, with similar worrisome numerical results as shown in our numerical investigations. As a matter of fact, the issue carries over to all split-forms that do not reduce to a central discretisation for the particular test cases considered in this work. However, the central discretisation is not preferable either as it does not have a non-linear global bound.

\textcy{Our theoretical derivations and numerical investigations helped us to rule out energy stability. However, the derivations also show that the anti-dissipation is proportional to the mesh size $h$ and thus the schemes might be Lax stable. We have thus investigated the impact of grid refinement on the spectra and the exponential growth. Grid refinement steps do not change the magnitude of the faulty eigenvalues, but do shift them to higher frequencies. Furthermore, the faulty exponential growths of smooth fluctuations get pushed to later times when refining the grid. These investigations indicate that (i) for a given simulation end time, it is possible to choose a grid resolution where the erroneous exponential growth is not dominating the results yet; (ii) the entropic (including split-forms) schemes might be Lax stable and hence the scheme might be converging for smooth solutions as the mesh size $h\rightarrow 0$. For practical simulations with finite grid size $h>0$, however, it is clearly not a desirable behaviour of any high-fidelity scheme, that the grid resolution has to be adjusted solely to push the faulty behaviour to later times such that they do not affect the results.}

It is clear that the erroneous behaviour of the high-order entropic split-form schemes needs further investigations and a proper fix. We tested several mechanisms to introduce dissipation through the volume, but it seems that it is not straightforward to find a strategy that does provide enough dissipation for \textcy{local energy} stability without overwhelming the high-order accuracy with excessive dissipation. The goal is to construct a scheme that is \textcy{locally energy} stable, has a non-linear entropy bound and is still high-order accurate. So, far, we are only able to get two of these properties when using the high-order SBP split-form approach. Thus, further research is needed and this work is a first documentation of the issues. Without a fix, results obtained with these schemes should be interpreted carefully. Due to the \textcy{locally energy} unstable behaviour, nonphysical structures might get introduced and subsequently evolved by the fluid equation 
to appear as 'meaningful' vortices and turbulent-like features, while being purely a numerical artefact. 

\begin{acknowledgements}
The authors thank Dr. David Kopriva, Dr. Hendrik Ranocha and Dr. Juan Manzanero for the discussion on linear stability of split-forms, linearisation and non-linear stability. We further thank Dr. Hendrik Ranocha for the extension of Trixi.jl with automatic differentiation capabilities that we use for the 2D compressible Euler results.
\end{acknowledgements}

\section*{Declaration}
\subsection*{Funding}
Gregor Gassner was supported by the European Research Council (ERC) under the European Union's Eights Framework Program Horizon 2020 with the research project Extreme, ERC Grant Agreement No. 714487.
 \subsection*{Conflict of interest}
 The authors declare that they have no conflict of interest.
\subsection*{Availability of data and material}
Not applicable
\subsection*{Code availability}
For the numerical results, prototype simulation software was developed and is available for everyone upon request from the authors if not available as open source software.

\bibliographystyle{spmpsci}      
\bibliography{ggBib}   


\appendix
\normalsize
\section{Additional simulation results for the 2D density wave\label{sec:app}}

In this section, we summarise additional simulation results obtained for the nonlinear compressible Euler equations without source terms and the \emph{same} two-dimensional density wave \eqref{eq:entropywave2D} of Section~\ref{sec:euler}, where the density is initialised as a wave with an amplitude $A=0.98$ on a $4\times 4$ grid and a polynomial degree $N=5$. In particular, since the exact solution is known, we compare the time evolution of the $L_2$ error in density up to $T=200$ for all discretisation choices of volume and surface fluxes presented in Section~\ref{sec:euler}, namely the purely central scheme ('central') and the central scheme with Rusanov surface flux ('central+Rus.'), the entropy and kinetic energy preserving flux by Chandrashekar ('EC'), either with the same flux at the surface or with Rusanov ('EC+Rus.'), and also the split-form Kennedy and Gruber  flux with Rusanov surface flux ('KG+Rus.').

The evolution of the errors for all schemes is shown in Fig.~\ref{fig:euler_densitywave_L2_A098}, together with the results using a finer grid of $8\times 8$.
Additionally, the same analysis is repeated for two 'weaker' density amplitudes $A=0.75$, shown in Fig.~\ref{fig:euler_densitywave_L2_A075} and $A=0.50$, shown in Fig.~\ref{fig:euler_densitywave_L2_A050}.

The results clearly show that for
the central scheme with a central or a Rusanov surface flux, the error initially changes but remains on the same level until $T=200$. The EC scheme shows a large initial increase of the error and crashes, due to negative density, for all amplitudes $A$; the larger the amplitude, the earlier it crashes. The EC scheme with Rusanov surface flux suffers from the same initial error  increase as the EC scheme, but can recover the wave with a much larger error, except in the case of $A=0.98$. This is in perfect agreement with the finding that the anti-dissipation of the EC scheme is related to the density gradient, and the Rusanov dissipation acts at later times to counteract the growth of the solution. Finally, the KG scheme with Rusanov flux has the same issues as the EC scheme and either crashes or the errors become very large.

Regarding the results on the finer $8\times8$ grid, we observe the same behaviour, but on a lower error level. The EC and KG schemes show again a strong initial error increase and either crash or introduce significantly higher error levels.

For completeness, in Fig.~\ref{fig:euler_density_t200}, we plot the density distribution  at  the end of the simulation ($t=200$), discarding crashed simulations, for the central, EC and KG scheme with Rusanov surface flux, and for both grids and the three choices of the amplitude $A$. The large error of the KG scheme is clearly visible, for the EC scheme, the wave still looks intact, but as reported, has a much larger error.

\newcommand\figwidth{0.45\textwidth}

\begin{figure}[!htbp]
\centering
\includegraphics[width=\figwidth]{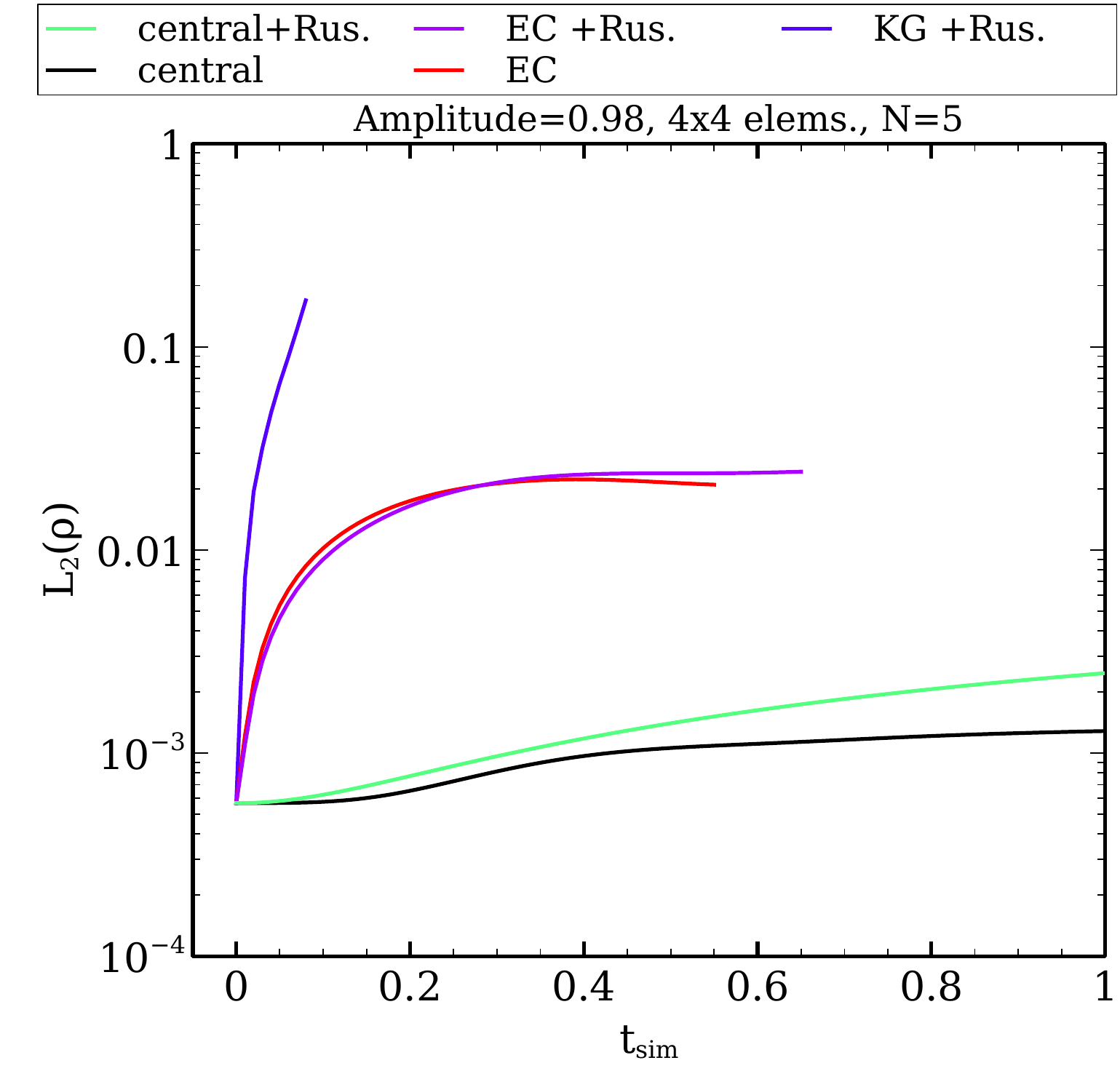} \includegraphics[width=\figwidth]{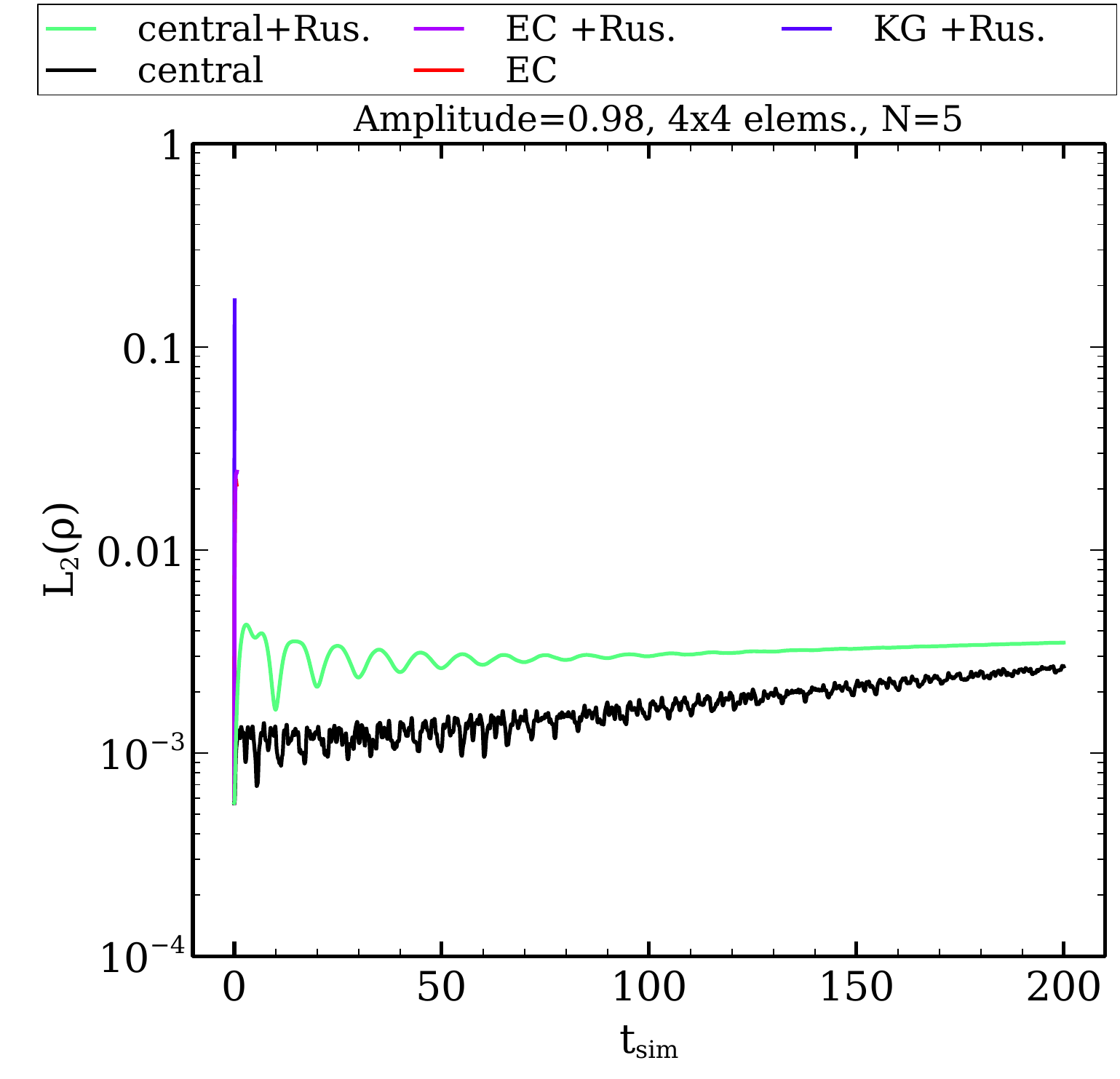} \\
\includegraphics[trim=0 0 0 40,clip,width=\figwidth]{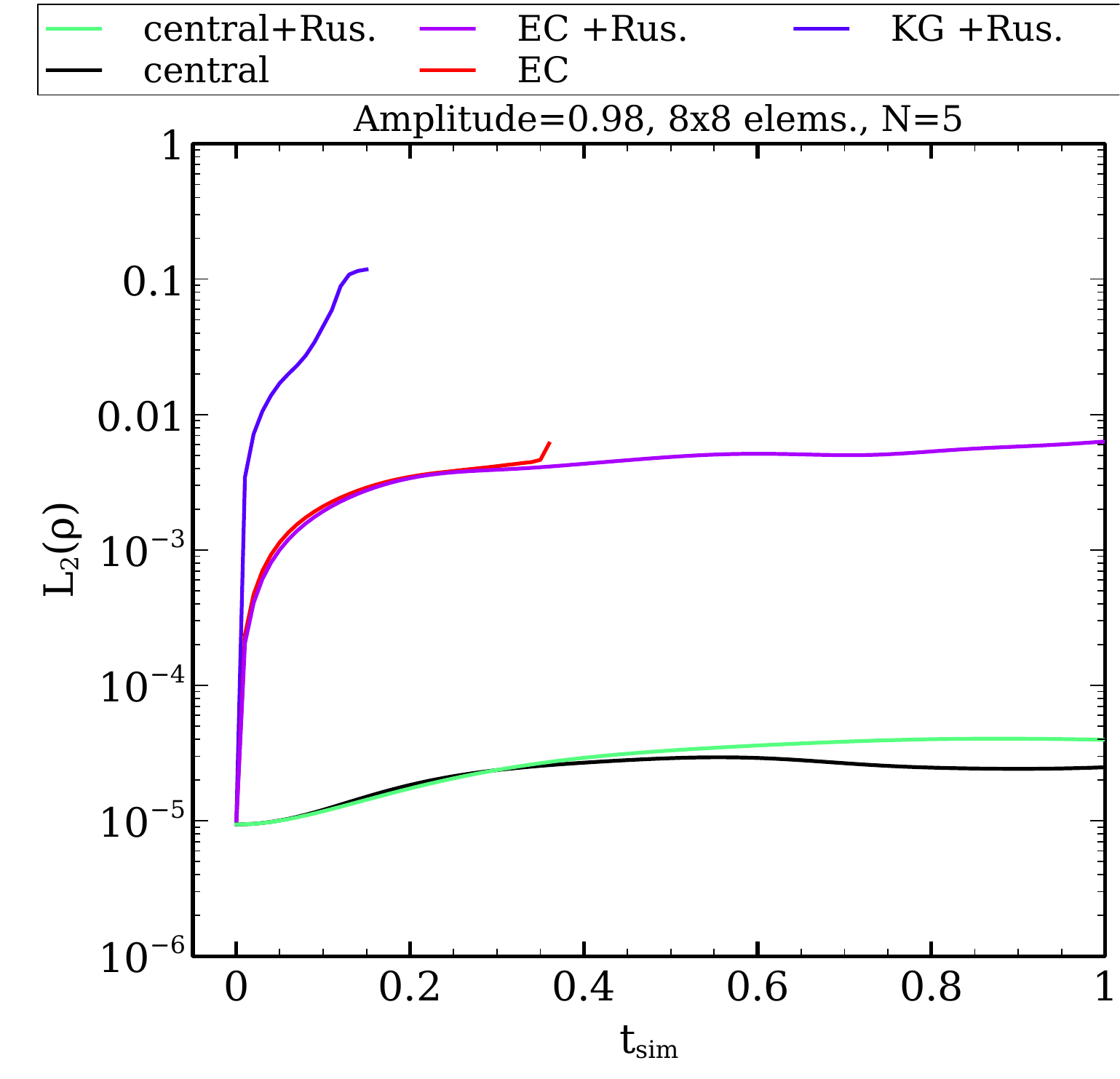}
\includegraphics[trim=0 0 0 40,clip,width=\figwidth]{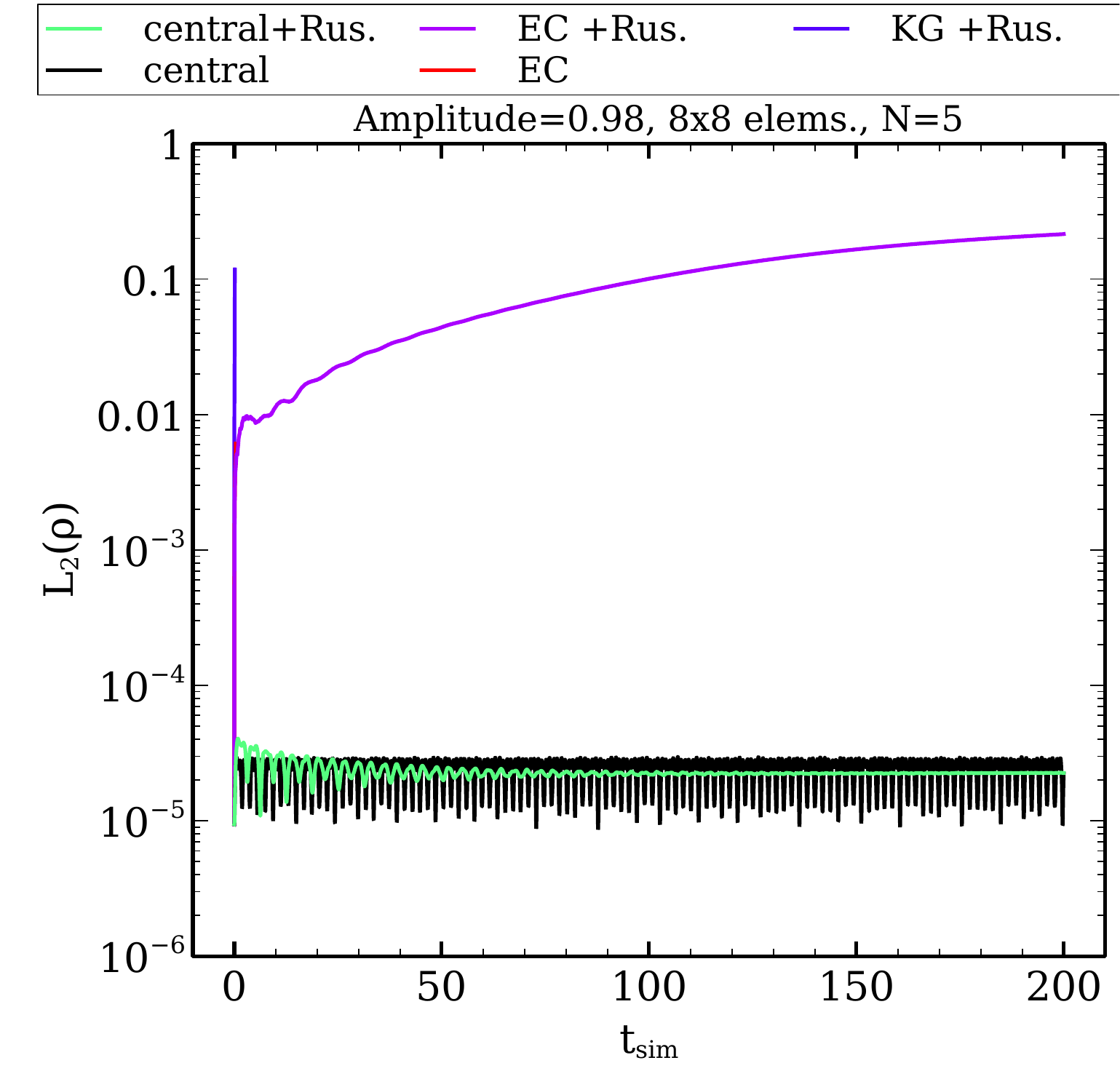}
\caption{\label{fig:euler_densitywave_L2_A098} $L_2$ error of density over simulation time, for initial density amplitude $A=0.98$ on $4\times4$ grid (top row) and $8\times8$ grid (bottom row). Left the initial phase $t_{sim}\leq1$ is shown, and $t_{sim}\leq200$ on the right.}
\end{figure}

\newpage
\begin{figure}[!htbp]
\centering
\includegraphics[width=\figwidth]{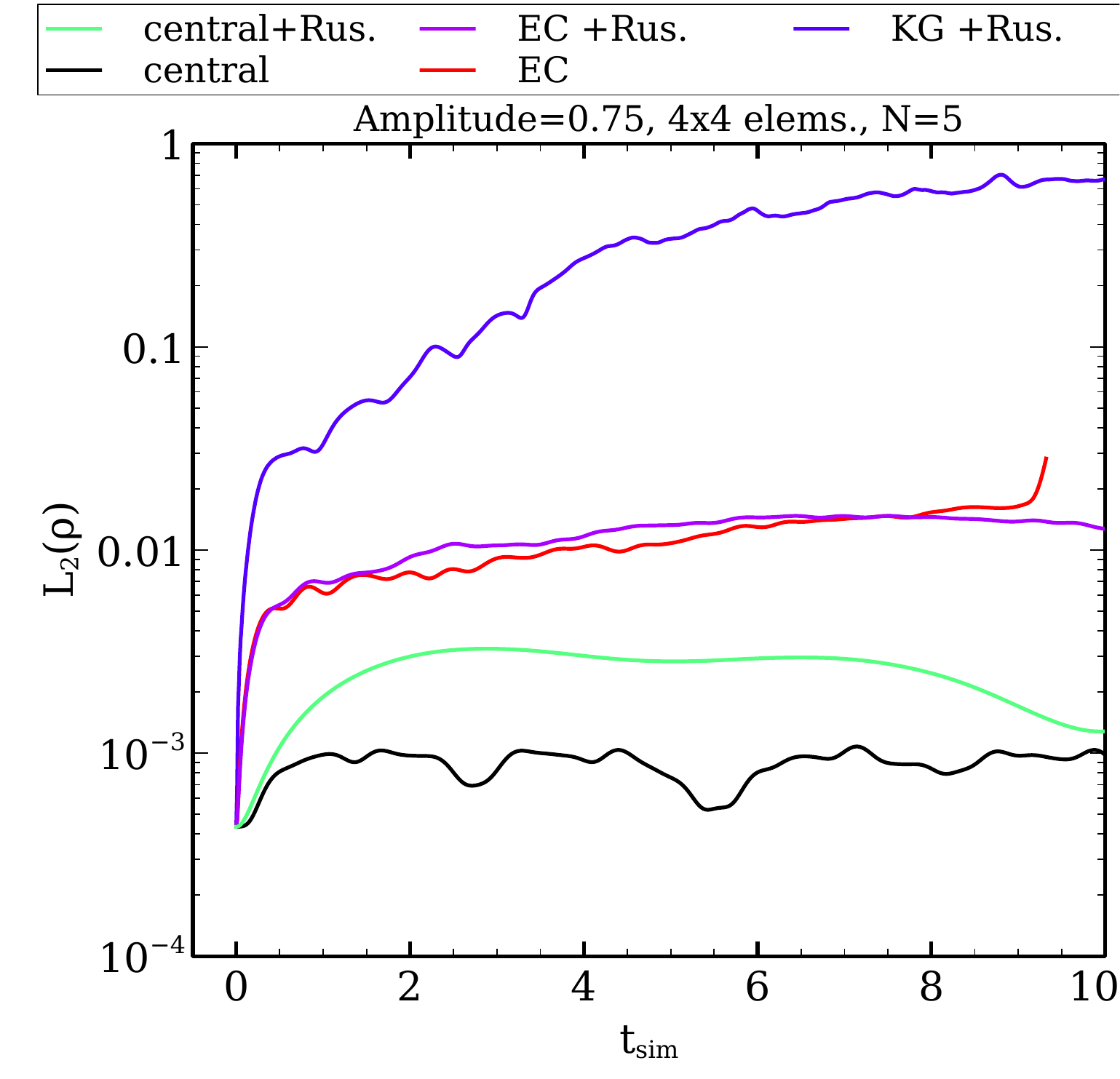} 
\includegraphics[width=\figwidth]{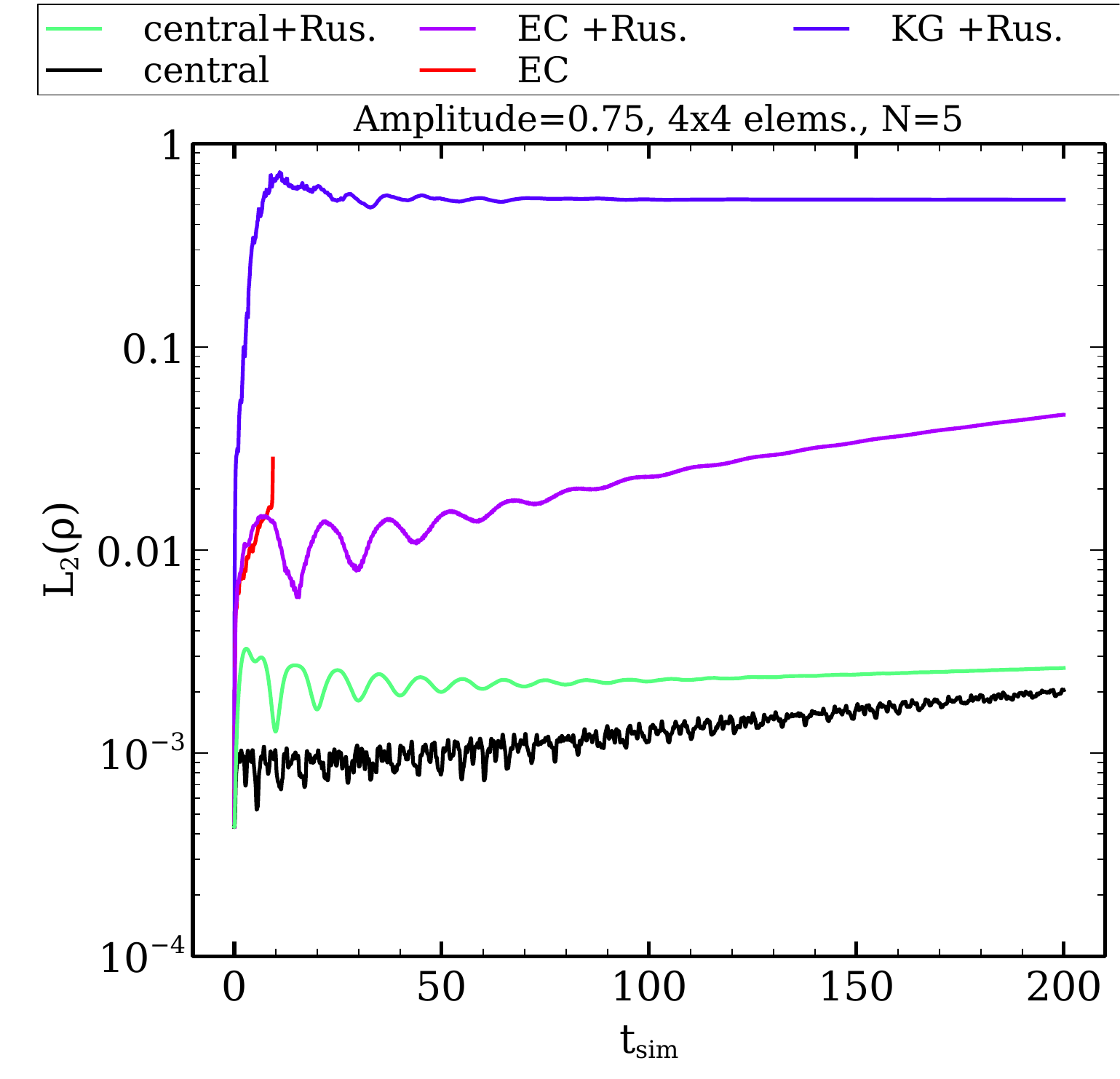} \\
\includegraphics[trim=0 0 0 40,clip,width=\figwidth]{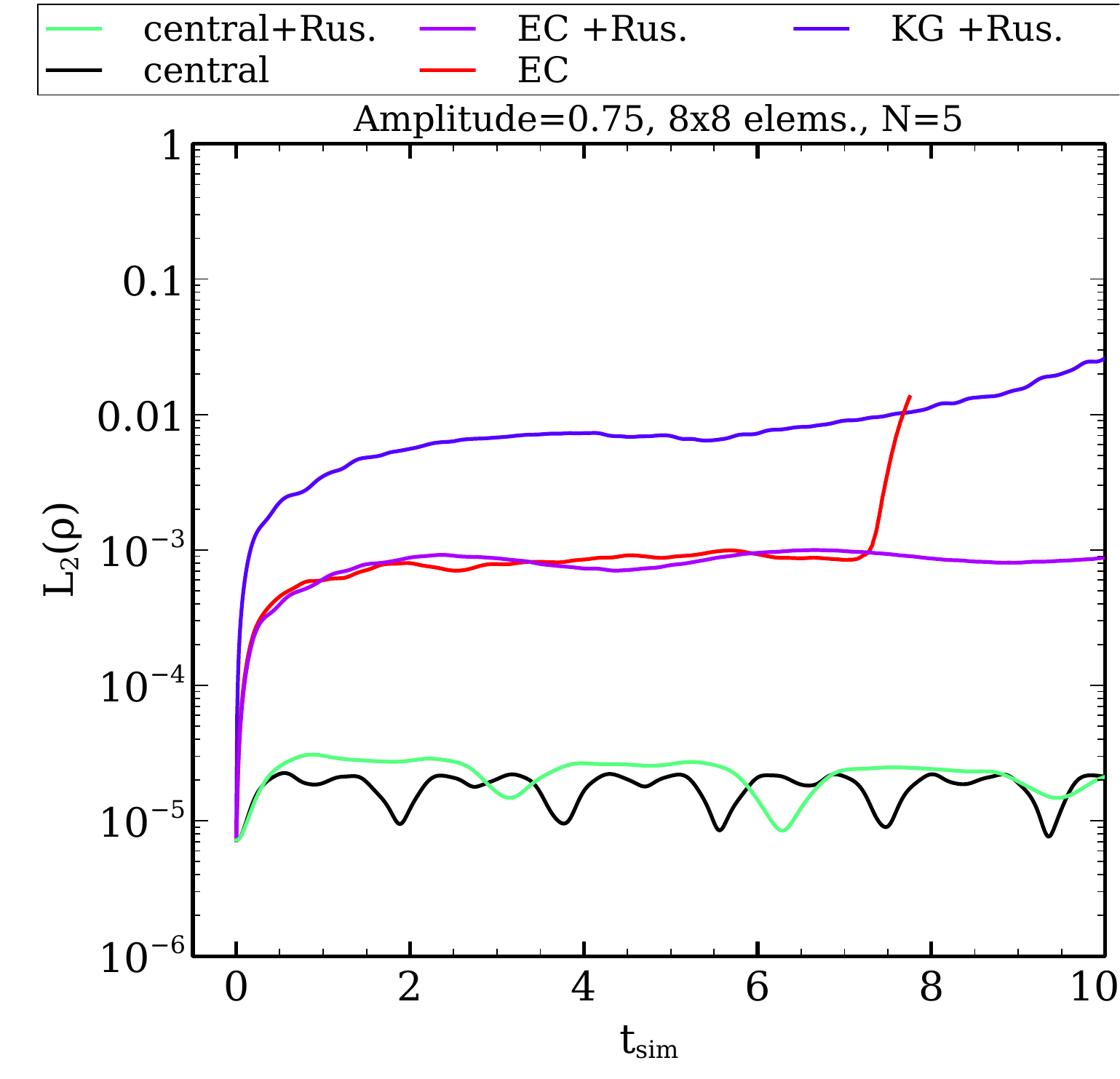}
\includegraphics[trim=0 0 0 40,clip,width=\figwidth]{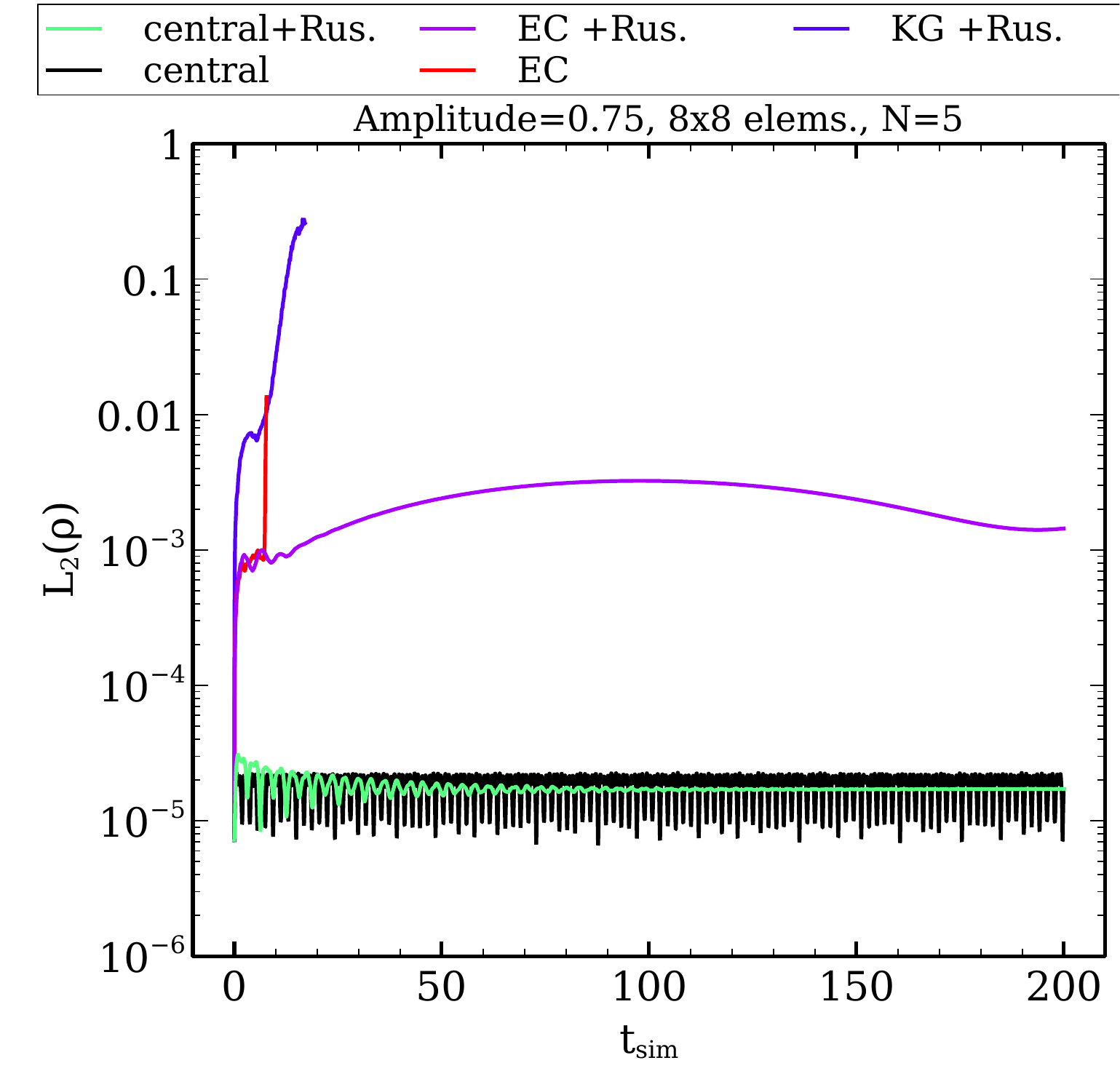}
\caption{\label{fig:euler_densitywave_L2_A075} $L_2$ error of density over simulation time, for initial density amplitude $A=0.75$ on $4\times4$ grid (top row) and $8\times8$ grid (bottom row). Left the initial phase $t_{sim}\leq10$ is shown, and $t_{sim}\leq200$ on the right.}
\end{figure}

\newpage
\begin{figure}[!htbp]
\centering
\includegraphics[width=\figwidth]{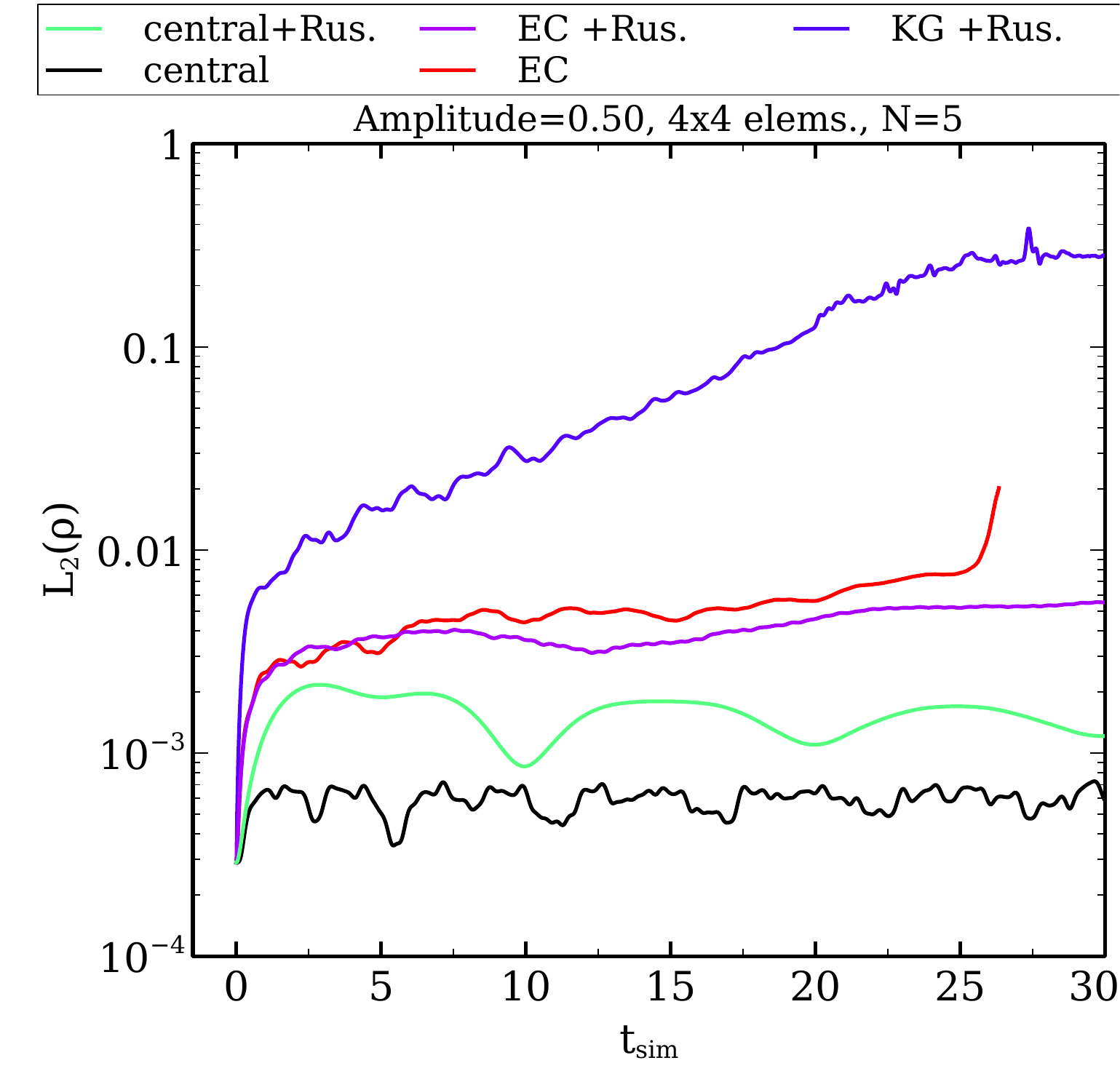} 
\includegraphics[width=\figwidth]{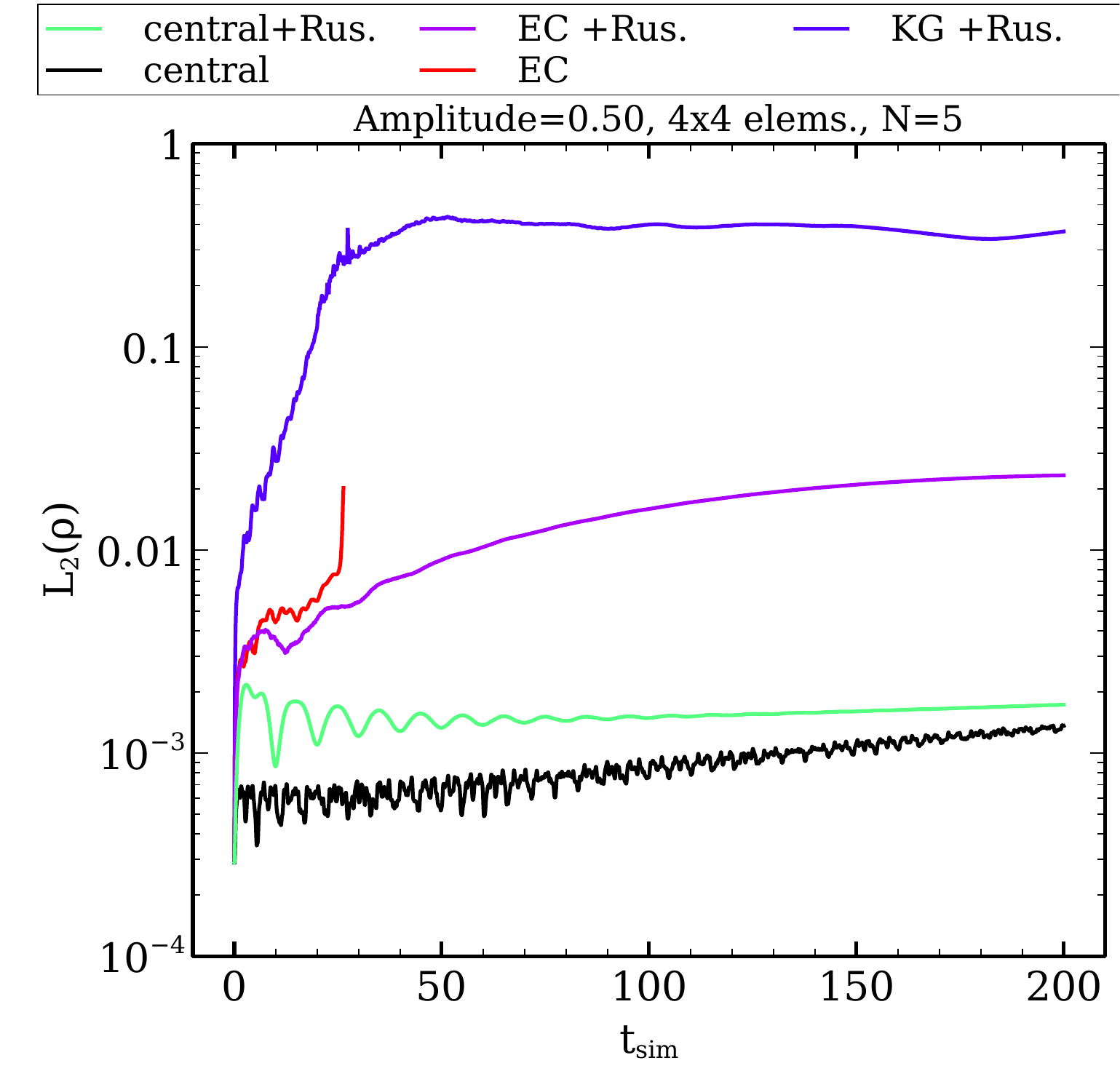} \\
\includegraphics[trim=0 0 0 40,clip,width=\figwidth]{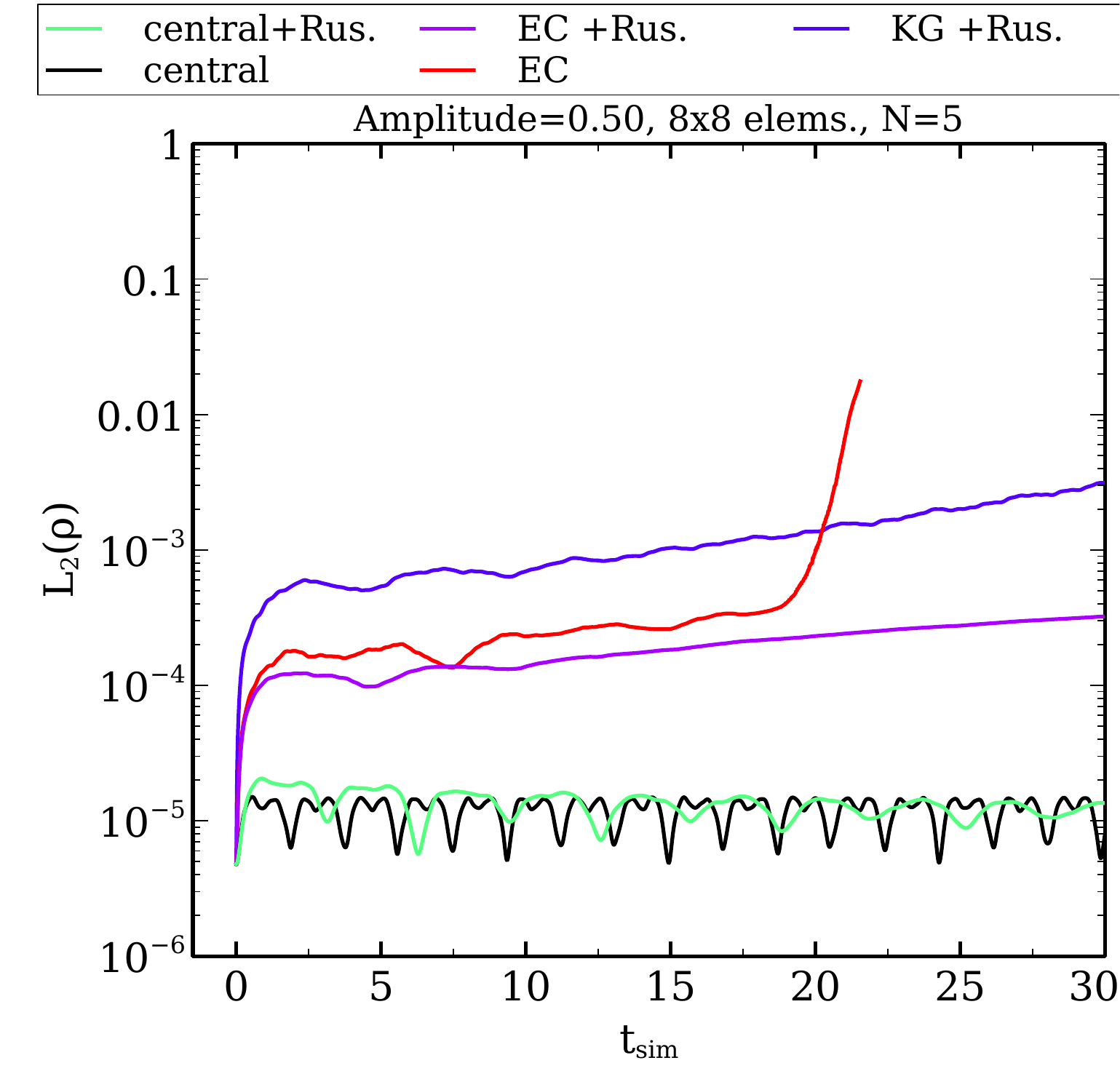}
\includegraphics[trim=0 0 0 40,clip,width=\figwidth]{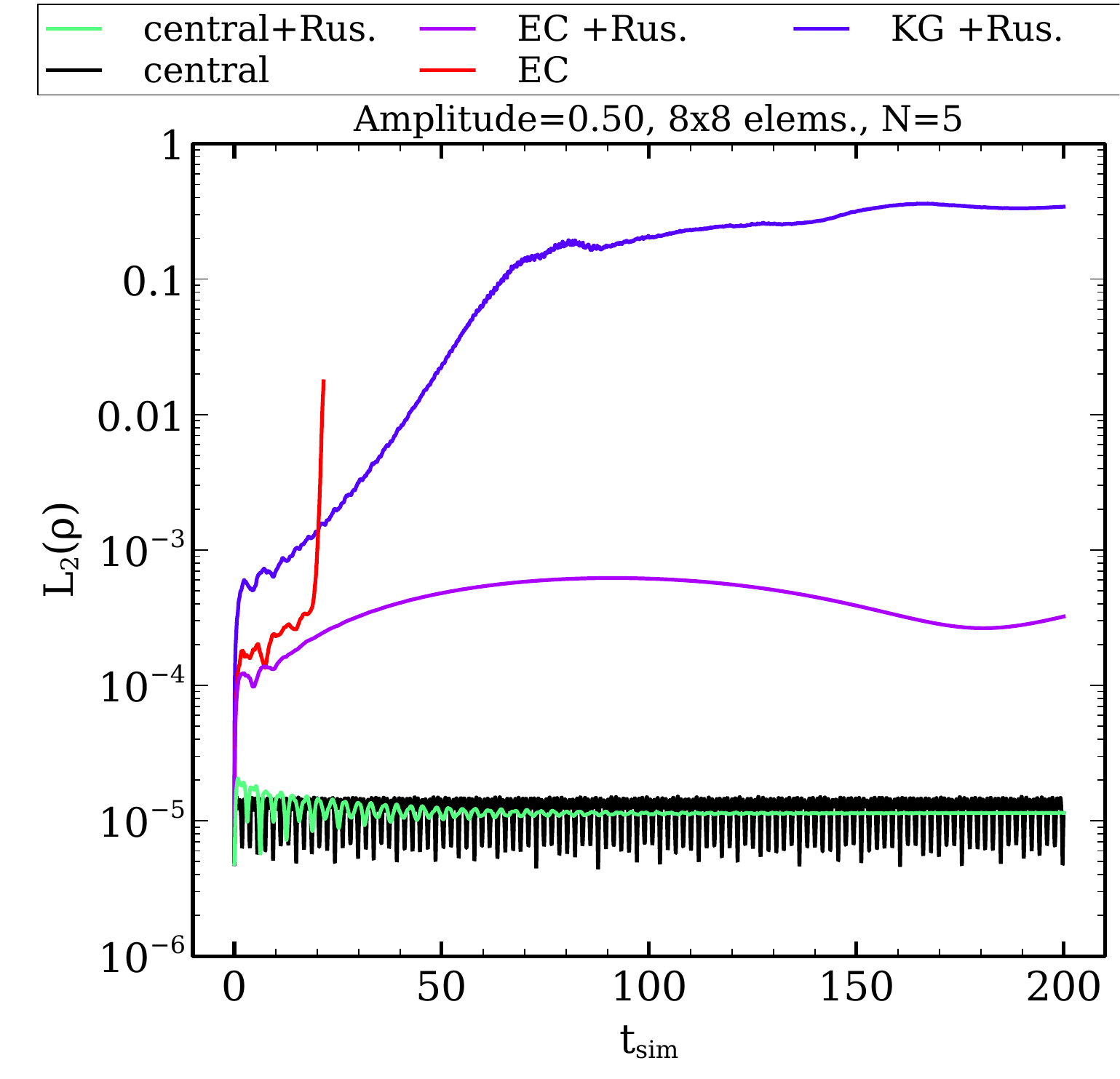}
\caption{\label{fig:euler_densitywave_L2_A050} $L_2$ error of density over simulation time, for initial density amplitude $A=0.5$ on $4\times4$ grid (top row) and $8\times8$ grid (bottom row). Left the initial phase $t_{sim}\leq30$ is shown, and $t_{sim}\leq200$ on the right.}
\end{figure}

\newcommand\figwidthird{0.19\textwidth}

\newpage
\begin{figure}[!htbp]
\centering
\begin{tabular}{|c|c|c|c|}\hline\rule{0pt}{3ex}
{\large $A$} & central+Rus. & EC+Rus. & KG+Rus. \\[1ex] \hline \rule{0pt}{0.14\textheight}
{\large $0.98$} &  
\includegraphics[width=\figwidthird]{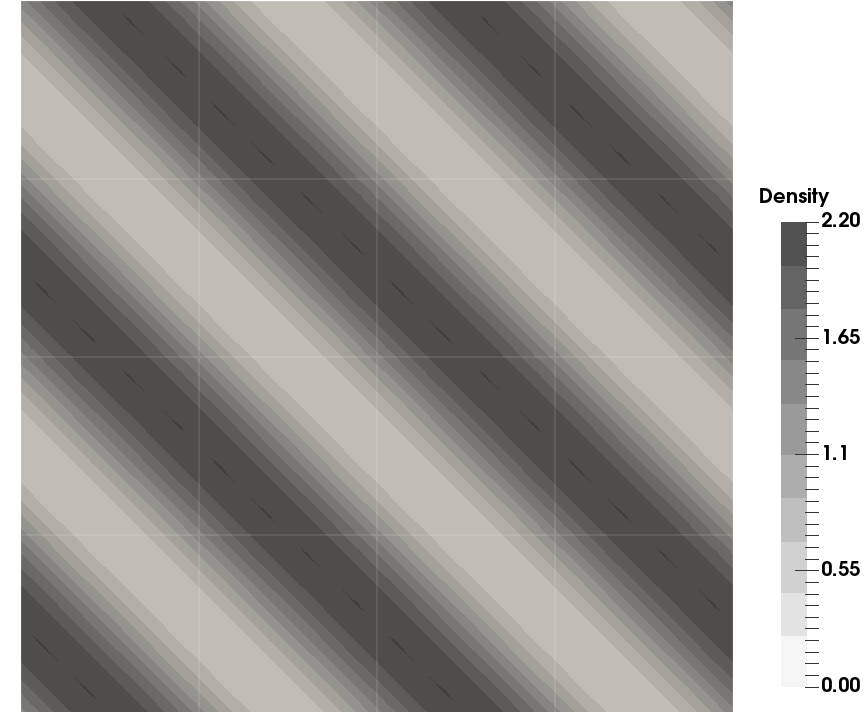} & 
\includegraphics[width=\figwidthird]{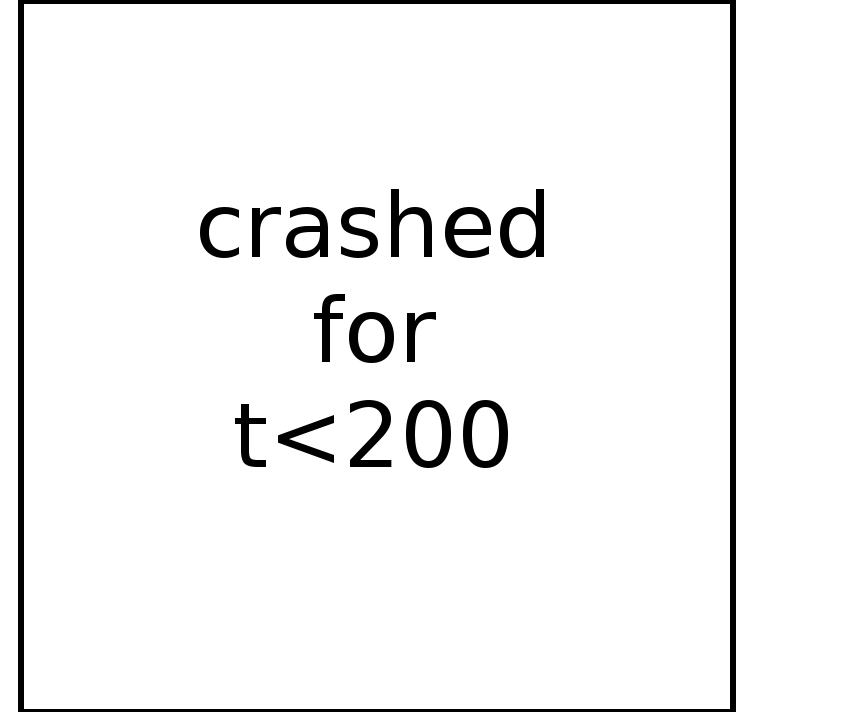} & \includegraphics[width=\figwidthird]{pics/density_negative_t_le_200.png} \\ &
\includegraphics[width=\figwidthird]{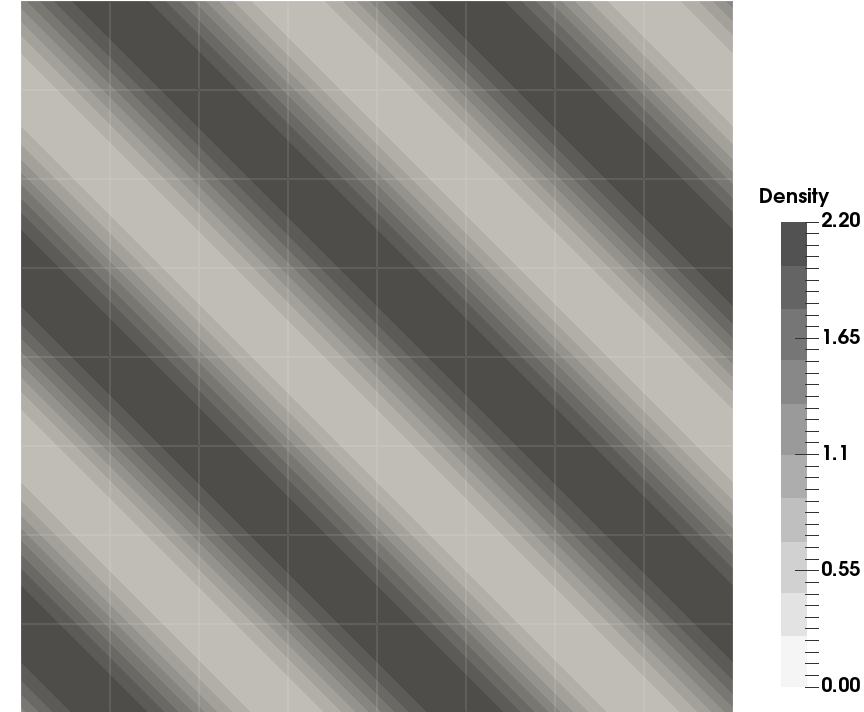} &
\includegraphics[width=\figwidthird]{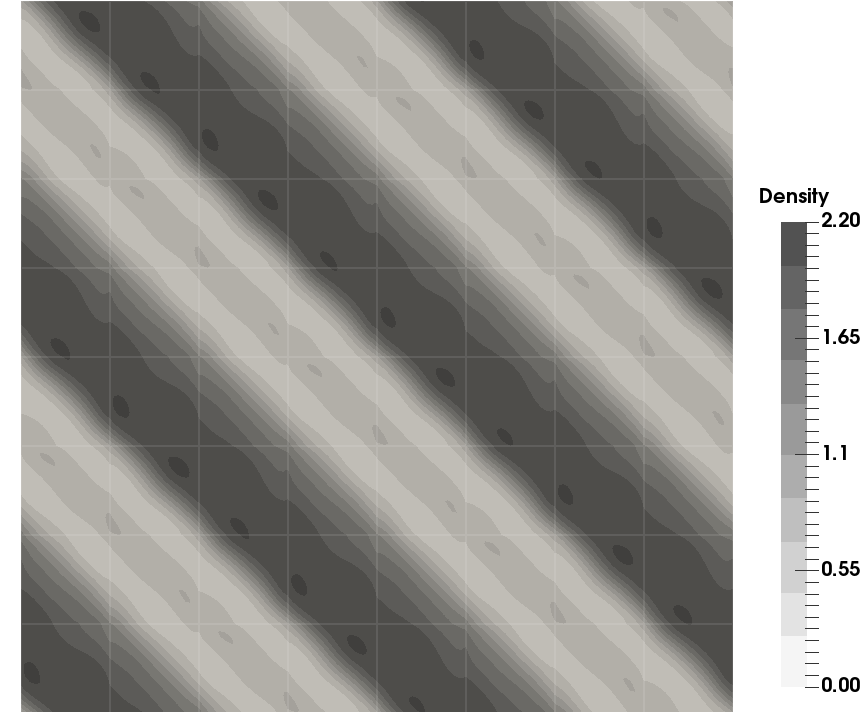} &
\includegraphics[width=\figwidthird]{pics/density_negative_t_le_200.png} \\[1ex]\hline\rule{0pt}{0.14\textheight}
 {\large $0.75$} &
\includegraphics[width=\figwidthird]{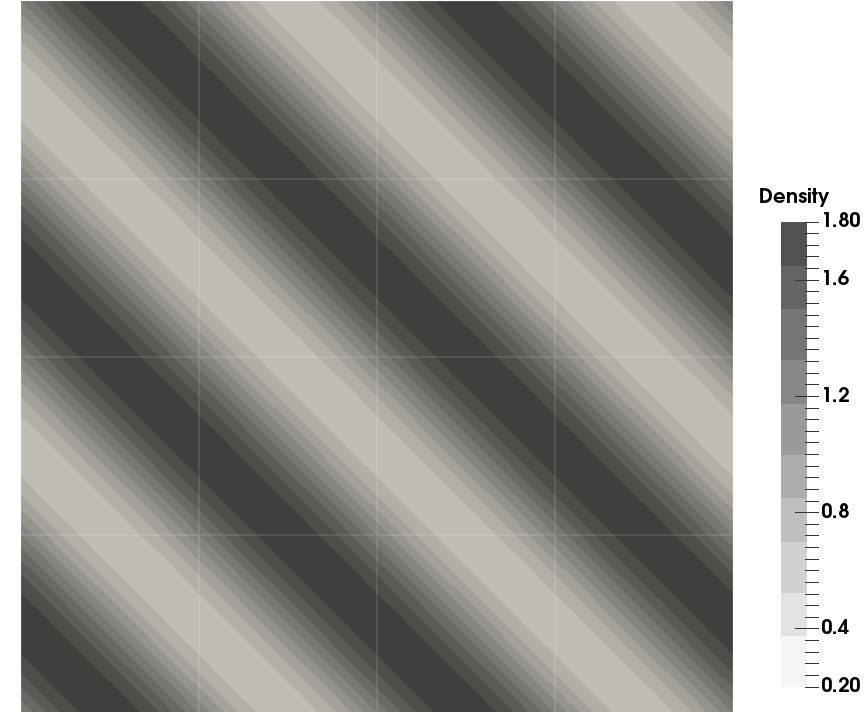} &
\includegraphics[width=\figwidthird]{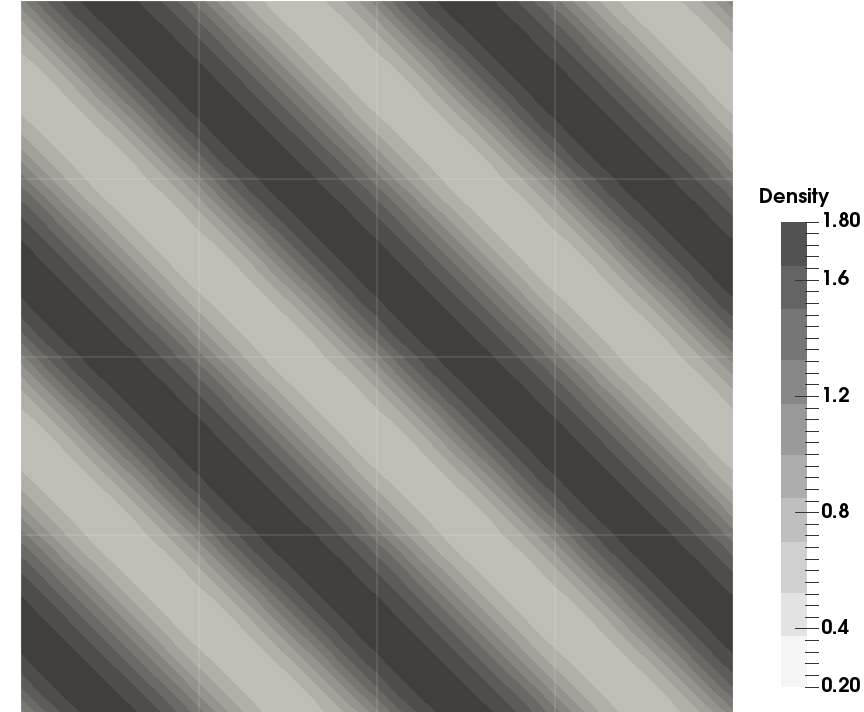} & 
\includegraphics[width=\figwidthird]{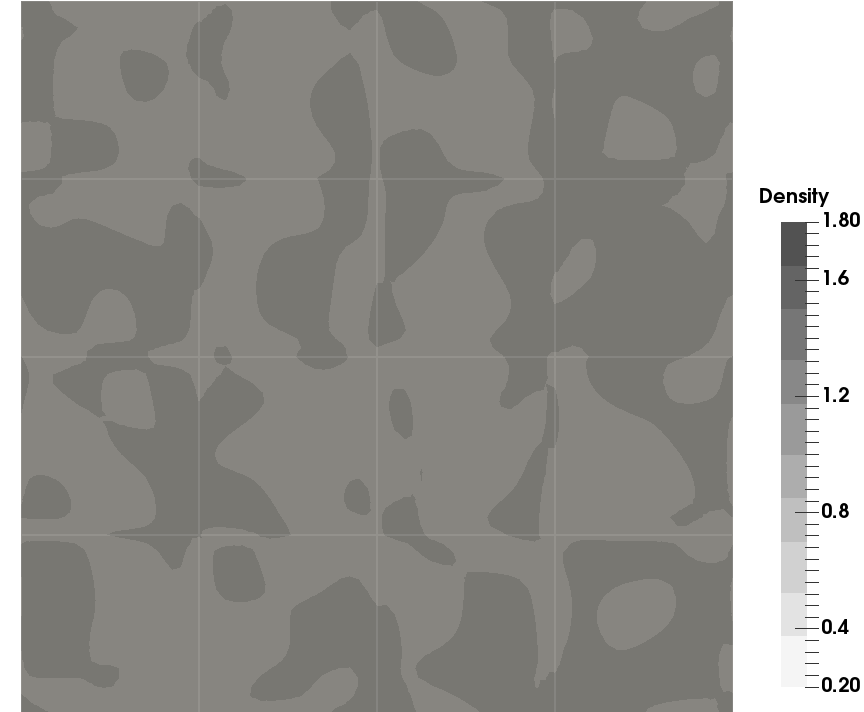} \\ &
\includegraphics[width=\figwidthird]{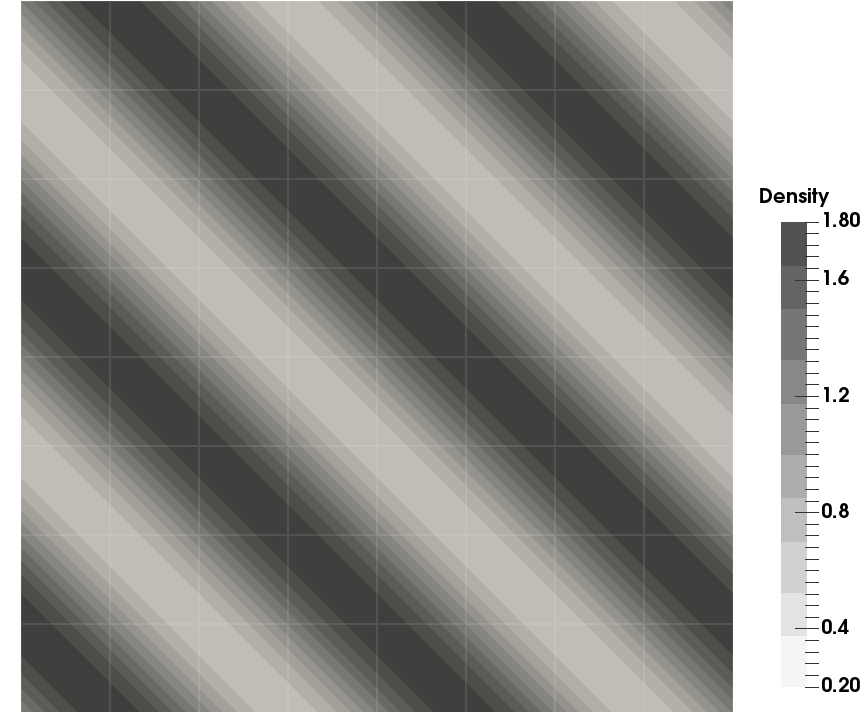} &
\includegraphics[width=\figwidthird]{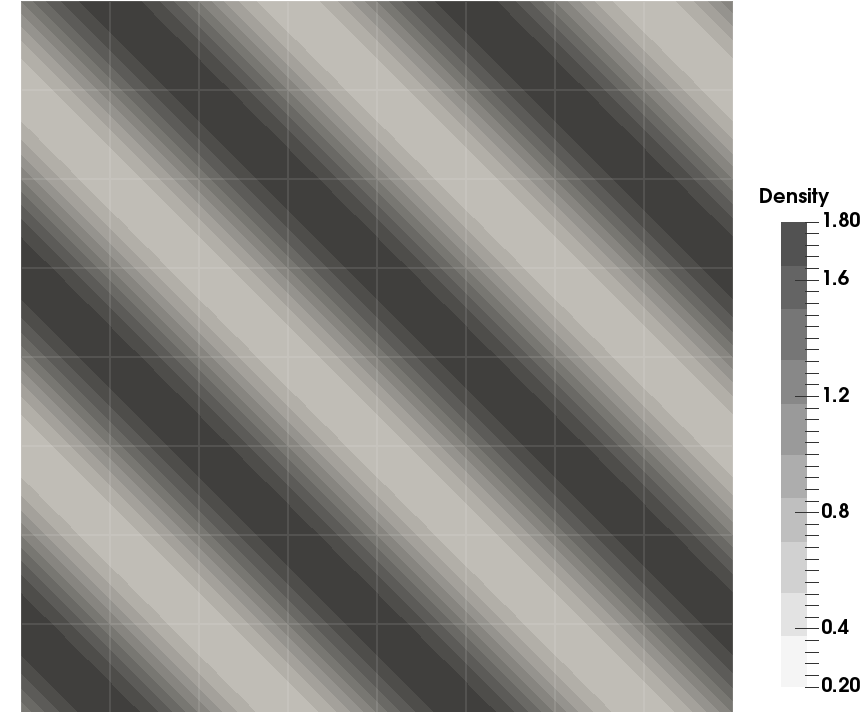} &  
\includegraphics[width=\figwidthird]{pics/density_negative_t_le_200.png} \\[1ex]\hline\rule{0pt}{0.14\textheight}
{\large $0.50$} & 
\includegraphics[width=\figwidthird]{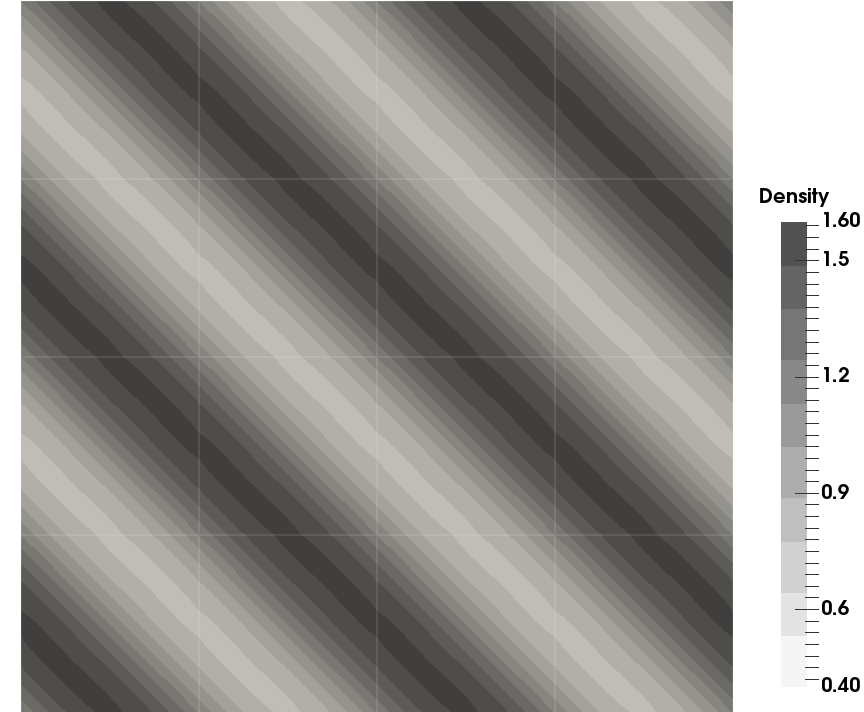} &
\includegraphics[width=\figwidthird]{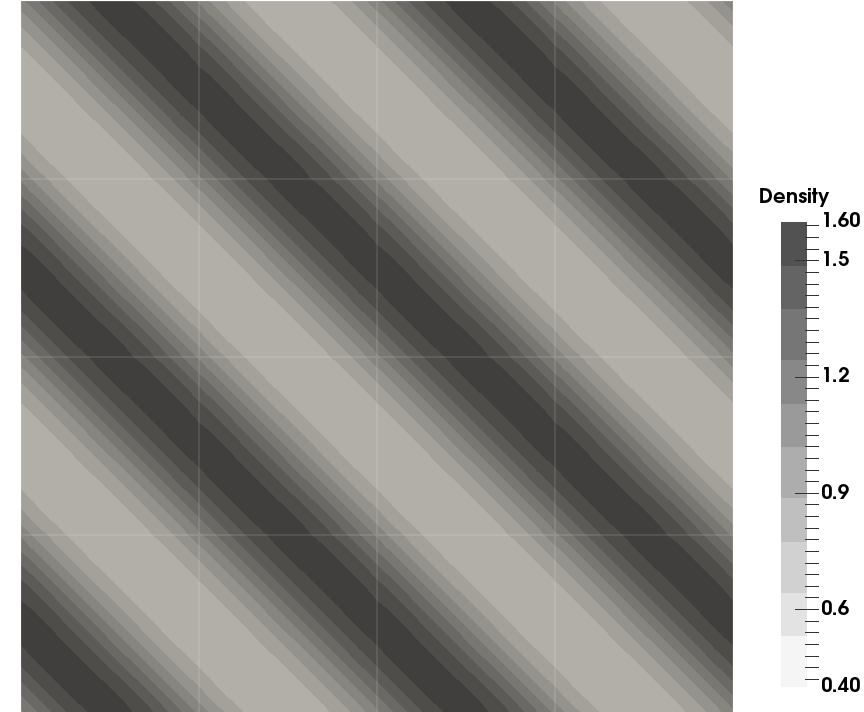} & \includegraphics[width=\figwidthird]{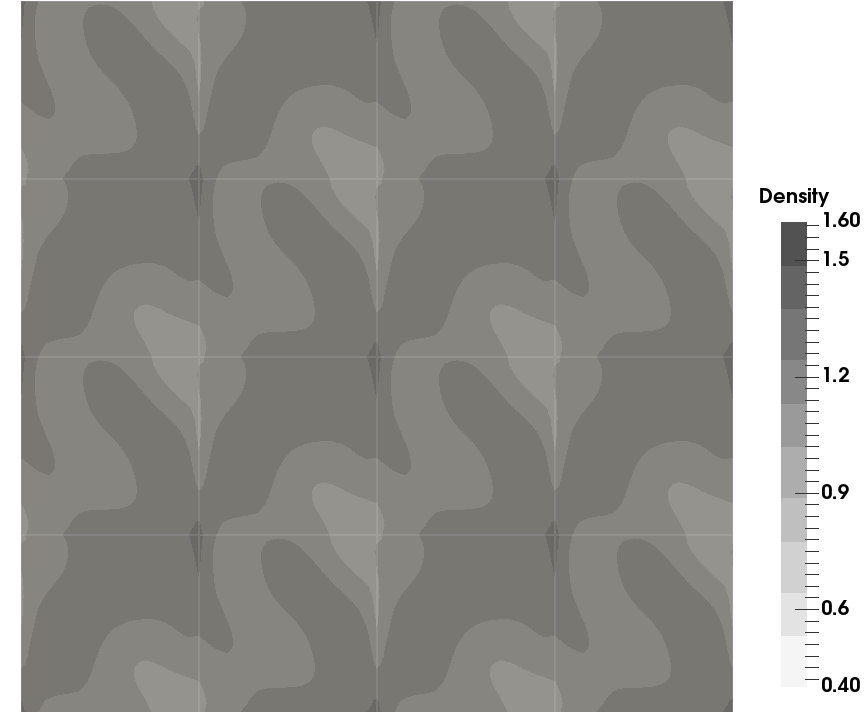} \\ &
\includegraphics[width=\figwidthird]{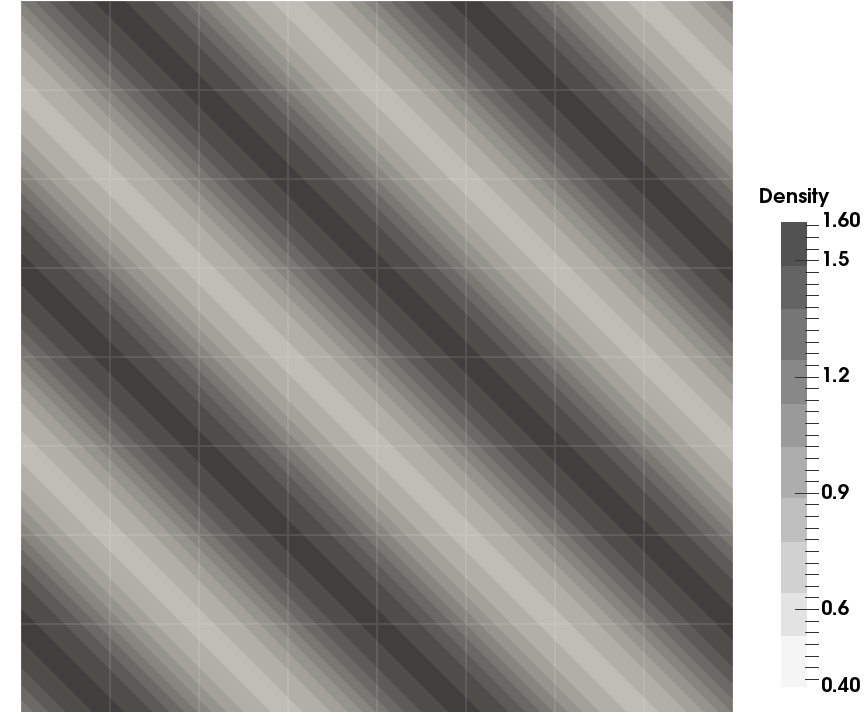} &
\includegraphics[width=\figwidthird]{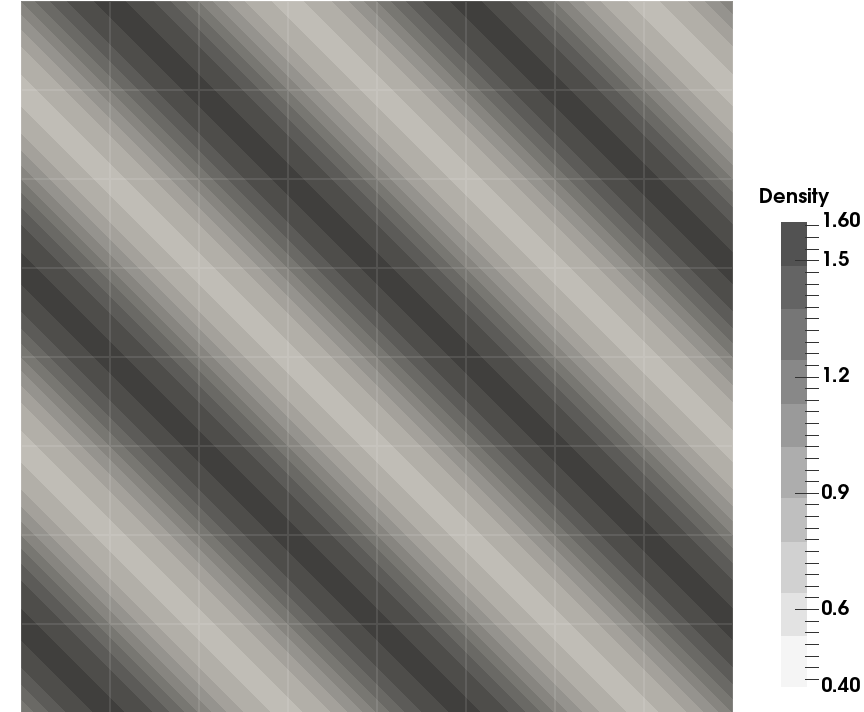} & \includegraphics[width=\figwidthird]{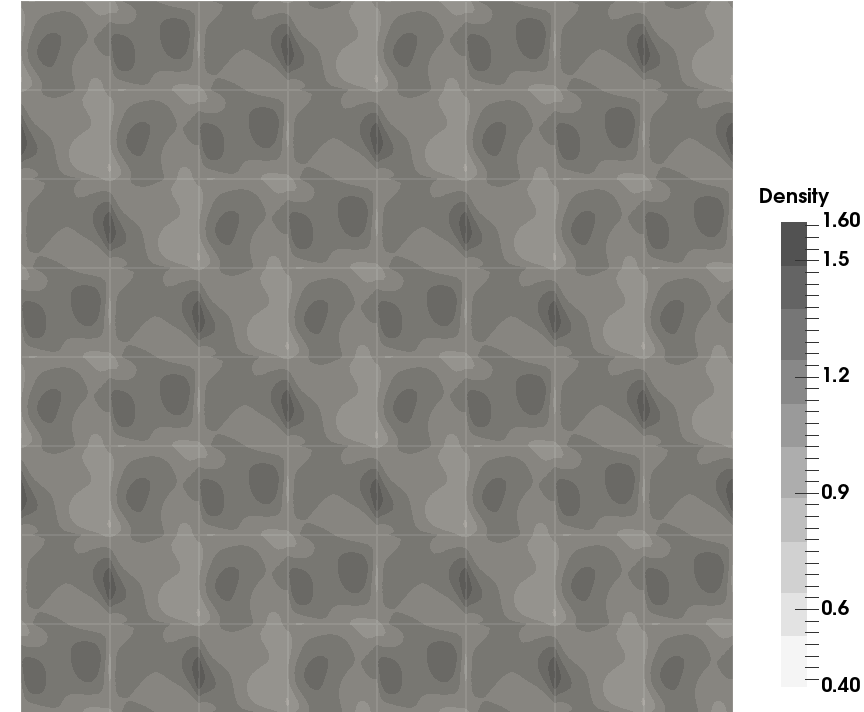} \\[1ex]\hline 
\end{tabular}
\caption{\label{fig:euler_density_t200} Final density distribution at $T=200$ for the initial density amplitudes $A=0.98/0.75/0.5$, with a Rusanov surface flux and the volume flux being either central (left), EC (middle) or KG (right), for a $4\times 4$ grid (upper row) and $8\times8$ grid (lower row). }
\end{figure}


\end{document}